\documentclass[twoside,leqno,onecolumn]{article}

\usepackage[letterpaper]{geometry}
\usepackage[T1]{fontenc}
\usepackage{amsfonts}
\usepackage{graphicx}
\usepackage{epstopdf}
\usepackage{enumitem}
\usepackage{algorithmic}
\usepackage{graphicx} 
\usepackage{tikz}
\usetikzlibrary{fit}
\usetikzlibrary {arrows.meta}
\usetikzlibrary{shapes.geometric}  
\usepackage[normalem]{ulem} 
\usepackage{float}
\usepackage{placeins}
\usepackage{flafter}

\makeatother

\usepackage{float}
\usepackage{mathtools}
\usepackage{amssymb}
\usepackage{booktabs}
\usepackage{multirow}
\usepackage{array}
\usepackage{longtable}
\usepackage{siunitx}
\usepackage{xurl}
\usepackage{tikz}
\usepackage{pgfplots}
\usepackage{enumitem}
\usepackage{amsthm}
\usepackage{framed}
\usepackage{xcolor}
\usepackage{hyperref}

\DeclarePairedDelimiter{\abs}{\lvert}{\rvert}
\DeclarePairedDelimiter{\norm}{\lVert}{\rVert}
\DeclareMathOperator{\diam}{diam}
\DeclareMathOperator{\conv}{conv}
\DeclareMathOperator{\MST}{MST}
\DeclareMathOperator{\med}{med}
\DeclareMathOperator{\sgm}{sgm}

\newcommand{\best}[1]{\ifmmode\bm{#1}\else\bfseries #1\fi}
\newcommand{\bestmath}[1]{{\boldmath$#1$}}

\newcommand{\Mone}{\textup{M1}}
\newcommand{\Mtwo}{\textup{M2}}

\newcommand{\ExactBCLT}{\textup{BCLT}}
\newcommand{\Identity}{\hyperlink{Identitytarget}{\mbox{M2}}}
\newcommand{\Idenone}{\hyperlink{Idenonetarget}{\mbox{M1}}}

\pgfplotsset{compat=1.18}
\usepackage{tikz}
\usepackage{pgfplots}
\pgfplotsset { compat=1.18 }
\usetikzlibrary
  {
    shapes.geometric,
    positioning,
    matrix,
    calc,
  }
\usepgfplotslibrary{groupplots}

\definecolor{btterminal}{RGB}{29,78,216}
\definecolor{btsteiner}{RGB}{217,119,6}
\definecolor{btedge}{RGB}{2,132,199}
\definecolor{btbottleneck}{RGB}{185,28,28}
\definecolor{btbox}{RGB}{156,163,175}
\definecolor{btwirecol}{RGB}{209,213,219}
\definecolor{btframe}{RGB}{107,114,128}

\tikzset
  {
    terminal/.style =
      {
        circle,
        fill = btterminal,
        draw = btterminal,
        inner sep = 0pt,
        minimum size = 3.4pt,
      },
    steiner/.style =
      {
        regular polygon,
        regular polygon sides = 3,
        fill = btsteiner,
        draw = btsteiner,
        inner sep = 0pt,
        minimum size = 5pt,
      },
    btedge/.style =
      {
        draw = btedge,
        line width = 0.9pt,
      },
    btbottleneck/.style =
      {
        draw = btbottleneck,
        line width = 1.6pt,
      },
    btedge3d/.style =
      {
        draw = btedge,
        line width = 0.9pt,
        preaction = { draw = white, line width = 2.6pt },
      },
    btbottleneck3d/.style =
      {
        draw = btbottleneck,
        line width = 1.6pt,
        preaction = { draw = white, line width = 3.4pt },
      },
    btangle/.style =
      {
        draw = btbottleneck,
        line width = 1.1pt,
      },
    bthull/.style =
      {
        draw = btbox,
        dashed,
        line width = 0.5pt,
      },
    btwire/.style =
      {
        draw = btwirecol,
        line width = 0.4pt,
      },
    btbox/.style =
      {
        draw = btbox,
        dashed,
        line width = 0.4pt,
      },
    btdisc/.style =
      {
        draw = btbox,
        line width = 0.4pt,
        opacity = 0.30,
      },
    nodelabel/.style =
      {
        font = \small,
        btterminal,
      },
    elabel/.style =
      {
        font = \footnotesize,
        btbox,
        fill = white,
        inner sep = 0.6pt,
      },
    legendbox/.style =
      {
        draw = gray!60,
        fill = white,
        rounded corners = 1pt,
        inner sep = 2.5pt,
        font = \small,
      },
  }

\pgfplotsset
  {
    steineraxis/.style =
      {
        width = 0.82\linewidth,
        height = 0.45\linewidth,
        grid = major,
        grid style = { line width = .1pt, draw = gray!20 },
        major grid style = { line width = .2pt, draw = gray!50 },
        legend style = { font = \small, draw = gray!60, fill = white },
        tick label style = {font = \small},
        label style = { font = \small },
        every axis plot/.append style = { thick },
      },
    steinergeo/.style =
      {
        scale only axis,
        xtick = \empty,
        ytick = \empty,
        enlargelimits = false,
        axis line style = { thin, btframe },
      },
    steineraxis3d/.style =
      {
        view = {120}{22},
        scale only axis,
        unit vector ratio = 1 1 1,
        xtick = \empty,
        ytick = \empty,
        ztick = \empty,
        enlargelimits = false,
        axis line style = { thin, btframe },
        legend style = { font = \small, draw = gray!60, fill = white },
      },
  }

\theoremstyle{plain}
\newtheorem{theorem}{Theorem}[section]
\newtheorem{proposition}[theorem]{Proposition}
\newtheorem{property}[theorem]{Property}
\newtheorem{lemma}[theorem]{Lemma}
\newtheorem{corollary}[theorem]{Corollary}
\theoremstyle{definition}

\newtheorem{model}{Model}

\colorlet{shadecolor}{black!10}
\AddToHook{env/model/before}{\shaded}
\AddToHook{env/model/after}{\endshaded}

\usepackage[authoryear,longnamesfirst]{natbib}
\setcitestyle{authoryear,open={[},close={]},aysep={,}}

\usepackage{xurl}
\AtEndPreamble{\hypersetup{breaklinks=true}}

\newtheorem{remark}[theorem]{Remark}

\definecolor{purple}{RGB}{120,40,170}

\begin{document}

\thispagestyle{empty}

\newcommand\relatedversion{}

\title{\Large Mixed-Integer Second-Order Cone Formulations for the \\Euclidean Bottleneck Steiner Tree Problem in $\mathbb{R}^d$}
    \author{Gabriel P. Braun\thanks{Universidade Federal do Rio de Janeiro, braun@cos.ufrj.br
    }
    \and Marcia Fampa\thanks{Universidade Federal do Rio de Janeiro,
  fampa@cos.ufrj.br
  }
  \and Nelson Maculan\thanks{Universidade Federal do Rio de Janeiro,
  maculan@cos.ufrj.br
  }
  \and Renan V. Pinto\thanks{Universidade Federal do Rio de Janeiro, renanvp@im.ufrj.br
  }
  }

\date{\today}

\maketitle

\pagenumbering{arabic}

\begin{abstract}
  Given a set of terminals in $\mathbb{R}^d$ and an integer $k\ge0$, the Euclidean bottleneck Steiner tree problem (EBSTP) asks for a tree spanning the terminals with at most $k$ Steiner points that minimizes its longest edge. Although the planar case has been widely studied, exact optimization methods for arbitrary dimension remain scarce. We develop two exact mixed-integer second-order cone programming (MISOCP) formulations for the EBSTP: a single-commodity flow formulation and a formulation with connectivity cuts. Both formulations share a geometric core based on edge-conditional second-order cone constraints and are strengthened by problem-specific preprocessing, edge-dependent big-$M$ constants, symmetry-breaking constraints, constructive upper bounds, and new global lower bounds on the optimal value. We establish the correctness of the formulations and prove the validity of all strengthening techniques. Computational experiments on 1050 benchmark instances in dimensions $d\in\{2,3,4\}$, with up to $50$ terminals and $9$ Steiner points, show that the single-commodity flow formulation consistently outperforms the formulation with connectivity cuts, with an advantage { in solved instances} that grows with dimension. The proposed global lower bounds substantially strengthen the root relaxation, while the single-commodity flow formulation provides the most robust exact solution method across all test sets. These results establish the single-commodity flow formulation as an effective exact approach for the EBSTP in arbitrary fixed dimension.
\end{abstract}

\section{Introduction}

The Euclidean Steiner tree problem in $\mathbb{R}^d$ is to find a network of minimum length interconnecting a set of given terminals, optionally adding auxiliary points (\emph{Steiner points}) to reduce the total length~\citep*{hwang1992,smith1992}. The classical objective is additive: each edge contributes to the total cost. The bottleneck variant studied in this paper replaces the additive objective by a min-max one: we minimize the length of the longest edge. The bottleneck criterion arises in geometric design problems where the limiting factor is the maximum admissible span of an individual link rather than the total length of the network. An example is very-large-scale integration routing, where the bottleneck value captures the worst-case signal delay across a routed net~\citep*{sarrafzadeh1992}; the same criterion arises in settings such as wireless relay placement, where it is the maximum transmission distance between adjacent links, or geometric infrastructure design under a worst-case link reliability constraint.

We study the following problem: given $p$ terminals $a_1,\dots,a_p\in\mathbb{R}^d$ and an integer $k\ge 0$, find a tree spanning all terminals using at most $k$ Steiner points and minimizing the longest edge. Bottleneck Steiner trees in $\mathbb{R}^d$ are ill-posed as $k\to\infty$, since any edge can be subdivided arbitrarily finely, causing the bottleneck value to approach zero. The budget $k$ gives the problem a finite-dimensional search space. We refer to this problem as the \emph{Euclidean bottleneck Steiner tree problem (EBSTP)}. \citet*{bae2009exact} use the same budget restriction in the plane. For $k=0$, the EBSTP reduces to the (polynomially solvable) bottleneck minimum spanning tree problem on the terminals; imposing degree bounds on that spanning tree already makes it \textsf{NP}-hard in the plane~\citep*{andersen2016minimum}.

{ The literature on bottleneck Euclidean Steiner trees is concentrated in the plane. \citet*{sarrafzadeh1992} and \citet*{ganley1996} gave algorithms for a fixed topology, \citet*{wang2002} proved the planar problem \textsf{NP}-hard and gave approximation algorithms, as did \citet*{wang2002bsk} for the $k$-Steiner variant, and \citet*{bae2009exact,bae2011exact} (henceforth \ExactBCLT{}) gave the first exact algorithm for fixed $k$ in any $L_q$ metric in the plane, with running time $O(f(k)\cdot(p^k + p\log p))$, which we use as a baseline in Section~\ref{sec:planarexp}. More recently, \citet{bandyapadhyay2024} showed that the planar problem is fixed-parameter tractable in $k$, \citet*{li2024line} solved a line-constrained variant, and \citet*{brazil2012bottleneck,brazil2015exact,brazil2013generalised} studied the bottleneck $2$-connected $k$-Steiner network problem in normed planes. These exact methods are combinatorial and planar, and none extends directly to $\mathbb{R}^d$ with $d\ge 3$.}

{ By contrast, mixed-integer nonlinear programming (MINLP) formulations for the additive Euclidean Steiner tree problem cover general dimension. \citet*{maculan2000} introduced a nonconvex MINLP formulation (MMX), and \citet*{fampa2001,fampa2004} proposed a convex MINLP formulation (FM) with edge-length variables and a big-$M$ term that deactivates the length constraint of unselected edges, recast in conic form; throughout this paper, a MINLP formulation is termed \emph{convex} if relaxing its integrality constraints yields a convex optimization problem. These formulations were progressively strengthened by \citet*{fampa2008,vanlaarhoven2013,dambrosio2015,fampa2016branch,ouzia2022,pinto2023}. The surveys of \citet*{fampa2016}, \citet*{fampa2019} and \citet*{maculan2025} review this line of work. \citet*{blanco2025fixed} apply similar formulations to fixed-topology trees with neighborhoods in $\mathbb{R}^d$.}

We adopt the modeling structure of FM and adapt it to the bottleneck objective, which is a min-max rather than a sum{ ; Section~\ref{sec:relation} discusses what changes in the formulation and in its relaxation as a consequence}. To the best of our knowledge, no exact mixed-integer second-order cone programming (MISOCP) { formulation} for the EBSTP in $\mathbb{R}^d$ has been studied previously.

The main contributions of this work are the following.

\begin{itemize}[leftmargin=*,itemsep=0.3em]
  \item We develop the first exact mixed-integer second-order cone programming
    (MISOCP) formulations for the EBSTP in arbitrary fixed dimension{ : a
    single-commodity flow formulation (\Mone{}) and a formulation with
    dynamically separated connectivity cuts (\Mtwo{}). Both share a geometric core in
    which a single epigraph variable is linked to the selected edges through
    edge-conditional second-order cone constraints. The modeling elements are standard;
    the results specific to the bottleneck objective are the structural properties that
    justify them and the analysis of their relaxation
    (Sections~\ref{sec:problem}, \ref{sec:formulations} and~\ref{sec:relation})}.

  \item { We derive strengthening techniques that apply to both formulations:
    a constructive upper bound, used for a partial mixed-integer programming (MIP) start and for edge elimination; edge-type-specific big-$M$ constants; the restriction of Steiner coordinates to the terminal box; and a symmetry-breaking ordering of the Steiner labels.}
    {  We also show that the terminal--terminal candidate edges can be restricted to a minimum
spanning tree of the terminals in every dimension. 
\citet*{bae2010exact,bae2011exact} prove
this restriction in the plane, and 
their exchange argument extends it to $\mathbb{R}^d$.
 The same exchange, run in the opposite direction, shows that the edges of this tree shorter than a lower bound on the optimal value can be fixed to one.}
    { We also identify the main challenges these techniques leave unresolved (Section~\ref{sec:strengthenings}).}

{ \item We derive two global lower bounds on the optimal value, including a new
clustering bound $L^\star_{\le k}\ge\tfrac12 e_k$ for $k\ge1$, where $e_1\ge e_2\ge\dots\ge e_{p-1}$ are the
edge lengths of a minimum spanning tree of the terminals sorted in nonincreasing order,  and whose constant $\tfrac12$ is best possible (Section~\ref{sec:diam}).}

  { \item We analyze the continuous relaxation of the flow formulation. We show that it always has an optimal solution in which all Steiner points coincide, characterize exactly the instances on which its bound is zero, which never happens for $k=1$, and prove that a valid inequality $L\ge\ell$ raises the relaxation value to $\max\{\ell,v_k\}$, where $v_k$ is the value without it. Hence the global lower bounds are the only strengthening techniques that raise the root bound on such instances, and on them the root bound equals the larger of the two bounds (Section~\ref{sec:relaxation}).}

  { \item We present a computational study in dimensions $d\in\{2,3,4\}$ and on structured instances in dimensions up to $6$, including a comparison with \ExactBCLT{} in the plane. The experiments compare the two formulations proposed, quantify the effect of each strengthening technique, and show that \Mone{} is the more effective and robust exact method, with an advantage in solved instances that grows with the dimension and the Steiner budget (Section~\ref{sec:computational}).}
\end{itemize}

The remainder of the paper is organized as follows. Section~\ref{sec:problem} formalizes the problem, lists basic properties of optimal solutions and develops the indexed representation used throughout. Section~\ref{sec:contrast} contrasts bottleneck Steiner trees with classic Euclidean Steiner minimal trees and shows that the classic structural regularities need not survive the change of objective. Section~\ref{sec:formulations} presents the formulations \Mone{}
and \Mtwo{} {   
 and relates them to the additive formulations MMX and FM}. Section~\ref{sec:strengthenings} develops the structural strengthening techniques common to
both formulations. Section~\ref{sec:computational} reports the computational study. Section~\ref{sec:conclusion} concludes with directions for further research.


\paragraph{Notation and symbols}

Table~\ref{tab:notation} collects the core notation, symbols and the
two formulation names used throughout the paper.

\begin{longtable}{ll}
  \caption{Notation}
  \label{tab:notation}\\
  \hline
  Symbol & Description\\
  \hline
  \endfirsthead
  \hline
  Symbol & Description\\
  \hline
  \endhead
  \hline
  \endfoot
  \hline
  \endlastfoot
  $\norm{\cdot}_2$ & Euclidean norm on $\mathbb{R}^d$\\
  $P$         & index set of terminals, $\abs{P}=p$\\
  $S$         & index set of Steiner points, $\abs{S}=k$\\
  $V=P\cup S$ & node index set\\
  $E$         & edge set of the complete graph $G:=G(V,E)$\\
  $E^\star$   & candidate edge set after elimination, $E^\star\subseteq E$\\
  $H$   & a minimum spanning tree of the terminals\\
  $E^\star_H$ & candidate edge set whose terminal--terminal edges belong to $H$, $E^\star_H\subseteq E^\star$\\
  $[i,j]\in E$ & edge of $G(V,E)$, with $i,j\in V$, $i<j$\\
  $E_T$& edge set of a geometric tree $T=(W,E_T)$\\
  $\{u,v\}\in E_T$&  edge of a geometric tree $T=(W,E_T)$;  unordered pair of embedded points in $u,v\in W$\\
  $n=p+k$     & total number of nodes\\
  $d$         & ambient dimension\\
  $a_i\in\mathbb{R}^d$       & coordinates of terminal $i\in P$\\
  $M$, $M_{ij}$ & big-$M$ constants\\
  $B$         & terminal box: smallest axis-aligned box containing the terminals\\
  $L^{\star}_{\le k}$ & optimal value of  the EBSTP\\
  \Mone{}, \Mtwo{} & single-commodity flow formulation and formulation with connectivity cuts\\
  $x_j\in\mathbb{R}^d$       & variable for all formulations, representing the coordinates of Steiner point $j\in S$\\
  $y_{ij}\in\{0,1\}$    & variable for all formulations, representing the edge-selection for  $[i,j]\in E$\\
  $L\in\mathbb{R}_+$         & variable for all formulations, representing the bottleneck (maximum edge length) value\\
\end{longtable}

\FloatBarrier

\section{Problem definition and basic properties}
\label{sec:problem}

Let $P:=\{1,\dots,p\}$ be the index set of terminals with distinct positions $a_i\in\mathbb{R}^d$, $i\in P$, where $p\ge 2$.\footnote{The recent conic-formulation literature of Maculan and Fampa~\citep*{pinto2023,maculan2025} denotes the number of terminals by~$n$ and the Steiner index set by $X=\{1,\dots,n-2\}$. We retain $p$ for the terminal count and reserve $n:=p+k$ for the total number of nodes.} A \emph{feasible geometric tree} is a finite tree $T=(W,E_T)$ embedded in $\mathbb{R}^d$ in which $W\subset\mathbb{R}^d$ is finite, $\{a_i:i\in P\}\subseteq W$, and at most $k$ points of $W$ are not terminals. The points in $W\setminus\{a_i:i\in P\}$ are called \emph{Steiner points}. The bottleneck value of a feasible tree is
\begin{equation}
 \textstyle  \beta(T)\;:=\;\max_{\{u,v\}\in E_T}\,\norm{u-v}_2,
  \label{eq:bottleneck}
\end{equation}
and the Euclidean bottleneck Steiner tree problem (EBSTP) is
\begin{equation}
  L^{\star}_{\le k}\;:=\;
  \min\bigl\{\beta(T)\,:\,T\text{ is a feasible geometric tree}\bigr\}.
  \label{eq:ebstp}\tag{EBSTP}
\end{equation}
Any feasible geometric tree attaining this minimum is called a \emph{bottleneck Steiner tree} of $P$.

\begin{proposition}\label{prop:monotone}
  The optimal value $L^{\star}_{\le k}$ is nonincreasing in the Steiner budget~$k$.
\end{proposition}

\begin{proof}
  The feasible set of~\eqref{eq:ebstp} grows with $k$, so its minimum cannot increase.
\end{proof}

The EBSTP is computationally challenging in general. When the Steiner budget $k$ is part of the input, the problem is \textsf{NP}-hard even in the Euclidean plane and cannot be approximated within a factor of $\sqrt{2}$ unless $\textsf{P}=\textsf{NP}$~\citep*{wang2002}. In contrast, for every fixed $k$, the planar problem can be solved in polynomial time~\citep*{bae2009exact,bae2011exact} and is fixed-parameter tractable with respect to~$k$~\citep*{bandyapadhyay2024}. The limiting case $k=0$ reduces to the bottleneck (min-max) spanning tree problem, for which any Euclidean minimum spanning tree of the terminals is optimal~\citep*{hwang1992}. This observation underlies the constructive upper bound in Section~\ref{sec:UH}.

Figure~\ref{fig:gallery} illustrates optimal bottleneck Steiner trees obtained by our formulations across different dimensions and instance families.

\begin{figure}[htbp]
  \centering
  \begin{minipage}[t]{0.48\linewidth}\centering
\begin{tikzpicture}
\begin{axis}
  [
    steinergeo,
    width = 6.0000cm,
    height = 4.8000cm,
    xmin = -0.1204, xmax = 1.2015,
    ymin = -0.0017, ymax = 1.0559,
  ]
  \addplot[btedge, no marks, forget plot] coordinates {(0.0991,0.7020) (0.2608,0.5946)};
  \addplot[btbottleneck, no marks, forget plot] coordinates {(0.5033,0.8758) (0.5237,0.5808)};
  \addplot[btedge, no marks, forget plot] coordinates {(0.9535,0.5078) (0.8482,0.6270)};
  \addplot[btedge, no marks, forget plot] coordinates {(0.9535,0.5078) (0.7042,0.3540)};
  \addplot[btedge, no marks, forget plot] coordinates {(0.7898,0.0713) (0.7042,0.3540)};
  \addplot[btedge, no marks, forget plot] coordinates {(0.1384,0.3932) (0.1258,0.1497)};
  \addplot[btedge, no marks, forget plot] coordinates {(0.1384,0.3932) (0.2608,0.5946)};
  \addplot[btedge, no marks, forget plot] coordinates {(0.1777,0.8442) (0.2608,0.5946)};
  \addplot[btedge, no marks, forget plot] coordinates {(0.8223,0.9794) (0.9820,0.9829)};
  \addplot[btedge, no marks, forget plot] coordinates {(0.9820,0.9829) (0.9661,0.8655)};
  \addplot[btedge, no marks, forget plot] coordinates {(0.5237,0.5808) (0.2608,0.5946)};
  \addplot[btedge, no marks, forget plot] coordinates {(0.5237,0.5808) (0.7042,0.3540)};
  \addplot[btedge, no marks, forget plot] coordinates {(0.8482,0.6270) (0.9661,0.8655)};
  \addplot[only marks, mark = *, mark size = 2.4pt, mark options = {fill = btterminal, draw = btterminal}, forget plot] coordinates {(0.0991,0.7020) (0.5033,0.8758) (0.9535,0.5078) (0.7898,0.0713) (0.1384,0.3932) (0.1777,0.8442) (0.8223,0.9794) (0.9820,0.9829) (0.5237,0.5808) (0.1258,0.1497) (0.8482,0.6270) (0.9661,0.8655) };
  \addplot[only marks, mark = triangle*, mark size = 3.4pt, mark options = {fill = btsteiner, draw = btsteiner}, forget plot] coordinates {(0.2608,0.5946) (0.7042,0.3540) };
\end{axis}
\end{tikzpicture}\\[1pt]
    {\small (a) uniform, $p=12$, $k=2$, $d=2$}
  \end{minipage}\hfill
  \begin{minipage}[t]{0.48\linewidth}\centering
\begin{tikzpicture}
\begin{axis}
  [
    steinergeo,
    width = 6.0000cm,
    height = 4.8000cm,
    xmin = -0.6351, xmax = 0.7025,
    ymin = -0.5239, ymax = 0.5462,
  ]
  \addplot[btedge, no marks, forget plot] coordinates {(0.2174,-0.4501) (0.2945,-0.4126)};
  \addplot[btedge, no marks, forget plot] coordinates {(0.3890,0.0766) (0.3786,-0.0622)};
  \addplot[btbottleneck, no marks, forget plot] coordinates {(0.3890,0.0766) (0.2118,0.2029)};
  \addplot[btedge, no marks, forget plot] coordinates {(0.0677,0.4717) (0.0997,0.4724)};
  \addplot[btedge, no marks, forget plot] coordinates {(0.0997,0.4724) (0.0346,0.3292)};
  \addplot[btedge, no marks, forget plot] coordinates {(0.2945,-0.4126) (0.3830,-0.3556)};
  \addplot[btedge, no marks, forget plot] coordinates {(0.3786,-0.0622) (0.5220,0.0713)};
  \addplot[btedge, no marks, forget plot] coordinates {(0.3786,-0.0622) (0.2201,-0.2113)};
  \addplot[btedge, no marks, forget plot] coordinates {(-0.3312,0.2274) (-0.2512,0.3528)};
  \addplot[btedge, no marks, forget plot] coordinates {(0.3830,-0.3556) (0.2201,-0.2113)};
  \addplot[btedge, no marks, forget plot] coordinates {(-0.3358,0.1523) (-0.3384,0.1620)};
  \addplot[btedge, no marks, forget plot] coordinates {(-0.3358,0.1523) (-0.2512,0.3528)};
  \addplot[btedge, no marks, forget plot] coordinates {(-0.4096,0.2819) (-0.4546,0.2677)};
  \addplot[btedge, no marks, forget plot] coordinates {(-0.4096,0.2819) (-0.2512,0.3528)};
  \addplot[btedge, no marks, forget plot] coordinates {(0.0346,0.3292) (-0.0479,0.4302)};
  \addplot[btedge, no marks, forget plot] coordinates {(0.0346,0.3292) (0.2118,0.2029)};
  \addplot[btedge, no marks, forget plot] coordinates {(-0.0479,0.4302) (-0.2512,0.3528)};
  \addplot[only marks, mark = *, mark size = 2.4pt, mark options = {fill = btterminal, draw = btterminal}, forget plot] coordinates {(0.2174,-0.4501) (0.3890,0.0766) (0.0677,0.4717) (0.0997,0.4724) (0.2945,-0.4126) (0.3786,-0.0622) (-0.3312,0.2274) (0.3830,-0.3556) (-0.3358,0.1523) (-0.4096,0.2819) (0.0346,0.3292) (-0.4546,0.2677) (-0.3384,0.1620) (-0.0479,0.4302) (0.5220,0.0713) };
  \addplot[only marks, mark = triangle*, mark size = 3.4pt, mark options = {fill = btsteiner, draw = btsteiner}, forget plot] coordinates {(-0.2512,0.3528) (0.2118,0.2029) (0.2201,-0.2113) };
\end{axis}
\end{tikzpicture}\\[2pt]
    {\small (b) clustered, $p=15$, $k=3$, $d=2$}
  \end{minipage}

  \vspace{1em}
  \begin{minipage}[t]{0.48\linewidth}\centering
\begin{tikzpicture}
\begin{axis}
  [
    steineraxis3d,
    width = 6.0000cm,
    height = 6.0000cm,
    xmin = 0.0688, xmax = 1.0133,
    ymin = 0.0531, ymax = 0.9976,
    zmin = 0.0879, zmax = 1.0324,
  ]
  \addplot3[btedge3d, no marks, forget plot] coordinates {(0.3932,0.1777,0.8442) (0.1258,0.1497,0.8482)};
  \addplot3[btedge3d, no marks, forget plot] coordinates {(0.0991,0.7020,0.5033) (0.3464,0.5038,0.6862)};
  \addplot3[btedge3d, no marks, forget plot] coordinates {(0.3003,0.2442,0.7720) (0.3464,0.5038,0.6862)};
  \addplot3[btedge3d, no marks, forget plot] coordinates {(0.3932,0.1777,0.8442) (0.3464,0.5038,0.6862)};
  \addplot3[btedge3d, no marks, forget plot] coordinates {(0.2314,0.8285,0.8121) (0.3464,0.5038,0.6862)};
  \addplot3[btbottleneck3d, no marks, forget plot] coordinates {(0.7898,0.0713,0.1384) (0.8253,0.2515,0.4561)};
  \addplot3[btedge3d, no marks, forget plot] coordinates {(0.3464,0.5038,0.6862) (0.6581,0.6860,0.6322)};
  \addplot3[btedge3d, no marks, forget plot] coordinates {(0.9829,0.5237,0.5808) (0.8253,0.2515,0.4561)};
  \addplot3[btedge3d, no marks, forget plot] coordinates {(0.9829,0.5237,0.5808) (0.6581,0.6860,0.6322)};
  \addplot3[btedge3d, no marks, forget plot] coordinates {(0.6270,0.9661,0.8655) (0.6581,0.6860,0.6322)};
  \addplot3[btedge3d, no marks, forget plot] coordinates {(0.9829,0.5237,0.5808) (0.8995,0.6566,0.2667)};
  \addplot3[btedge3d, no marks, forget plot] coordinates {(0.8758,0.9535,0.5078) (0.6581,0.6860,0.6322)};
  \addplot3[btedge3d, no marks, forget plot] coordinates {(0.8223,0.9794,0.9820) (0.6270,0.9661,0.8655)};
  \addplot3[only marks, mark = *, mark size = 2.4pt, mark options = {fill = btterminal, draw = btterminal}, forget plot] coordinates {(0.0991,0.7020,0.5033) (0.8758,0.9535,0.5078) (0.7898,0.0713,0.1384) (0.3932,0.1777,0.8442) (0.8223,0.9794,0.9820) (0.9829,0.5237,0.5808) (0.1258,0.1497,0.8482) (0.6270,0.9661,0.8655) (0.8253,0.2515,0.4561) (0.8995,0.6566,0.2667) (0.2314,0.8285,0.8121) (0.3003,0.2442,0.7720) };
  \addplot3[only marks, mark = triangle*, mark size = 3.4pt, mark options = {fill = btsteiner, draw = btsteiner}, forget plot] coordinates {(0.3464,0.5038,0.6862) (0.6581,0.6860,0.6322) };
\end{axis}
\end{tikzpicture}\\[1pt]
    {\small (c) uniform, $p=12$, $k=2$, $d=3$}
  \end{minipage}\hfill
  \begin{minipage}[t]{0.48\linewidth}\centering
\begin{tikzpicture}
\begin{axis}
  [
    steineraxis3d,
    width = 6.0000cm,
    height = 6.0000cm,
    xmin = -0.5314, xmax = 0.4781,
    ymin = -0.4434, ymax = 0.5661,
    zmin = -0.4886, zmax = 0.5209,
  ]
  \addplot3[btedge3d, no marks, forget plot] coordinates {(-0.4546,0.2677,0.0657) (-0.3709,0.1283,-0.0721)};
  \addplot3[btedge3d, no marks, forget plot] coordinates {(-0.4546,0.2677,0.0657) (-0.3264,0.1961,-0.0052)};
  \addplot3[btedge3d, no marks, forget plot] coordinates {(-0.3323,0.2655,0.0783) (-0.3709,0.1283,-0.0721)};
  \addplot3[btedge3d, no marks, forget plot] coordinates {(-0.0418,-0.3720,0.3354) (-0.1467,-0.3734,0.3985)};
  \addplot3[btedge3d, no marks, forget plot] coordinates {(-0.0424,-0.2264,0.4405) (-0.1467,-0.3734,0.3985)};
  \addplot3[btedge3d, no marks, forget plot] coordinates {(0.3863,-0.4240,-0.3455) (0.3697,-0.3974,-0.4083)};
  \addplot3[btedge3d, no marks, forget plot] coordinates {(0.3863,-0.4240,-0.3455) (0.2248,-0.2607,-0.0991)};
  \addplot3[btbottleneck3d, no marks, forget plot] coordinates {(-0.3323,0.2655,0.0783) (0.0021,0.2289,0.0936)};
  \addplot3[btbottleneck3d, no marks, forget plot] coordinates {(0.0636,-0.0979,0.1478) (0.2248,-0.2607,-0.0991)};
  \addplot3[btedge3d, no marks, forget plot] coordinates {(-0.0424,-0.2264,0.4405) (0.0636,-0.0979,0.1478)};
  \addplot3[btbottleneck3d, no marks, forget plot] coordinates {(0.0021,0.2289,0.0936) (0.0636,-0.0979,0.1478)};
  \addplot3[btbottleneck3d, no marks, forget plot] coordinates {(0.3010,0.3834,0.0775) (0.0021,0.2289,0.0936)};
  \addplot3[btedge3d, no marks, forget plot] coordinates {(0.3010,0.3834,0.0775) (0.3545,0.5436,0.0645)};
  \addplot3[btedge3d, no marks, forget plot] coordinates {(0.4012,0.5467,0.0809) (0.3545,0.5436,0.0645)};
  \addplot3[only marks, mark = *, mark size = 2.4pt, mark options = {fill = btterminal, draw = btterminal}, forget plot] coordinates {(-0.0424,-0.2264,0.4405) (0.3863,-0.4240,-0.3455) (0.3010,0.3834,0.0775) (0.4012,0.5467,0.0809) (-0.0418,-0.3720,0.3354) (0.3697,-0.3974,-0.4083) (-0.4546,0.2677,0.0657) (-0.1467,-0.3734,0.3985) (-0.3323,0.2655,0.0783) (-0.3264,0.1961,-0.0052) (0.3545,0.5436,0.0645) (-0.3709,0.1283,-0.0721) };
  \addplot3[only marks, mark = triangle*, mark size = 3.4pt, mark options = {fill = btsteiner, draw = btsteiner}, forget plot] coordinates {(0.0021,0.2289,0.0936) (0.0636,-0.0979,0.1478) (0.2248,-0.2607,-0.0991) };
\end{axis}
\end{tikzpicture}\\[2pt]
    {\small (d) clustered, $p=12$, $k=3$, $d=3$}
  \end{minipage}
  \caption{  
  Optimal bottleneck Steiner trees on representative instances, computed with
formulation~\Mone{} (Section~\ref{sec:formulations}). Filled blue circles are terminals,
orange triangles are Steiner points, and thick red edges are the bottleneck edges (all of
them when several tie for the maximum). Instances (a) and (b) are planar. In the three-dimensional instances (c) and (d),
a nearer edge is drawn over a farther one at every crossing. The budgets are chosen for
visual clarity and need not coincide with those of the computational study.
  }
  \label{fig:gallery}
\end{figure}

\subsection{Indexed representation on exactly \texorpdfstring{$k$}{k} Steiner points}
\label{sec:indexed}

The continuous Steiner positions and the discrete tree topology must be modeled jointly. We model both through a fixed indexed node set. Let
\begin{equation}
  S\;:=\;\{p+1,\dots,p+k\},\qquad V\;:=\;P\cup S,
\end{equation}
and let $G=(V,E)$ be the complete undirected graph on~$V$, with \[E:=\{[i,j]:i,j\in V,\,i<j\}.\] Throughout the paper, we assume that the terminal coordinates $x_i:=a_i$ for $i\in P$ are fixed; the Steiner coordinates $x_j\in\mathbb{R}^d$ for $j\in S$ are decision variables. For each coordinate $q=1,\dots,d$, set $\ell_q:=\min_{i\in P}(a_i)_q$ and $u_q:=\max_{i\in P}(a_i)_q$, and let 
\[\textstyle B:=\prod_{q=1}^d[\ell_q,u_q]
\] 
be the terminal box, the smallest axis-aligned box containing all terminals; the formulations restrict every Steiner coordinate to $B$, a restriction whose validity follows from Proposition~\ref{prop:box} below. Modeling the problem on this fixed node set requires committing to a single value of $k$, but an optimal tree may use fewer Steiner points. The next two results justify using a fixed node set: any optimal tree can be padded to exactly $k$ Steiner points without increasing the bottleneck value (Lemma~\ref{lem:padding}), so the EBSTP is equivalent to its restriction to indexed trees on~$V$ in which every Steiner point has degree at least~$2$ (Proposition~\ref{prop:indexed}). An analogous reduction for the additive Euclidean Steiner problem is due to \citet*{maculan2000}; the bottleneck objective requires a separate argument, which we give below.

\FloatBarrier

\begin{lemma}\label{lem:padding}
  For every tree $T'$ feasible for~\eqref{eq:ebstp} using $k'\le k$ Steiner points, there exists a tree $T$ feasible for~\eqref{eq:ebstp} using exactly $k$ Steiner points, with no Steiner leaves, and satisfying $\beta(T)\le\beta(T')$.
\end{lemma}

\begin{proof}
  Removing any Steiner leaf of $T'$ cannot disconnect terminals or increase the bottleneck value; iterating yields a feasible tree $T_0$ with no Steiner leaves and with $k_0\le k'$ Steiner points and $\beta(T_0)\le\beta(T')$. If $k_0=k$, set $T=T_0$. Otherwise, choose any edge $\{u,v\}\in E_{T_0}$ (one exists since $p\ge 2$) and insert $\alpha:=k-k_0$ Steiner points
  \[
   \textstyle x_j \;=\; u+\frac{j}{\alpha +1}(v-u), \qquad j=1,\dots,\alpha,
  \]
  splitting $\{u,v\}$ into a path of length $\alpha+1$. Each new edge has length $\norm{u-v}_2/(\alpha+1)\le\norm{u-v}_2$, so the bottleneck value does not increase; each inserted point has degree $2$, so no Steiner leaf is created. The resulting tree $T$ has exactly $k$ Steiner points and $\beta(T)\le\beta(T_0)\le\beta(T')$.
\end{proof}

\begin{proposition}\label{prop:indexed}
  The optimal value of~\eqref{eq:ebstp} equals the optimal value obtained by optimizing the bottleneck value~\eqref{eq:bottleneck} over trees on the indexed node set $V=P\cup S$ in which every Steiner point has degree at least~$2$.
\end{proposition}

\begin{proof}
  Any tree on $V$ with no Steiner leaves uses exactly $k$ Steiner points and is feasible for~\eqref{eq:ebstp}{ \ (if some of its nodes coincide, identify them and delete one edge from each resulting cycle; this does not increase the bottleneck value and leaves at most $k$ Steiner points)}, so the indexed optimum is no smaller than $L^{\star}_{\le k}$. The reverse inequality follows from Lemma~\ref{lem:padding}: any optimal tree for~\eqref{eq:ebstp} can be padded to a tree on $V$ with exactly $k$ Steiner points, no Steiner leaves, and the same bottleneck value.
\end{proof}

Under Proposition~\ref{prop:indexed} the no-Steiner-leaf condition can be imposed as a pre-modeling assumption, eliminating the need for activation variables on the Steiner points.

A further structural reduction constrains \emph{where} the Steiner points may be assumed to lie. Classically, some optimum keeps every Steiner point inside the convex hull $\conv(P)$ of the terminals, where $\conv(P):=\conv\{a_i:i\in P\}$; we prove this convex-hull containment for the bottleneck objective, as it both justifies the terminal box restriction used by the formulations and yields the tighter big-$M$ constants of Section~\ref{sec:bigM}.

\begin{proposition}\label{prop:box}
  There exists an optimal solution of~\eqref{eq:ebstp} in which every Steiner point lies in $\conv(P)$.
\end{proposition}

\begin{proof}
  Let $T^\star$ be an optimal tree. The Euclidean projection $\Pi$ onto the closed convex set $\conv(P)$ is nonexpansive~\citep*[\S 8.1]{boyd2004}: $\norm{\Pi(x)-\Pi(y)}_2\le\norm{x-y}_2$ for all $x,y$. Every terminal lies in $\conv(P)$, so $\Pi(a_i)=a_i$ for all $i\in P$. Project every Steiner point of $T^\star$ onto $\conv(P)$ while keeping the terminals fixed and the tree topology unchanged. By nonexpansiveness, the length of each edge weakly decreases, so the resulting tree has bottleneck value no larger than that of $T^\star$. Hence it is also optimal. The projection may map distinct vertices to the same point. As the terminals remain fixed and distinct, every coincident pair contains at least one Steiner point. Repeatedly identify any such pair. If the two vertices are adjacent, contract the resulting zero-length edge. Otherwise, identifying them creates a unique cycle, from which we delete one edge. In either case the graph remains a connected spanning tree, and the bottleneck value cannot increase: contraction removes only a zero-length edge, while deleting an edge can only remove a candidate attaining the maximum in ~\eqref{eq:bottleneck}.
  After all coincident vertices have been merged, we obtain a tree with $k'\le k$ distinct Steiner points, all lying in $\conv(P)$. Finally, apply Lemma~\ref{lem:padding} to insert $k-k'$ additional degree-two Steiner points by subdividing existing edges. Each inserted point is a convex combination of the endpoints of the subdivided edge and therefore also belongs to $\conv(P)$. By Lemma~\ref{lem:padding}, this operation does not increase the bottleneck value.
  The resulting tree has exactly $k$ Steiner points, all lying in $\conv(P)$, and attains the same optimal value as $T^\star$.
\end{proof}

As $\conv(P)\subseteq B$, Proposition~\ref{prop:box} in particular justifies restricting every Steiner coordinate to the terminal box~$B$, a restriction applied throughout the formulations of Section~\ref{sec:formulations}.

\FloatBarrier

\section{Structural differences between bottleneck and classic Euclidean Steiner trees}
\label{sec:contrast}

The reductions developed in Section~\ref{sec:problem} show that bottleneck Steiner trees have a geometry that differs from that of classical Euclidean Steiner minimal trees. While both problems allow the introduction of Steiner points, the change from a min-sum objective to a min-max objective alters the structure of optimal solutions.

In the classical Euclidean Steiner tree problem, the objective is to minimize total edge length. This objective imposes geometric restrictions: Steiner points have degree three, the incident edges meet at $120^\circ$, full Steiner topologies play a central role, and every Steiner minimal tree satisfies strong local optimality conditions (see, e.g., \citep*{fampa2016,fampa2019}). These properties underlie many of the specialized geometric algorithms developed for the problem.

The EBSTP is driven by a different objective. Instead of minimizing the sum of edge lengths, it minimizes the length of the longest edge in the tree. The optimization therefore seeks to balance long edges across the network rather than achieving local length efficiency. Thus, bottleneck-optimal trees need not satisfy the classical regularity properties of Euclidean Steiner minimal trees. In particular, they may violate the $120^\circ$ angle condition, the degree-three property of Steiner points, and the full-Steiner topology assumption under which a Steiner minimal tree on $p$ terminals contains at most $p-2$ Steiner points.

 The EBSTP shares the same broad complexity profile as the classical Euclidean Steiner tree problem \citep*{garey1977complexity}: both are \textsf{NP}-hard when Steiner points are allowed, and both reduce to a polynomially solvable spanning-tree problem when $k=0$ (Section~\ref{sec:problem}). Although the problems remain in the same computational complexity class, the geometry of their optimal solutions changes. Bottleneck-optimal trees should therefore be studied as a different object, not as a variant of Euclidean Steiner minimal trees.

In this section, we illustrate these differences on a few small instances. The bottleneck-optimal trees shown here are obtained using formulation \Mone{}, defined in Section~\ref{sec:formulations}; we use it only to generate optimal solutions and use no properties of \Mone{} in this section. Whether this additional structural flexibility translates into computational advantages at the full Steiner budget $k=p-2$ remains an empirical question, which we do not address in this illustrative discussion.

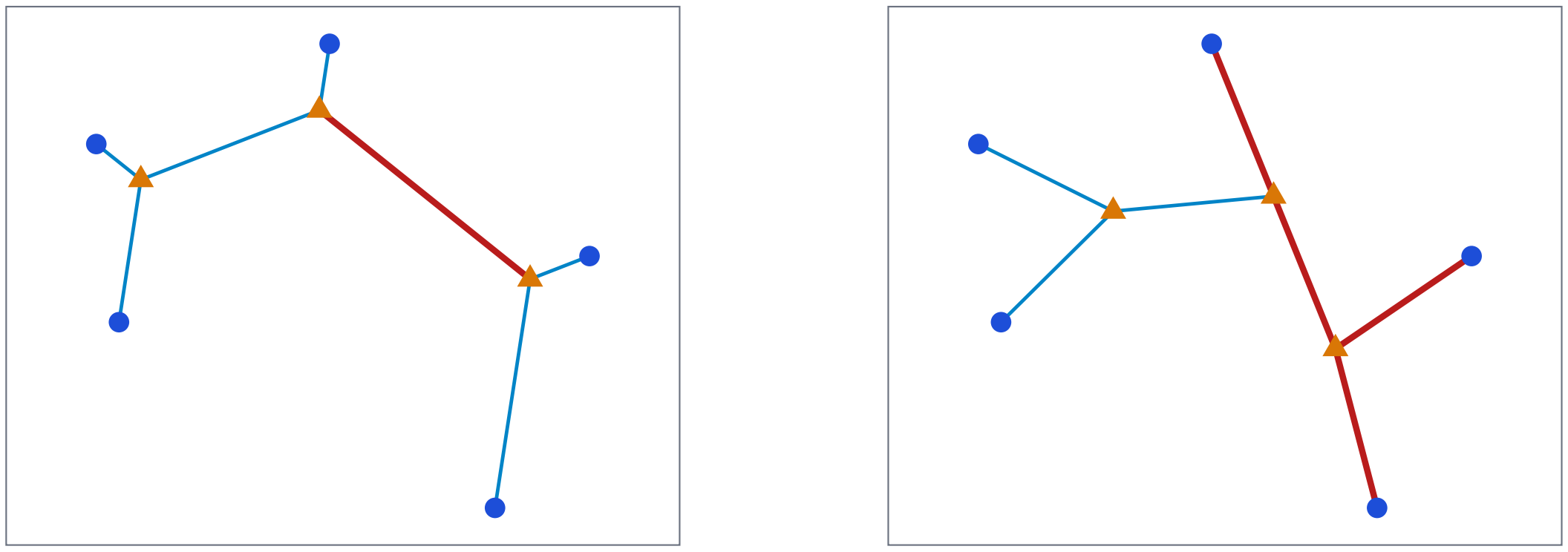
\begin{figure}[htbp]
  \centering
  \begin{minipage}[t]{0.48\linewidth}\centering
\begin{tikzpicture}
\begin{axis}
  [
    steinergeo,
    width = 6.0000cm,
    height = 4.8000cm,
    xmin = -0.0570, xmax = 1.1096,
    ymin = 0.0069, ymax = 0.9402,
  ]
  \addplot[btedge, no marks, forget plot] coordinates {(0.0991,0.7020) (0.1763,0.6402)};
  \addplot[btedge, no marks, forget plot] coordinates {(0.5033,0.8758) (0.4856,0.7606)};
  \addplot[btedge, no marks, forget plot] coordinates {(0.9535,0.5078) (0.8506,0.4678)};
  \addplot[btedge, no marks, forget plot] coordinates {(0.7898,0.0713) (0.8506,0.4678)};
  \addplot[btedge, no marks, forget plot] coordinates {(0.1384,0.3932) (0.1763,0.6402)};
  \addplot[btedge, no marks, forget plot] coordinates {(0.1763,0.6402) (0.4856,0.7606)};
  \addplot[btbottleneck, no marks, forget plot] coordinates {(0.4856,0.7606) (0.8506,0.4678)};
  \addplot[only marks, mark = *, mark size = 2.4pt, mark options = {fill = btterminal, draw = btterminal}, forget plot] coordinates {(0.0991,0.7020) (0.5033,0.8758) (0.9535,0.5078) (0.7898,0.0713) (0.1384,0.3932) };
  \addplot[only marks, mark = triangle*, mark size = 3.4pt, mark options = {fill = btsteiner, draw = btsteiner}, forget plot] coordinates {(0.1763,0.6402) (0.4856,0.7606) (0.8506,0.4678) };
\end{axis}
\end{tikzpicture}\\[2pt]
    {\small (a) minimum-total-length tree}
  \end{minipage}\hfill
  \begin{minipage}[t]{0.48\linewidth}\centering
\begin{tikzpicture}
\begin{axis}
  [
    steinergeo,
    width = 6.0000cm,
    height = 4.8000cm,
    xmin = -0.0570, xmax = 1.1096,
    ymin = 0.0069, ymax = 0.9402,
  ]
  \addplot[btedge, no marks, forget plot] coordinates {(0.0991,0.7020) (0.3327,0.5855)};
  \addplot[btbottleneck, no marks, forget plot] coordinates {(0.5033,0.8758) (0.6105,0.6116)};
  \addplot[btbottleneck, no marks, forget plot] coordinates {(0.9535,0.5078) (0.7177,0.3472)};
  \addplot[btbottleneck, no marks, forget plot] coordinates {(0.7898,0.0713) (0.7177,0.3472)};
  \addplot[btedge, no marks, forget plot] coordinates {(0.1384,0.3932) (0.3327,0.5855)};
  \addplot[btedge, no marks, forget plot] coordinates {(0.3327,0.5855) (0.6105,0.6116)};
  \addplot[btbottleneck, no marks, forget plot] coordinates {(0.6105,0.6116) (0.7177,0.3472)};
  \addplot[only marks, mark = *, mark size = 2.4pt, mark options = {fill = btterminal, draw = btterminal}, forget plot] coordinates {(0.0991,0.7020) (0.5033,0.8758) (0.9535,0.5078) (0.7898,0.0713) (0.1384,0.3932) };
  \addplot[only marks, mark = triangle*, mark size = 3.4pt, mark options = {fill = btsteiner, draw = btsteiner}, forget plot] coordinates {(0.3327,0.5855) (0.6105,0.6116) (0.7177,0.3472) };
\end{axis}
\end{tikzpicture}\\[2pt]
    {\small (b) bottleneck Steiner tree ($k:=p-2$) }
  \end{minipage}
  \caption{
  The two objectives on the same $5$ terminals with $k=p-2=3$ Steiner points: (a) the
minimum-total-length (classic Steiner) tree, computed with an FM-type MISOCP on the indexed
node set; (b) the bottleneck optimum computed with \Mone{}. Both trees use all three Steiner
points, so the objectives are compared at the same Steiner-point count. 
  The additive
    optimum minimizes total length but concentrates its length in a single
    long edge of length $0.4680$; the bottleneck optimum rebalances the tree
    so that its longest edge shrinks to $0.2852$, a $39\%$ reduction, with
    several edges tying for the new maximum. Conventions as in Figure~\ref{fig:gallery}.}
  \label{fig:hero}
\end{figure}

Beyond the rebalancing visible in Figure~\ref{fig:hero}, both optima use the same $p-2=3$ Steiner points and differ only in their objective: the bottleneck optimum is not a reweighting of the additive optimum but a structurally distinct tree.

\begin{figure}[htbp]
  \centering
  \begin{minipage}[t]{0.32\linewidth}\centering
\begin{tikzpicture}
\begin{axis}
  [
    steinergeo,
    width = 4.6000cm,
    height = 3.6800cm,
    xmin = -0.2250, xmax = 1.2250,
    ymin = -0.0800, ymax = 1.0800,
  ]
  \addplot[btbottleneck, no marks, forget plot] coordinates {(0.0000,0.0000) (0.5000,0.5000)};
  \addplot[btbottleneck, no marks, forget plot] coordinates {(1.0000,0.0000) (0.5000,0.5000)};
  \addplot[btbottleneck, no marks, forget plot] coordinates {(1.0000,1.0000) (0.5000,0.5000)};
  \addplot[btbottleneck, no marks, forget plot] coordinates {(0.0000,1.0000) (0.5000,0.5000)};
  \addplot[only marks, mark = *, mark size = 2.4pt, mark options = {fill = btterminal, draw = btterminal}, forget plot] coordinates {(0.0000,0.0000) (1.0000,0.0000) (1.0000,1.0000) (0.0000,1.0000) };
  \addplot[only marks, mark = triangle*, mark size = 3.4pt, mark options = {fill = btsteiner, draw = btsteiner}, forget plot] coordinates {(0.5000,0.5000) };
\end{axis}
\end{tikzpicture}\\[2pt]
    {\small (a) degree-four Steiner point}
  \end{minipage}\hfill
  \begin{minipage}[t]{0.32\linewidth}\centering
\begin{tikzpicture}
\begin{axis}
  [
    steinergeo,
    width = 4.6000cm,
    height = 3.6800cm,
    xmin = -0.0570, xmax = 1.1096,
    ymin = 0.0069, ymax = 0.9402,
  ]
  \addplot[btedge, no marks, forget plot] coordinates {(0.0991,0.7020) (0.3481,0.5706)};
  \addplot[btbottleneck, no marks, forget plot] coordinates {(0.5033,0.8758) (0.6105,0.6115)};
  \addplot[btbottleneck, no marks, forget plot] coordinates {(0.9535,0.5078) (0.7177,0.3472)};
  \addplot[btbottleneck, no marks, forget plot] coordinates {(0.7898,0.0713) (0.7177,0.3472)};
  \addplot[btedge, no marks, forget plot] coordinates {(0.1384,0.3932) (0.3481,0.5706)};
  \addplot[btedge, no marks, forget plot] coordinates {(0.3481,0.5706) (0.6105,0.6115)};
  \addplot[btbottleneck, no marks, forget plot] coordinates {(0.6105,0.6115) (0.7177,0.3472)};
  \draw[btangle] (axis cs:0.2574,0.6184) arc[start angle = 152.1686, end angle = 220.2309, radius = 0.1025];
  \addplot[only marks, mark = *, mark size = 2.4pt, mark options = {fill = btterminal, draw = btterminal}, forget plot] coordinates {(0.0991,0.7020) (0.5033,0.8758) (0.9535,0.5078) (0.7898,0.0713) (0.1384,0.3932) };
  \addplot[only marks, mark = triangle*, mark size = 3.4pt, mark options = {fill = btsteiner, draw = btsteiner}, forget plot] coordinates {(0.3481,0.5706) (0.6105,0.6115) (0.7177,0.3472) };
\end{axis}
\end{tikzpicture}\\[2pt]
    {\small (b) minimum angle $68.1^\circ$}
  \end{minipage}\hfill
  \begin{minipage}[t]{0.32\linewidth}\centering
\begin{tikzpicture}
\begin{axis}
  [
    steinergeo,
    width = 4.6000cm,
    height = 3.6800cm,
    xmin = -0.0796, xmax = 0.9179,
    ymin = 0.1339, ymax = 0.9319,
  ]
  \draw[bthull] (axis cs:0.0203,0.2998) -- (axis cs:0.3412,0.1918) -- (axis cs:0.4215,0.1889) -- (axis cs:0.8180,0.2669) -- (axis cs:0.7177,0.7861) -- (axis cs:0.4246,0.8768) -- cycle;
  \addplot[btedge, no marks, forget plot] coordinates {(0.7177,0.7861) (0.6605,0.6565)};
  \addplot[btbottleneck, no marks, forget plot] coordinates {(0.7177,0.7861) (0.4246,0.8768)};
  \addplot[btbottleneck, no marks, forget plot] coordinates {(0.6605,0.6565) (0.6823,0.3505)};
  \addplot[btedge, no marks, forget plot] coordinates {(0.8180,0.2669) (0.6823,0.3505)};
  \addplot[btbottleneck, no marks, forget plot] coordinates {(0.4215,0.1889) (0.1147,0.1889)};
  \addplot[btbottleneck, no marks, forget plot] coordinates {(0.4215,0.1889) (0.6823,0.3505)};
  \addplot[btedge, no marks, forget plot] coordinates {(0.3412,0.1918) (0.1147,0.1889)};
  \addplot[btedge, no marks, forget plot] coordinates {(0.0203,0.2998) (0.1147,0.1889)};
  \addplot[only marks, mark = *, mark size = 2.4pt, mark options = {fill = btterminal, draw = btterminal}, forget plot] coordinates {(0.7177,0.7861) (0.6605,0.6565) (0.8180,0.2669) (0.4246,0.8768) (0.4215,0.1889) (0.3412,0.1918) (0.0203,0.2998) };
  \addplot[only marks, mark = triangle*, mark size = 3.4pt, mark options = {fill = btsteiner, draw = btsteiner}, forget plot] coordinates {(0.1147,0.1889) (0.6823,0.3505) };
\end{axis}
\end{tikzpicture}\\[2pt]
    {\small (c) Steiner point outside $\conv(P)$}
  \end{minipage}
  \caption{Counterexamples to the classic Steiner regularities, computed with \Mone{}: (a) the unit
square with $k=1$, where a degree-four Steiner point at the center with four $90^\circ$
angles violates the degree-three, $120^\circ$-angle and full-topology conditions; (b) a random
instance ($p=5$, $k=3$) whose minimum Steiner angle, marked by the arc, is $68.1^\circ$;
(c) a random instance ($p=7$, $k=2$) whose returned optimum places a Steiner point outside
$\conv(P)$ (dashed), by about $9\%$ of the largest per-axis extent of the terminals.
Conventions as in Figure~\ref{fig:gallery}.}
  \label{fig:propviol}
\end{figure}
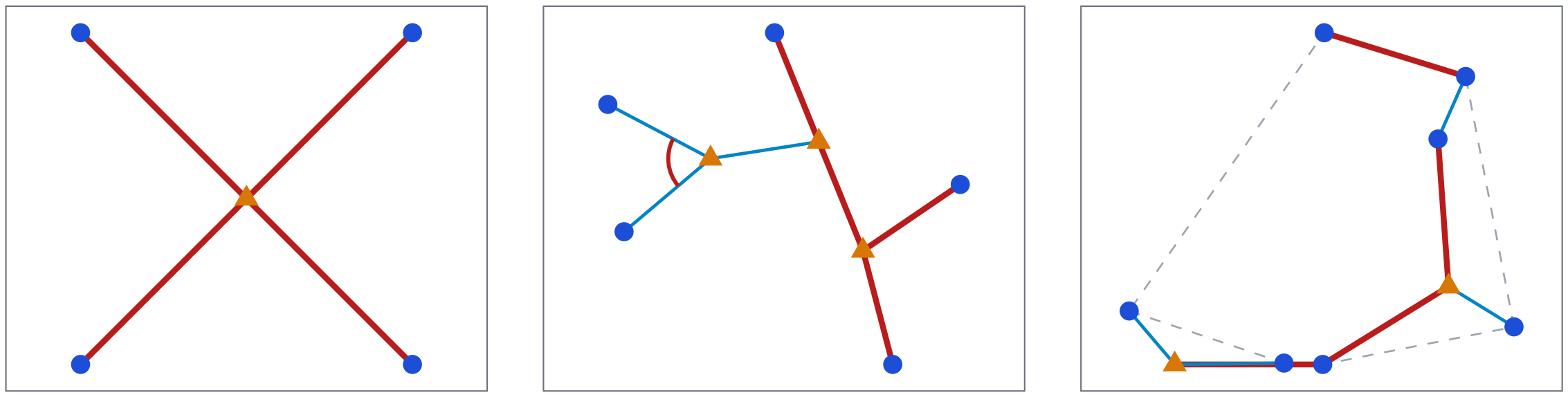

Figure~\ref{fig:propviol} gives a counterexample to each of the classic regularities, one in each of (a)--(c). On the unit square (Figure~\ref{fig:propviol}(a)) the bottleneck optimum places a single Steiner point of degree four at the center, with four $90^\circ$ angles: the maximum-degree-three and $120^\circ$ angle conditions both fail, and the tree is not a full Steiner topology (a full topology on $p=4$ terminals would use two Steiner points of degree three). On the random instance of Figure~\ref{fig:propviol}(b) the minimum angle at a Steiner point is $68.1^\circ$, far below $120^\circ$. In Figure~\ref{fig:propviol}(c) an optimal solution returned by the solver places a Steiner point outside $\conv(P)$. This does not contradict Proposition~\ref{prop:box}, whose projection argument establishes the existence of at least one optimal solution with all Steiner points in $\conv(P)$, but does not require every optimal solution to satisfy this property. The Steiner point shown in the figure nevertheless remains inside the terminal box $B\supseteq\conv(P)$. The departures from classic Steiner structure confirms that the bottleneck problem should be viewed as a distinct combinatorial optimization problem rather than a reweighted version of the additive one.

\FloatBarrier

\section{Mixed-integer second-order cone formulations}
\label{sec:formulations}

{ As shown in Section~\ref{sec:contrast}, optimal bottleneck trees need not satisfy the angle, degree, and topology conditions of Steiner minimal trees. The formulations in this section therefore impose none of these conditions.} We can regard these formulations as bottleneck adaptations of the conic Steiner formulation FM of \citet*{fampa2001,fampa2004}, which itself is a convexification of the nonconvex MINLP formulation MMX of \citet*{maculan2000}.

MMX uses a graph $G(V,E)$ whose vertex set $V$ contains $p$ given terminals and $p-2$ Steiner points. The set of edges $E$ connects each terminal to all Steiner points, and all Steiner points to each other. Terminal-Steiner and Steiner-Steiner edges carry binary selection variables $y_{ij}$. The position of Steiner point $j$ is represented by the vector variable $x_j\in\mathbb{R}^d$. The additive objective $\sum_{[i,j]\in E}\norm{x_i-x_j}_2\,y_{ij}$ is minimized subject to degree-one constraints on terminals and degree-three constraints on Steiner points, producing a full Steiner tree on $p$ terminals with exactly $p-2$ Steiner points~\citep*{maculan2000}. The MINLP formulation is nonconvex because each objective term couples the Euclidean distance $\norm{x_i-x_j}_2$, a nonlinear function of the continuous Steiner coordinates, with the binary selection variable $y_{ij}$.

In FM, Fampa and Maculan introduce auxiliary edge-length variables $d_{ij}\ge 0$ and replace each such coordinate-binary product by the linear-in-$d$ objective term $d_{ij}$ together with the second-order cone inequality
$   
  \norm{x_i-x_j}_2 \;\le\; d_{ij} + M\,(1-y_{ij})$, for $[i,j]\in E,
$
which is equivalent to the original disjunction at $y_{ij}\in\{0,1\}$ and convex when relaxed~\citep*{fampa2004}. The constant $M$ bounds the Euclidean diameter of the terminal set.

Two adaptations are required to obtain a bottleneck-objective formulation for the EBSTP defined in Section~\ref{sec:problem}. First, the additive objective $\sum d_{ij}$ is replaced by a single epigraph variable $L\ge 0$ acting on every selected edge through
$  
  \norm{x_i-x_j}_2 \;\le\; L + M\,(1-y_{ij})$, for $[i,j]\in E,
$ 
making the formulation a min-max one in the epigraph variable $L$, which represents the height of the epigraph of the longest-edge function. The auxiliary variables $d_{ij}$ of FM are no longer needed. The edge-specific constants $M_{ij}$ that replace $M$ in the common block are derived in Section~\ref{sec:bigM}. Second, the degree-one/degree-three structure of MMX, which forces exactly $p-2$ Steiner points and a full Steiner tree, is dropped, since the EBSTP allows any tree topology with at most $k$ Steiner points; we replace it by the cardinality constraint $\sum y_{ij}=p+k-1$ together with the no-Steiner-leaf constraints justified by Proposition~\ref{prop:indexed}. The two formulations presented below share this geometric core and differ only in how connectivity of the selected edge set is enforced.

\subsection{Common block}

The formulations developed in this section rely on a single structural requirement that links the geometric results of Section~\ref{sec:problem} to the optimization formulations. Rather than proving correctness separately for each modeling ingredient, we isolate the property that every valid candidate edge set and collection of big-$M$ constants must satisfy.

Let $E^\star\subseteq E$ be a candidate edge set and let $M_{ij}\ge 0$ be constants associated with the candidate edges. We require that they satisfy the following property.

\begin{property}\label{propP}
  There exists an optimal tree $T^\star$ for~\eqref{eq:ebstp} represented over the indexed node
  set $V$ with node positions $\bar x_i$, for $i\in V$, such that
  \begin{itemize}
    \item all Steiner points of $T^\star$ lie in $\conv(P)$;
      \vspace{-5pt}
    \item { every Steiner point of $T^\star$ has degree at least two in $T^\star$;}
      \vspace{-5pt}
    \item every edge of $T^\star$ belongs to $E^\star$;
      \vspace{-5pt}
    \item every candidate edge $[i,j]\in E^\star$ satisfies
      $\norm{\bar x_i-\bar x_j}_2\le L^{\star}_{\le k}+M_{ij}$.
  \end{itemize}
\end{property}

Property~\ref{propP} is the key correctness condition of the paper. It guarantees that the modeling restrictions introduced to strengthen the formulations do not remove all optimal solutions. In particular, every preprocessing step developed later (bounding Steiner coordinates, eliminating candidate edges, and tightening the big-$M$ constants) is designed to preserve at least one optimal tree satisfying Property~\ref{propP}. Once this property holds, the exactness proofs of both formulations reduce to showing that connectivity and tree structure are modeled correctly.

We can satisfy Property~\ref{propP} by taking $E^\star=E$ and $M_{ij}=\diam(B)$, the Euclidean diameter of the terminal box $B$, for every edge. This choice yields correct but weak formulations. The purpose of the strengthening techniques developed in Section~\ref{sec:strengthenings} is to replace this naive choice by substantially tighter candidate edge sets and edge-type-specific constants while maintaining Property~\ref{propP}.

The geometric results of Section~\ref{sec:problem} justify preserving at least one optimal solution. Proposition~\ref{prop:indexed} ensures the existence of an optimal solution without Steiner leaves, while Proposition~\ref{prop:box} guarantees the existence of an optimal solution whose Steiner points lie in $\conv(P)$. 
These results establish the existence of an optimal tree satisfying the structural requirements
encoded in Property~\ref{propP}{ : the tree constructed in the proof of
Proposition~\ref{prop:box} has exactly $k$ Steiner points, all in $\conv(P)$, and no Steiner
leaves, as it is obtained by applying Lemma~\ref{lem:padding}}.
All formulations share the variables $L\in\mathbb{R}$, $y_{ij}\in\{0,1\}$, for all $[i,j]\in E^\star$, and $x_j\in\mathbb{R}^d$, for all $j\in S$, and the following constraints. Recall that, for notational convenience, we assume throughout that
$x_i:=a_i$ for all $i\in P$.

\begin{align}
  & \norm{x_i-x_j}_2 \le L + M_{ij}(1-y_{ij}),         && [i,j]\in E^\star, \label{eq:soc}\\
  & \sum_{[i,j]\in E^\star} y_{ij} \;=\; p+k-1, \label{eq:edgecount}\\
  &\sum_{i:[i,j]\in E^\star} y_{ij} + \sum_{i:[j,i]\in E^\star} y_{ji} \;\ge\; 2,     && j\in S, \label{eq:noleaf} \\
  & L\ge 0, \\
  & y_{ij}\in\{0,1\},                                  && [i,j]\in E^\star,\\
  & x_j\in\mathbb{R}^d,\quad \ell_q\le(x_j)_q\le u_q,  && j\in S,\; q=1,\dots,d. \label{eq:x}
\end{align}
We minimize $L$ in both formulations. Each of the cone constraints~\eqref{eq:soc} is represented as a standard cone constraint~\citep*{lubin2023}. The variable bounds in~\eqref{eq:x} confine each Steiner point to the terminal box $B$ introduced in Section~\ref{sec:indexed}, in accordance with Proposition~\ref{prop:box}. Constraint~\eqref{eq:edgecount} fixes the cardinality of the selected edge set to $\abs{V}-1$; together with connectivity, this guarantees that the selected subgraph is a tree. Constraint~\eqref{eq:noleaf} enforces the no-Steiner-leaf property established in Proposition~\ref{prop:indexed}. We refer to the system~\eqref{eq:soc}--\eqref{eq:x} as the \emph{common block} and denote it by~$\mathcal{C}$.

\subsection{Formulation \Mone{}: single-commodity flow formulation}\label{sec:m1}

Fix a root terminal $r\in P$. Let $A:=\{(i,j):[i,j]\in E^\star\}\cup\{(j,i):[i,j]\in E^\star\}$ be the bidirected arc set induced by $E^\star$.
Let $f_{ij}\ge 0$ represent the flow on arc $(i,j)\in A$. The root supplies one unit of flow to each of the remaining indexed nodes, for a total supply of $p+k-1$ units. Connectivity is enforced through the standard single-commodity flow formulation~\citep*{magnanti1995}, which routes this flow from the root to every other node.

\begin{model}\label{mod:m1}
  \begin{align}
    (\hypertarget{Idenonetarget}{\Idenone})\qquad \min\;\; & L \notag\\
    \text{s.t.}\;\; & \mathcal{C}, \notag\\
    & \sum_{i:(r,i)\in A} f_{ri} \;-\; \sum_{i:(i,r)\in A} f_{ir} \;=\; p+k-1, \label{eq:m1root}\\
    & \sum_{i:(i,j)\in A} f_{ij} \;-\; \sum_{i:(j,i)\in A} f_{ji} \;=\; 1,
    && j\in V\setminus\{r\}, \label{eq:m1bal}\\
    & 0\;\le\; f_{ij} \;\le\; (p+k-1)\, y_{ij},
    && [i,j]\in E^\star, \label{eq:m1cap1}\\
    & 0\;\le\; f_{ji} \;\le\; (p+k-1)\, y_{ij},
    && [i,j]\in E^\star. \label{eq:m1cap}
  \end{align}
\end{model}

\begin{proposition}\label{prop:m1}
  The optimal value of {\rm \Idenone} equals $L^{\star}_{\le k}$.
\end{proposition}

\begin{proof}
  Let $(x,y,f,L)$ be an integer-feasible solution of {\rm \Idenone}, and define $E(y):=\{[i,j]\in E^\star : y_{ij}=1\}$.

  We first show that $(V,E(y))$ is connected. Suppose, to the contrary, that it is disconnected. Then some connected component $C\subseteq V\setminus\{r\}$ is separated from the root. As no edge in the cut $\delta(C):=\{[i,j]\in E^\star:\abs{\{i,j\}\cap C}=1\}$ is selected, $y_{ij}=0$ for all $[i,j]\in \delta(C)$, and therefore $f_{ij}=f_{ji}=0$ by~\eqref{eq:m1cap1}--\eqref{eq:m1cap}. Hence no flow can enter $C$, contradicting the balance constraints~\eqref{eq:m1bal}, which require each node in $V\setminus\{r\}$ to receive one unit of net inflow. Thus $(V,E(y))$ is connected. Together with~\eqref{eq:edgecount},
  which enforces $|E(y)|=|V|-1$, this implies that $(V,E(y))$ is a tree.
  Constraint~\eqref{eq:noleaf} further guarantees the no-Steiner-leaf property. Finally, for every selected edge $[i,j]\in E(y)$, constraint~\eqref{eq:soc} yields $\norm{x_i-x_j}_2\le L$. Therefore $(x,y)$ defines a feasible indexed Steiner tree with bottleneck value at most $L$, and Proposition~\ref{prop:indexed} implies $L\ge L^{\star}_{\le k}$. Consequently, the optimal value of {\rm \Idenone} is at least $L^{\star}_{\le k}$.

    For the reverse inequality, let $T^\star$ be the optimal tree guaranteed by
  Property~\ref{propP}. Thus $T^\star$ has bottleneck value $L^{\star}_{\le k}$, all Steiner
  points lie in $\conv(P)${ \ and have degree at least two}, all edges belong to
  $E^\star$, and $\norm{x_i-x_j}_2\le L^{\star}_{\le k}+M_{ij}$ for every $[i,j]\in E^\star$.
  Set $y$ equal to the incidence vector of $T^\star$, { let $f$ send one unit of flow
  from $r$ to every other node along the unique path of $T^\star$, so that $f$ uses only
  selected edges and carries at most $p+k-1$ units on each,} and let $L=L^{\star}_{\le k}$.
  Constraints~\eqref{eq:m1root}--\eqref{eq:m1cap} are then satisfied{ ;
  \eqref{eq:edgecount} holds because $T^\star$ has $p+k-1$ edges, \eqref{eq:noleaf} because no
  Steiner point of $T^\star$ is a leaf, and~\eqref{eq:x} because $\conv(P)\subseteq B$}.
   If $[i,j]\in E(y)$, so that $y_{ij}=1$, constraint~\eqref{eq:soc} reduces to
  $\norm{x_i-x_j}_2\le L=L^{\star}_{\le k}$,
  which holds because $T^\star$ has bottleneck value $L^{\star}_{\le k}$. If $[i,j]\notin E(y)$, so that $y_{ij}=0$, constraint~\eqref{eq:soc} becomes
  $\norm{x_i-x_j}_2\le L^{\star}_{\le k}+M_{ij}$,
  which is precisely the metric condition in Property~\ref{propP}. Hence $(x,y,f,L)$ is feasible for {\rm \Idenone} with objective value $L^{\star}_{\le k}$. Therefore the optimal value of {\rm \Idenone} is at most $L^{\star}_{\le k}$.
  Combining the two inequalities yields the result.
\end{proof}

\subsection{Formulation \Mtwo{}: connectivity cuts}
\label{sec:m2}

Connectivity can also be enforced via cuts. For $C\subseteq V$, let $\delta(C):=\{[i,j]\in E^\star:\abs{\{i,j\}\cap C}=1\}$ be the cutset induced by $C$.

\begin{model}\label{mod:m2}
  \begin{align}
    (\hypertarget{Identitytarget}{\Identity})\qquad \min\;\; & L \notag\\
    \text{s.t.}\;\; & \mathcal{C}, \notag\\
    & \sum_{[i,j]\in\delta(C)} y_{ij} \;\ge\; 1,
    && \forall\, \emptyset\ne C\subseteq V\setminus\{r\}. \label{eq:m2cut}
  \end{align}
\end{model}

\begin{proposition}\label{prop:m2}
  The optimal value of {\rm \Identity} equals $L^{\star}_{\le k}$.
\end{proposition}

\begin{proof}
  Let $(x,y,L)$ be an integer-feasible solution for {\rm \Identity} and let $E(y):=\{[i,j]\in E^\star: y_{ij}=1\}$. If $(V,E(y))$ were disconnected, some component $C\subseteq V\setminus\{r\}$ would satisfy $\sum_{[i,j]\in\delta(C)} y_{ij}=0$, violating~\eqref{eq:m2cut}. Hence $(V,E(y))$ is connected; together with~\eqref{eq:edgecount} it is a tree on~$V$. The no-Steiner-leaf condition~\eqref{eq:noleaf} is enforced and every selected edge satisfies $\norm{x_i-x_j}_2\le L$ by~\eqref{eq:soc}, so by Proposition~\ref{prop:indexed} the {\rm \Identity} optimum is at least $L^{\star}_{\le k}$. Conversely, the construction used in the proof of Proposition~\ref{prop:m1} (set $y$ to the incidence vector of the indexed representation of the optimal tree~$T^{\star}$ and $L:=L^{\star}_{\le k}$) yields a feasible {\rm \Identity} solution: cardinality holds by construction, and every cut $\delta(C)$ contains at least one edge of~$T^{\star}$.
\end{proof}

\paragraph{Relation between \Mone{} and \Mtwo{}} Propositions~\ref{prop:m1} and~\ref{prop:m2} establish that both {\rm \Idenone} and {\rm \Identity} are exact formulations of the EBSTP. Moreover, their integer feasible solutions are in
one-to-one correspondence under the projection $(x,y,f,L)\mapsto(x,y,L)$, so both formulations differ only in their representation of connectivity.
The formulations differ in their continuous relaxations. Let $P_{\rm{SCF}}$ denote the projection onto the $(x,y,L)$-space of the single-commodity flow relaxation of {\rm \Idenone} and let $P_{\rm{CUT}}$ denote the relaxation of {\rm \Identity} obtained by including all connectivity cuts~\eqref{eq:m2cut}. Classical polyhedral results imply that $P_{\rm{CUT}}\subseteq P_{\rm{SCF}}$, that is, the connectivity-cut relaxation is at least as strong as the single-commodity flow relaxation~\citep*{magnanti1995}. Consequently, if all cuts were enforced explicitly, {\rm \Identity} would provide a bound no weaker than that of {\rm \Idenone}.

In practice, however, the cut family~\eqref{eq:m2cut} is exponential in size and is therefore generated dynamically. The relaxation available at a branch-and-bound node contains only the subset of cuts separated up to that point and may thus be weaker than the full connectivity-cut relaxation.

{ 
A multi-commodity flow formulation, with one unit commodity from the root to every other node,
has a continuous relaxation whose projection onto the $(x,y,L)$-space coincides with
$P_{\rm CUT}$~\citep*{magnanti1995}, at the cost of $2m(p+k-1)$ flow variables. For the EBSTP
this stronger connectivity relaxation barely changes the bound on $L$: on the $360$ instances of
test set~C with $p\le25$ and $k\le5$ the median
root gap~\eqref{eq:rootgap} is $0.424$ for {\rm \Idenone} and $0.421$ for the
multi-commodity formulation at $k=1$, and identical to three decimals at $k=3$ and $k=5$, as the
root bound of both formulations is determined by the global lower bounds of
Section~\ref{sec:diam} on these instances (Section~\ref{sec:relaxation}). The multi-commodity formulation raises the root gap  by more
than $0.01$ on only one of these instances, with $k=1$. In branch-and-bound on test set~B the multi-commodity formulation
explores a median of $1.7$ times fewer nodes but is slower by a median factor of $3.6$ and solves one instance fewer.
}

We separate the connectivity cuts~\eqref{eq:m2cut} as follows. Given a fractional solution $y$, a most-violated cut among those separating a terminal from the root is detected by solving, for each terminal $t\ne r$, a maximum-flow problem on the support graph with edge capacities $y_{ij}$; if the maximum-flow value is below $1-\varepsilon$, the source-side breadth-first search on the residual graph yields a violated cut. At integer $y$ the connected components of the support graph provide a complete separation, which is what exactness requires. {  We add the cuts as lazy constraints during branch-and-cut; the separation routine is called both
at integer solutions and at optimal fractional solutions of the node relaxations. Restricting it
to integer solutions makes {\rm \Identity} slower by a median factor of $1.21$ on test set~B and
leaves $14$ instances unsolved that the default solves.} Throughout, the symbol {\rm \Identity} refers to Formulation~{\rm \Identity} solved in this way, with connectivity cuts separated lazily; this is the variant compared to {\rm \Idenone} in Section~\ref{sec:mainexp}. The implementation details are described in Section~\ref{sec:implementation}.

{ The strength of the complete relaxation does not determine the comparison between the formulations. With both global
lower bounds, minimizing $L$ over the multi-commodity relaxation, and hence over $P_{\rm CUT}$, gives the same value as
the relaxation of {\rm \Idenone} on all $840$ instances of test set~C with $k\ge3$ (Section~\ref{sec:relaxation}). The
difference arises in the branch-and-bound search. On the instances of test set~C that both solve and that take {\rm \Idenone}, run on
one thread, at least $1$\,s, {\rm \Identity} explores about as many nodes at $d=2$ (median ratio $1.01$) but
spends a median of $1.9$ times as long per node. At $d=3$ and $d=4$ it explores a median of $4.1$ and $3.3$ times
as many nodes and spends $2.5$ and $2.8$ times as long per node (Section~\ref{sec:mainexp}). In the plane the per-node ratio grows with the
budget, from $1.1$ at $k=3$ to $4.5$ at $k=9$, while the median number of cuts grows from $407$ to \num{3300}. Each node of
{\rm \Identity} pays for the separation and for the lazy-constraint mode of the solver, which avoids reductions and
transformations that are incompatible with constraints added during the search~\citep*{gurobi2025}. The relaxation
at a node of {\rm \Identity} contains only the cuts separated so far.}

\subsection{Formulation dimensions}

\begin{table}[!h]
  \centering
  \caption{Dimensions of the initial formulations on the candidate edge set $E^\star$, where $m=\abs{E^\star}$ and $n=p+k$. The counts consider the complete candidate edge set. Edge elimination (Section~\ref{sec:elim}) reduces $m$ in practice. The bounds on the variables, the $k-1$ symmetry-breaking inequalities of scheme S1 (Section~\ref{sec:symm}) and the dynamically generated connectivity cuts (\Mtwo)
  are not included.}
  \label{tab:dimensions}
  \small
  \renewcommand*{\arraystretch}{1.3}
  \begin{tabular}{c c c c c}
    \toprule
    &\multicolumn{2}{c}{variables\phantom{vv}}& \multicolumn{2}{c}{constraints\phantom{vvvvv}}\\
    & binary & continuous & SOC & linear\\
    \midrule
    {\rm \Idenone}   & $m$ & $dk+1+2m$ & $m$ & $2m+n+k+1$\\
    {\rm \Identity}   & $m$ & $dk+1$    & $m$ & $k+1$\\
    \bottomrule
  \end{tabular}
\end{table}

Table~\ref{tab:dimensions} summarizes the sizes of the
two formulations proposed above as a function of the numbers of candidate edges $m=\abs{E^\star}$ and indexed nodes $n=p+k$.

Formulation {\rm \Idenone} augments the common block with $2m$ flow variables and $2m$ capacity constraints to obtain a static certificate of connectivity. Formulation {\rm \Identity} has the smallest initial formulation, with connectivity enforced through cuts generated dynamically during branch-and-cut.
The dimensions reported in Table~\ref{tab:dimensions} correspond to the initial formulations, before any connectivity
cuts are added during the solution process. They also do not include linear inequalities derived from bounds on the variables.

{ 
\subsection{Relation to the formulations MMX and FM}
\label{sec:relation}

The formulations of this section reuse the modeling structure of MMX and FM: continuous Steiner
coordinates $x_j$, binary edge-selection variables $y_{ij}$, and second-order cone constraints
that are switched on and off by a big-$M$ term. The epigraph representation of the bottleneck
objective, the cardinality constraint~\eqref{eq:edgecount}, the single-commodity flow
constraints and the connectivity cuts are standard modeling tools as well, and we do not claim
novelty for any of them individually. We chose them because standard MISOCP solvers support
them directly. The differences are the following.

\paragraph{Objective coupling}
{ In FM, each edge has its own length variable $d_{ij}$, the objective is
$\sum d_{ij}$, and constraint~\eqref{eq:soc} uses $d_{ij}$ instead of $L$. In {\rm \Idenone} and
{\rm \Identity}, all edges share the single variable $L$; the only difference between the two
inequalities is the replacement of $d_{ij}$ by $L$, but the big-$M$ term now relaxes the length
constraint of every edge against the same variable, and the consequences for the continuous
relaxation are different.} In both cases the relaxation is weak, because at fractional
$y_{ij}$ the big-$M$ term $M_{ij}(1-y_{ij})$ loosens the cone constraint. In FM, each $d_{ij}$ is
bounded below by its own cone constraint, and the objective sums these bounds. In the bottleneck
formulation, $L$ is bounded below by the cone constraints of all edges simultaneously, and
Section~\ref{sec:relaxation} shows that this bound can be exactly zero: the relaxation always
has an optimal solution in which all Steiner points coincide
(Proposition~\ref{prop:symrelax}), and Theorem~\ref{thm:degenerate}
characterizes the instances on which it attains $L=0$, referred to henceforth as
\emph{degenerate instances}. This occurs in all $840$ instances of test set~C with $k\ge3$, but in
none of those with $k=1$. The global lower bounds of Section~\ref{sec:diam} raise the relaxation
value directly. Proposition~\ref{prop:rootLB} shows that they determine
the value exactly on degenerate instances. The ordering and the
terminal--terminal edge elimination leave the zero-valued construction
feasible, as shown in Section~\ref{sec:relaxation}.

\paragraph{Topology and connectivity}
MMX and FM encode a full Steiner topology: exactly $p-2$ Steiner points, terminals of degree one
connected only to Steiner points, and Steiner points of degree three. Connectivity also uses the indexing constraints
\(\sum_{i\in S:\,i<j}y_{ij}=1\) for
\(j\in S\setminus\{p+1\}\), which make the Steiner subgraph a tree.
Together with the terminal degree-one constraints, these ensure
connectivity without flow variables~\citep*{maculan2000,fampa2001,fampa2004}. The EBSTP admits any Steiner budget $k$, including
$k>p-2$, terminal--terminal edges, and Steiner points of any degree at least two
(Sections~\ref{sec:problem} and~\ref{sec:contrast}), so none of these degree constraints can be
imposed. We therefore replace them by the cardinality constraint~\eqref{eq:edgecount}, the
no-Steiner-leaf constraints~\eqref{eq:noleaf} justified by Proposition~\ref{prop:indexed}, and
an explicit connectivity formulation, for which Sections~\ref{sec:m1} and~\ref{sec:m2} present
the two standard choices. With all cuts~\eqref{eq:m2cut} present, the relaxation of
{\rm \Identity} is at least as strong as that of {\rm \Idenone}, but the cuts are separated
lazily and {\rm \Idenone} is the faster formulation in our experiments
(Section~\ref{sec:mainexp}). On an instance degenerate for {\rm \Idenone}, a stronger connectivity
formulation leaves the bound with the global lower bounds unchanged
whenever its bare relaxation value does not exceed $\underline{L}:=\max\{\diam(P)/(p+k-1),\tfrac12 e_k\}$, where $\diam(P) = \max_{i',j'\in P}\norm{a_{i'}-a_{j'}}_2$.
This includes instances that are also degenerate for the stronger
formulation, and holds on all $840$ instances of test set~C with
$k\ge3$ in our computations (Section~\ref{sec:relaxation}).

\paragraph{Validity of preprocessing}
The constants~\eqref{eq:bigMdef} for terminal--Steiner and Steiner--Steiner edges coincide with
constants used for the additive problem (Section~\ref{sec:bigM}); the constant
$\norm{a_i-a_j}_2$ for terminal--terminal edges has no counterpart in FM, which excludes these
edges. For the additive problem the constants are valid because every Steiner minimal tree has
its Steiner points in $\conv(P)$. This is not true of every optimal solution of the EBSTP
(Figure~\ref{fig:propviol}(c)); Proposition~\ref{prop:box} shows that \emph{some} optimal
solution has this property, which is all that validity requires. The terminal box, the MST restriction and short-edge fixing of
Section~\ref{sec:elim}, and the constants of Section~\ref{sec:bigM}
are justified by preserving at least one optimal solution.
Basic length-based edge elimination~\eqref{eq:elim} preserves every
optimal tree. Property~\ref{propP} states the joint condition needed
for these choices to preserve the optimal value.
}


\section{Structural strengthening techniques}
\label{sec:strengthenings}

The formulations of Section~\ref{sec:formulations} are not yet computationally practical. A naive choice of big-$M$ constants together with a complete candidate edge set produces weak relaxations and unnecessarily large formulations. This section introduces the strengthening techniques incorporated into the baseline formulation.

\subsection{Terminal box for Steiner coordinates}
\label{sec:terminal-box}

The common block restricts each Steiner coordinate to the terminal box $B$ introduced in Section~\ref{sec:indexed} through the variable bounds~\eqref{eq:x}. Proposition ~\ref{prop:box} guarantees the existence of an optimal solution whose Steiner points all lie in $\conv(P)$, and as $\conv(P)\subseteq B$, these bounds preserve at least one optimal solution even though $\conv(P)$ itself generally cannot be represented by simple axis-aligned bounds.

Consequently, $B$ serves as the modeling region for the formulations, while the stronger convex-hull containment property of Proposition ~\ref{prop:box} is exploited in Section~\ref{sec:bigM} to derive tighter edge-type-specific big-$M$ constants.

{ Restricting the Steiner points to $\conv(P)$ itself, instead of the terminal box $B$, shrinks the
feasible region of the continuous relaxation, and we tested two ways of modeling it. The first,
using the convex-combination constraints, writes every Steiner point as a convex combination of the
terminals, $x_j=\sum_{i\in P}\lambda_{ij}a_i$ with $\lambda_{ij}\ge0$ and $\sum_{i\in P}\lambda_{ij}=1$
for $j\in S$. This extended formulation of $\conv(P)$ adds $kp$ variables and $k(d+1)$ equality
constraints. The second imposes the H-representation of $\conv(P)$, that is, every facet-defining
inequality, on every Steiner point; it adds no variables but $k$ inequalities per facet and, when
$\conv(P)$ is not full-dimensional, the equations of its affine hull, the smallest affine subspace that contains $P$. We computed the facets with
the monotone chain algorithm of \citet*{andrew1979} in the plane and, for $d\ge3$, with the double
description method \citep*{fukuda1996} as implemented in cddlib \citep*{cddlib}, called in exact
rational arithmetic through the Julia package Polyhedra.jl \citep*{legat2023polyhedral}. In the
worst case the number of facets is $O(p^{\lfloor d/2\rfloor})$, attained by cyclic polytopes
\citep*{mcmullen1970}, and it is linear in $p$ for $d\le3$.

Neither variant strengthens the root relaxation, and their effect on the search is small. With the convex-combination constraints, the root bounds of the $360$ instances of
test set~C with $p\le25$ and $k\le5$ agree with those of the terminal box within $8\times10^{-7}$. In branch-and-bound on test set~B, the convex-combination and the H-representation constraints each solve
all $240$ instances, one more than the terminal box, which leaves one instance at the time limit; all three reach the same
optimal values, and for each of them the median ratio of its node count to that of the terminal box is $1.00$.
With the convex-combination constraints, the median solve time is unchanged in the plane and
$5\%$ shorter at $d=3$; with the H-representation, it is unchanged in both. On the $58$ instances
that take the terminal box at least $10$\,s, the convex-combination constraints are faster on $36$,
by a median factor of $1.09$. Over the $239$ instances that the terminal box solves, they lower the
total solve time by $19\%$, and the H-representation raises it by $5\%$. At $d=4$, on the $120$ instances of test set~C
with $d=4$, $p\le25$ and $k\le5$, the H-representation reduces the median node count by a factor of $1.10$ and
the median solve time by $7\%$ (by $14\%$ on the instances that take at least $10$\,s), but it
solves one instance fewer than the terminal box ($117$ against $118$). We also added the H-representation as lazy constraints, with two settings of the solver's
\texttt{Lazy} attribute: with value~$1$ a constraint enters the model only when it cuts off a feasible
solution, and with value~$3$ also when it cuts off the relaxation solution at the root
node~\citep*{gurobi2025}. Both are worse than the terminal box. With values $1$ and $3$, the median
solve time grows by factors of $1.9$ and $1.8$ in the plane and of $3.1$ and $3.2$ at $d=3$, the
solver explores about four times as many nodes in the plane and $6.5$ times as many at $d=3$, and
two instances and one instance at $d=3$, respectively, are no longer solved. We keep the terminal box~$B$, whose constraints are the simplest of the three, in our final experiments discussed in Section \ref{sec:computational}.}

\subsection{Constructive upper bound for the optimal value}
\label{sec:UH}

We derive a simple constructive upper bound that produces both a feasible solution and a certificate for the optimal value. Let $T_{\MST}$ denote the Euclidean minimum spanning tree on the terminal set $P$. As the EBSTP with $k=0$ reduces to the Euclidean bottleneck spanning tree problem, and every Euclidean MST is a bottleneck spanning tree, we have
$
L^{\star}_{\le 0}=\beta(T_{\MST}).
$
Proposition~\ref{prop:monotone} therefore implies
$
L^{\star}_{\le k}\le \beta(T_{\MST})
$
for every Steiner budget $k\ge 0$.

The following construction yields a sharper instance-specific upper bound.
Starting from $T_{\MST}$, we perform $k$ successive bisections. At each iteration, select a longest edge $\{u,v\}$ of the current tree, introduce a Steiner point at its midpoint
$
q=(u+v)/2,
$
and replace $\{u,v\}$ by the two edges $\{u,q\}$ and $\{q,v\}$. After $k$ iterations, the resulting tree $T_H$ spans all terminals and contains exactly $k$ Steiner points, each of degree~2. We define
$
U_H:=\beta(T_H).
$

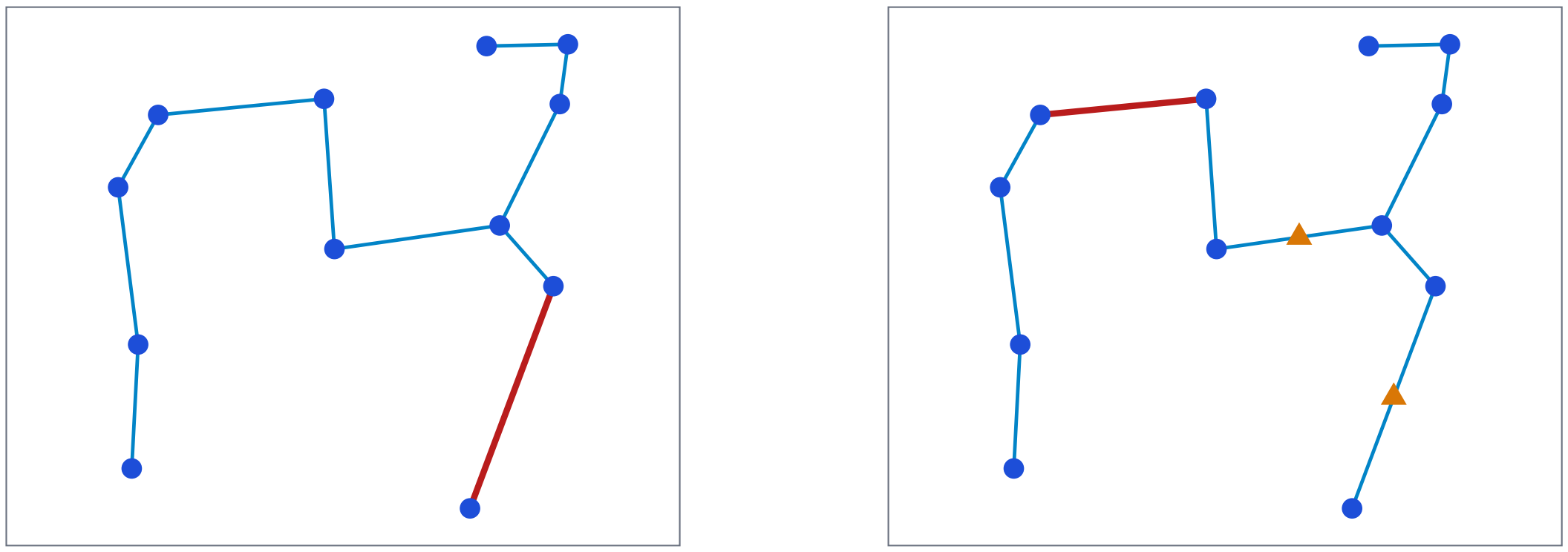
\begin{figure}[htbp]
  \centering
  \begin{minipage}[t]{0.48\linewidth}\centering
\begin{tikzpicture}
\begin{axis}
  [
    steinergeo,
    width = 6.0000cm,
    height = 4.8000cm,
    xmin = -0.1204, xmax = 1.2015,
    ymin = -0.0017, ymax = 1.0559,
  ]
  \addplot[btedge, no marks, forget plot] coordinates {(0.5033,0.8758) (0.1777,0.8442)};
  \addplot[btedge, no marks, forget plot] coordinates {(0.9535,0.5078) (0.8482,0.6270)};
  \addplot[btbottleneck, no marks, forget plot] coordinates {(0.9535,0.5078) (0.7898,0.0713)};
  \addplot[btedge, no marks, forget plot] coordinates {(0.0991,0.7020) (0.1384,0.3932)};
  \addplot[btedge, no marks, forget plot] coordinates {(0.0991,0.7020) (0.1777,0.8442)};
  \addplot[btedge, no marks, forget plot] coordinates {(0.8223,0.9794) (0.9820,0.9829)};
  \addplot[btedge, no marks, forget plot] coordinates {(0.9820,0.9829) (0.9661,0.8655)};
  \addplot[btedge, no marks, forget plot] coordinates {(0.5033,0.8758) (0.5237,0.5808)};
  \addplot[btedge, no marks, forget plot] coordinates {(0.1384,0.3932) (0.1258,0.1497)};
  \addplot[btedge, no marks, forget plot] coordinates {(0.5237,0.5808) (0.8482,0.6270)};
  \addplot[btedge, no marks, forget plot] coordinates {(0.8482,0.6270) (0.9661,0.8655)};
  \addplot[only marks, mark = *, mark size = 2.4pt, mark options = {fill = btterminal, draw = btterminal}, forget plot] coordinates {(0.0991,0.7020) (0.5033,0.8758) (0.9535,0.5078) (0.7898,0.0713) (0.1384,0.3932) (0.1777,0.8442) (0.8223,0.9794) (0.9820,0.9829) (0.5237,0.5808) (0.1258,0.1497) (0.8482,0.6270) (0.9661,0.8655) };
\end{axis}
\end{tikzpicture}\\[2pt]
    {\small (a) terminal MST, $\beta(T_{\MST})=0.4662$}
  \end{minipage}\hfill
  \begin{minipage}[t]{0.48\linewidth}\centering
\begin{tikzpicture}
\begin{axis}
  [
    steinergeo,
    width = 6.0000cm,
    height = 4.8000cm,
    xmin = -0.1204, xmax = 1.2015,
    ymin = -0.0017, ymax = 1.0559,
  ]
  \addplot[btbottleneck, no marks, forget plot] coordinates {(0.5033,0.8758) (0.1777,0.8442)};
  \addplot[btedge, no marks, forget plot] coordinates {(0.9535,0.5078) (0.8482,0.6270)};
  \addplot[btedge, no marks, forget plot] coordinates {(0.0991,0.7020) (0.1384,0.3932)};
  \addplot[btedge, no marks, forget plot] coordinates {(0.0991,0.7020) (0.1777,0.8442)};
  \addplot[btedge, no marks, forget plot] coordinates {(0.8223,0.9794) (0.9820,0.9829)};
  \addplot[btedge, no marks, forget plot] coordinates {(0.9820,0.9829) (0.9661,0.8655)};
  \addplot[btedge, no marks, forget plot] coordinates {(0.5033,0.8758) (0.5237,0.5808)};
  \addplot[btedge, no marks, forget plot] coordinates {(0.1384,0.3932) (0.1258,0.1497)};
  \addplot[btedge, no marks, forget plot] coordinates {(0.8482,0.6270) (0.9661,0.8655)};
  \addplot[btedge, no marks, forget plot] coordinates {(0.9535,0.5078) (0.8716,0.2895)};
  \addplot[btedge, no marks, forget plot] coordinates {(0.7898,0.0713) (0.8716,0.2895)};
  \addplot[btedge, no marks, forget plot] coordinates {(0.5237,0.5808) (0.6860,0.6039)};
  \addplot[btedge, no marks, forget plot] coordinates {(0.8482,0.6270) (0.6860,0.6039)};
  \addplot[only marks, mark = *, mark size = 2.4pt, mark options = {fill = btterminal, draw = btterminal}, forget plot] coordinates {(0.0991,0.7020) (0.5033,0.8758) (0.9535,0.5078) (0.7898,0.0713) (0.1384,0.3932) (0.1777,0.8442) (0.8223,0.9794) (0.9820,0.9829) (0.5237,0.5808) (0.1258,0.1497) (0.8482,0.6270) (0.9661,0.8655) };
  \addplot[only marks, mark = triangle*, mark size = 3.4pt, mark options = {fill = btsteiner, draw = btsteiner}, forget plot] coordinates {(0.8716,0.2895) (0.6860,0.6039) };
\end{axis}
\end{tikzpicture}\\[2pt]
    {\small (b) $T_H$ after $k=2$ bisections, $U_H=0.3271$}
  \end{minipage}
  \caption{Constructive upper bound obtained by iterative edge bisection on the terminal set of Figure~\ref{fig:gallery}(a) with Steiner budget $k=2$. Each bisection halves a longest
    edge by inserting a degree-two Steiner point at its midpoint; after the
    budget is spent, the bottleneck value of the resulting tree $T_H$ is the
    certified upper bound $U_H$. Filled blue circles represent terminals, orange
    triangles represent inserted Steiner points, and the thick red
  edges indicate the maximum-length (bottleneck) edges of the tree.}
  \label{fig:bisection}
\end{figure}

Figure~\ref{fig:bisection} illustrates the procedure on the instance of Figure~\ref{fig:gallery}(a).

\begin{lemma}\label{lem:bisection}
  For any feasible tree $T$ and any edge $\{u,v\}\in E_T$, replacing $\{u,v\}$ by the two-edge path $u\to x\to v$ with $x=(u+v)/2$ yields a feasible tree $T'$ in which the inserted midpoint $x$ has degree~$2$ and $\beta(T')\le\beta(T)$.
\end{lemma}

\begin{proof}
  This is the case $\alpha=1$ of the construction in the proof of Lemma~\ref{lem:padding}: the two new edges $[u,x]$ and $[x,v]$ each have length $\tfrac12\norm{u-v}_2\le\norm{u-v}_2$, so the bottleneck value does not increase, and the inserted midpoint $x$ has exactly two incident edges.
\end{proof}

\begin{proposition}\label{prop:UH}
  $L^{\star}_{\le k}\le U_H\le \beta(T_{\MST})$.
\end{proposition}

\begin{proof}
  $T_H$ is feasible for~\eqref{eq:ebstp}, so $L^{\star}_{\le k}\le U_H$. Each step applies Lemma~\ref{lem:bisection} to a longest edge, so the bottleneck value weakly decreases from $\beta(T_{\MST})$ to $U_H$.
\end{proof}

We use the pair $(U_H,T_H)$ for three purposes throughout the paper. First, it provides a partial MIP start for the MISOCP formulations: the epigraph variable is initialized at $L=U_H$, and the Steiner coordinates of $T_H$ are supplied as starting values, while the solver completes the edge-selection variables to obtain a full incumbent solution. Second, the upper bound $U_H$ is used in Section~\ref{sec:elim} to eliminate candidate edges from $E^\star$ that cannot belong to any optimal tree. Third, the upper bound provides the normalization scale used in the root-relaxation analysis of Section~\ref{sec:rootexp}.

{ 
On the $734$ instances of test set~C solved to optimality by {\rm \Idenone} or {\rm \Identity},
$U_H$ is optimal within $10^{-4}$ on $15.0\%$ of the instances, and its median relative
difference to the optimal value, measured by $(U_H-L^\star_{\le k})/U_H$, is $11.8\%$.
The median relative difference
grows with the budget, from $4.3\%$ at $k=1$ to $13.5\%$ at $k=5$ and $16.9\%$ at $k=9$, because
bisection only creates Steiner points of degree two on edges of the terminal MST.

A better upper bound would nevertheless help little, because most of the root gap lies on the
dual side: on the $360$ instances of test set~C
with $p\le25$ and $k\le5$, the primal part accounts for a
median of $28\%$ of the root gap~\eqref{eq:rootgap}. Two experiments
with a perfect upper bound show the most that a better heuristic could gain. In the first, we
give the solver the optimal value as an objective cutoff. This makes {\rm \Idenone} with the MST
restriction of Section~\ref{sec:elim} faster by a median factor of $1.46$ on test set~B and
solves the one instance that the partial start leaves at the time limit, but on the other $239$
instances the total solve time is $3\%$ higher. In the second, we give the solver the optimal
solution itself as a complete MIP start, with its coordinates, edge variables and $L$. This is
faster than the partial start from $T_H$ by a median factor of $1.62$, and of $1.51$ on the $56$
instances that take at least $10$\,s, with $1.6$ times fewer nodes; nevertheless, the total solve
time on test set~B grows from $4745$ to $4883$\,s, and one instance solved from the partial start
reaches the time limit when started from the optimum. A perfect upper bound thus makes the
median instance about $1.5$ times faster but does not reduce the total solve time.

We also tested a stronger heuristic. Starting from $T_H$, a local search repeatedly replaces one
Steiner point by the center of the minimum enclosing ball of two or three current points, as long
as the bottleneck value of the Euclidean MST of all points decreases. It takes a median of $0.05$\,s per
instance. On the $358$ of these $734$ instances with $p\le25$ and $k\le5$, it finds the optimal
value on $43.0\%$, against $15.1\%$ for $U_H$ on the same instances, and lowers the median relative
difference from $12.6\%$ to $2.7\%$. Used for the edge elimination and as the MIP start, however, it
leaves the median node count and the median solve time on test set~B unchanged and lowers the total
solve time by only $1\%$, so the formulations keep the simpler
bisection heuristic in our final experiments discussed in Section \ref{sec:computational}.
}

\FloatBarrier

\subsection{Edge elimination}
\label{sec:elim}

Edges with both endpoints in $P$ whose length exceeds $U_H$ cannot belong to any optimal tree, as selecting such an edge would induce a bottleneck value strictly larger than the value $U_H$ certified by~$T_H$. Formally, define
\begin{equation}
  E^\star := E \setminus \bigl\{[i,j]\in E\,:\,
  i,j\in P\;\text{and}\;\norm{a_i-a_j}_2 > U_H + \varepsilon\bigr\},
  \label{eq:elim}
\end{equation}
where $\varepsilon$ denotes a small numerical tolerance. Edges incident to Steiner points (i.e., Steiner-indexed vertices in the formulation) are retained because their lengths depend on the optimization decisions. Proposition~\ref{prop:UH} and the optimality of $L^{\star}_{\le k}$ yield:

\begin{corollary}\label{cor:elim}
  Every optimal tree $T^\star$ for~\eqref{eq:ebstp} satisfies $E_{T^\star}\subseteq E^\star$.
\end{corollary}

\begin{proof}
  Let $T^\star$ be any optimal tree for~\eqref{eq:ebstp}. Every terminal-terminal edge $[i,j]\in E_{T^\star}$ satisfies {$\norm{a_i-~a_j}_2 \le \beta(T^\star) = L^{\star}_{\le k} \le U_H$} by Proposition~\ref{prop:UH}, hence $[i,j]\in E^\star$. Edges of $T^\star$ incident to Steiner indices are retained in~$E^\star$ by definition.
\end{proof}

{ 
A much stronger restriction of the terminal--terminal edges holds in every dimension. In the
plane, some optimal tree has all its terminal--terminal edges in a minimum spanning tree of the
terminals: \citet*[Lemma~4]{bae2010exact} prove this for the Euclidean metric and
\citet*[Theorem~1, condition BST1]{bae2011exact} for every $L_p$ metric. 
The exchange argument in the proof of Lemma~2 of \citet*{bae2010exact} uses
only the cycle property of minimum spanning trees and therefore applies in $\mathbb{R}^d$; the
next proposition states the result in this generality.

\begin{proposition}\label{prop:mst}
  Let $H$ be any minimum spanning tree of the terminals. For every optimal tree $T$
for~\eqref{eq:ebstp} on the indexed node set $V$, there exists an optimal tree $T'$ with the
same Steiner points and the same Steiner-incident edges as $T$ in which every
terminal--terminal edge belongs to $H$.
\end{proposition}

\begin{proof}
  Let $\nu(T)$ be the number of terminal--terminal edges of $T$ outside $H$, and suppose
  $\nu(T)>0$, with $[u,v]$ such an edge. Removing $[u,v]$ splits $T$ into two subtrees containing
  $u$ and $v$, respectively. The path of $H$ between $u$ and $v$ consists of terminals, so it has
  an edge $[a,b]$ whose endpoints lie in different subtrees. As $[u,v]\notin H$, it is a longest
  edge of the unique cycle of $H+[u,v]$, so $\norm{a-b}_2\le\norm{u-v}_2$. As $T-[u,v]$ has no edge
  between the two subtrees, $[a,b]\notin E_T$, and $T':=T-[u,v]+[a,b]$ is a tree on $V$ with the
  same Steiner points and Steiner-incident edges as $T$ and $\beta(T')\le\beta(T)$. As
  $[a,b]\in H$, $\nu(T')=\nu(T)-1$, and at most $p-1$ exchanges give the result.
\end{proof}

As the Steiner points are unchanged and every inserted edge has length at most
$L^\star_{\le k}\le U_H$, Property~\ref{propP} holds for the candidate edge set
\begin{equation}
  E^\star_H := \bigl\{[i,j]\in E^\star : \{i,j\}\cap S\neq\emptyset \text{ or } [i,j]\in H\bigr\},\quad \mbox{(the MST restriction)}
  \label{eq:elimmst}
\end{equation}
which keeps at most $p-1$ terminal--terminal edges in every dimension (an edge of $H$ longer than
$U_H$ is already removed by~\eqref{eq:elim}). More generally, every optimal tree on $V$ without
Steiner leaves, with all Steiner points in $\conv(P)$ and with all terminal--terminal edges in
$H$ satisfies Property~\ref{propP} for $E^\star_H$: its terminal--terminal edges have length at
most $L^\star_{\le k}\le U_H$ and therefore survive rule~\eqref{eq:elim}, and the
constants~\eqref{eq:bigMdef} are valid for Steiner points in $\conv(P)$
(Proposition~\ref{prop:bigM}). Section~\ref{sec:implementation} states which
experiments use $E^\star_H$ and which use $E^\star$.

The same exchange, run in the opposite direction, fixes the short edges of $H$.

\begin{proposition}\label{prop:fixone}
  Let $\tilde{L}\le L^\star_{\le k}$ and $F:=\{[i,j]\in H:\norm{a_i-a_j}_2<\tilde{L}\}$.
  There exists an optimal tree satisfying Property~\ref{propP} for $E^\star_H$ that contains every
  edge of $F$. Hence $y_{ij}=1$ may be imposed for every $[i,j]\in F$.
\end{proposition}

\begin{proof}
  Call a tree on $V$ \emph{admissible} if it is optimal, has no Steiner leaves, has all its
  Steiner points in $\conv(P)$, and has all its terminal--terminal edges in $H$. Admissible trees
  satisfy Property~\ref{propP} for $E^\star_H$, as noted above, and one exists by
  Propositions~\ref{prop:indexed}, \ref{prop:box} and~\ref{prop:mst}. For an admissible tree $T$,
  let $\varphi(T):=|F\setminus E_T|$. We show that if $\varphi(T)>0$, there is an admissible
  tree $T'$ with $E_{T'}\cap F\supsetneq E_T\cap F$; as $\varphi$ then decreases by at least one,
  at most $|F|$ repetitions give an admissible tree containing every edge of $F$.

  Let $f=[u,v]\in F\setminus E_T$. Every edge $g$ on the path of $T$ between $u$ and $v$ is
  either incident to a Steiner point or a terminal--terminal edge of $H$; not all of them can be
  of the second kind, since together with $f\in H$ they would form a cycle in $H$. Choose such a
  path edge $g\notin F$; in particular $g$ is not a terminal--terminal edge of $F$.
  Then $T_1:=T+f-g$ is a spanning tree of $V$ whose terminal--terminal edges still lie in $H$,
  and $\beta(T_1)\le\max\{\beta(T),\norm{a_u-a_v}_2\}=L^\star_{\le k}$, because
  $\norm{a_u-a_v}_2<\tilde{L}\le L^\star_{\le k}$.

  Removing $g$ may have created Steiner leaves. Delete Steiner leaves repeatedly; each deletion
  removes one Steiner point and one Steiner-incident edge, so no edge of $F$ is lost, and the
  result is a tree $T_2$ with $k'\le k$ Steiner points, none of them a leaf, and
  $\beta(T_2)\le L^\star_{\le k}$. To restore exactly $k$ Steiner points, choose an edge $g'$
  of $T_2$ not in $F$: if $k'\ge1$, any Steiner-incident edge; if $k'=0$, then $T_2$ is a
  spanning tree of the terminals with all edges in $H$, so $T_2=H$, and we take a longest edge
  of $H$, which is not in $F$ because its length is $\beta(T_{\MST})\ge L^\star_{\le k}\ge
  \tilde{L}$ (Proposition~\ref{prop:UH}). Subdividing $g'$ by $k-k'$ equally spaced points, as
  in the proof of Lemma~\ref{lem:padding}, inserts degree-two Steiner points on the segment
  $g'$, hence in $\conv(P)$, without increasing the bottleneck. The resulting tree $T'$ is
  admissible, contains every edge of $F$ that $T$ contains, and contains $f$.
\end{proof}

The strict inequality in the definition of $F$ is needed. On the dodecahedron with $k=1$, every
MST edge has length $L^\star_{\le 1}$ (Section~\ref{sec:namedexp}), and fixing all of them to one
would leave a single edge for the Steiner point, which needs at least two. The fixings have
little effect on the relaxation: with $y_{ab}=1$ on an edge of $F$, constraint~\eqref{eq:soc}
gives $L\ge\norm{a-b}_2$, so the bare relaxation value becomes at least the length of the
longest fixed edge, which is below $\tilde{L}$ by the definition of $F$; by
Proposition~\ref{prop:rootLB}, the root bound with the global lower bounds is therefore
unchanged. Consistently, on test set~B, imposing the fixings changes no optimal value found by
{\rm \Idenone} and reduces its number of nodes by a median factor of only $1.08$; it solves the one
instance that {\rm \Idenone} otherwise leaves at the time limit, but raises its total solve time on the
other $239$ instances by $13\%$. For {\rm \Identity} the median node ratio is $1.00$, and five instances
that it solves without the fixings reach the time limit with them. The experiments of
Section~\ref{sec:computational} therefore do not impose the fixings.
}

\subsection{Edge-type-specific \texorpdfstring{big-$M$}{big-M} constants}
\label{sec:bigM}

A uniform big-$M$ choice such as $M=\diam(B)$ is valid but overly conservative for every edge type. The convex-hull containment result of Proposition~\ref{prop:box} allows us to restrict attention to distances within $\conv(P)$, yielding  tighter bounds that can be specialized to each edge type.

For terminal-terminal edges $[i,j]$ with $i,j\in P$, both endpoints are fixed, and we set $M_{ij} := \norm{a_i - a_j}_2$. For terminal-Steiner edges $[i,j]$ with $i\in P$, $j\in S$, the Steiner endpoint can be assumed to lie in $\conv(P)$, as at least one optimal solution satisfies this assumption (Proposition~\ref{prop:box}). We therefore use the distance from $a_i$ to the farthest point of $\conv(P)$,
\begin{equation}
  \textstyle \Delta_i \;:=\; \max_{t\in P}\norm{a_i - a_t}_2.
  \label{eq:DeltaiB}
\end{equation}
As the function $x\mapsto\norm{a_i-x}_2$ is convex, its maximum over the polytope $\conv(P)$ is attained at an extreme point, and every extreme point of $\conv(P)$ is a terminal. Hence, $\Delta_i = \max_{x\in\conv(P)}\norm{a_i-x}_2$. For Steiner-Steiner edges $[i,j]$ with $i,j\in S$, both endpoints may be assumed to lie in $\conv(P)$, as at least one optimal solution satisfies this assumption (Proposition~\ref{prop:box}), and we use $M_{ij} := \diam(P) = \max_{i',j'\in P}\norm{a_{i'}-a_{j'}}_2$, which equals the diameter of $\conv(P)$, i.e., the maximum distance between any two points in $\conv(P)$.

{ 
\citet*{fampa2004} use the uniform constant $M=\diam(P)$ for all edges. Constants of this form were proposed for the classical Euclidean Steiner tree problem in \citep*[Proposition 2]{fampa2016}, where their validity follows from the fact that every Steiner minimal tree has its Steiner points in $\conv(P)$. In the bottleneck setting this property does not hold for every optimal solution (Section~\ref{sec:contrast}), and the validity of the constants requires Proposition~\ref{prop:box}.
}

\begin{proposition}\label{prop:bigM}
  The constants
  \begin{equation}
    M_{ij} =
    \begin{cases}
      \norm{a_i-a_j}_2, & i,j\in P,\\[2pt]
      \Delta_i,         & i\in P,\,j\in S,\\[2pt]
      \diam(P),         & i,j\in S,
    \end{cases}
    \label{eq:bigMdef}
  \end{equation}
  preserve the optimal value $L^{\star}_{\le k}${ , and the bound given by the
  continuous relaxation with these constants is never smaller than the bound given by the
  uniform big-$M$ choice $M=\diam(B)$.}
\end{proposition}

\begin{proof}
  { We show that the constants~\eqref{eq:bigMdef} preserve the optimal value, and that
  the relaxation bound they give is never smaller than with $\diam(B)$.}

  \emph{(i) Tightening does not decrease the optimal value.}
  For any edge $[i,j]$, when $y_{ij}=0$ constraint~\eqref{eq:soc} reduces to $
  \|x_i-x_j\|_2 \le L + M_{ij}$.
  Hence, decreasing any $M_{ij}$ restricts the feasible set of the relaxation and can only (weakly) increase the optimal value. In particular, replacing $\diam(B)$ by the smaller constants in \eqref{eq:bigMdef} yields an objective value no smaller than that of the terminal box formulation.

  \emph{(ii) Feasibility at an optimal solution.}
  By Proposition~\ref{prop:box}, there exists an optimal solution in which all Steiner points lie in $\conv(P)$. Consider such a solution. If both endpoints are terminals, then $
  \|a_i-a_j\|_2 \le M_{ij}$
  holds by definition. If $i\in P$ and $j\in S$, then $x_j\in\conv(P)$, and by definition of $\Delta_i$,
  $
  \|a_i-x_j\|_2 \le \max_{x\in\conv(P)} \|a_i-x\|_2 = \Delta_i$.
  If $i,j\in S$, then $x_i,x_j\in\conv(P)$, implying $
  \|x_i-x_j\|_2 \le \diam(\conv(P)) = \diam(P)$.
  Thus all constraints remain satisfied at this optimal solution, so the optimal value is preserved.

  \emph{(iii) { Comparison of the relaxation bounds.}}
  As $\conv(P)\subseteq B$, we have
  $
  \Delta_i \le \max_{x\in B}\|a_i-x\|_2$, and $\diam(P) \le \diam(B)$,
  with strict inequality in general. Therefore, { by the monotonicity noted in (i), the
  feasible region of the continuous relaxation induced by~\eqref{eq:bigMdef} is contained in
  that obtained using the uniform bound $\diam(B)$, so the relaxation bound with the
  constants~\eqref{eq:bigMdef} is never smaller than with $\diam(B)$.}
\end{proof}

\subsection{Symmetry breaking on Steiner labels}
\label{sec:symm}

{ 
The indexed representation of Section~\ref{sec:indexed} introduces a permutation symmetry
$\mathfrak{S}_k$ on Steiner labels. This symmetry acts on continuous variables, and the automatic
symmetry detection of the solver used in our numerical experiments does not find it, even at its most aggressive setting
(Section~\ref{sec:symmexp}). We eliminate this symmetry by ordering the Steiner points along the first
coordinate,
\begin{equation}
\textstyle   (x_{p+j})_1 \;\le\; (x_{p+j+1})_1, \qquad j=1,\dots,k-1.
  \label{eq:s1}
\end{equation}
These $k-1$ linear inequalities keep at least one representative of every permutation orbit,
which is all that correctness requires. They do not eliminate the symmetry between Steiner points
that share their first coordinate. A lexicographic order on two coordinates would require binary
variables forcing equal first coordinates, and its big-$M$ relaxation, in which these binaries
may relax the second-coordinate constraint, has the same projection onto the Steiner coordinates
as~\eqref{eq:s1}. We refer to~\eqref{eq:s1} as the ordering S1.
}

\subsection{Global lower bounds on the optimal value}
\label{sec:diam}

The strengthening techniques developed in Sections~5.1--5.5 either tighten the integer feasible region or reduce the candidate edge set. None, however, addresses the principal weakness of the formulations: the continuous relaxation of the epigraph variable~$L$. At fractional edge-selection variables, the big-$M$ terms may effectively deactivate the cone constraints, allowing the trivial bound $L=0$. To overcome this limitation, we derive two valid global lower bounds on the optimal value that act directly on~$L$. The first is an elementary diameter bound obtained from the triangle inequality. The second is a new clustering bound based on the separation induced by removing the longest edges of the Euclidean minimum spanning tree. Besides strengthening the formulations, this result provides a new global lower bound for the EBSTP that also reveals a connection between the optimal value and the clustering structure of the terminal set. Together, the two inequalities replace the trivial relaxation by a strictly positive estimate while preserving the exact feasible region.

\begin{proposition}\label{prop:diam}
  $L^{\star}_{\le k}\ge\diam(P)/(p+k-1)$, where $\diam(P)=\max_{i,j\in P}\norm{a_i-a_j}_2$.
\end{proposition}

\begin{proof}
  In any feasible tree, every pair of terminals $a_i,a_j$ is joined by a simple path of at most $p+k-1$ edges. By the triangle inequality the sum of those edge lengths is at least $\norm{a_i-a_j}_2$, so the longest edge, and hence the bottleneck value, is at least $\norm{a_i-a_j}_2/(p+k-1)$. Maximizing over terminal pairs gives the bound.
\end{proof}

Accordingly, we add the valid inequality
\begin{equation}
\textstyle   L \;\ge\; \frac{\diam(P)}{p+k-1}.
  \label{eq:diamcut}
\end{equation}
The bound is computed in $O(p^2d)$ time during preprocessing. As
$\diam(P)>0$ whenever the terminal set contains at least two distinct
points, \eqref{eq:diamcut} eliminates the trivial relaxation value
$L=0$, which may otherwise satisfy~\eqref{eq:soc} when the edge-selection
variables are fractional.

The denominator $p+k-1$ corresponds to the maximum possible number of
edges in a terminal-to-terminal path. A seemingly stronger inequality,
\[
 \textstyle  L \ge \frac{\|a_i-a_j\|_2}{k+1},
\]
is generally invalid, because the path between two terminals may traverse
other terminals in addition to the at most $k$ Steiner points. For
example, consider four collinear terminals at unit spacing and $k=0$.
An optimal tree, the path through consecutive terminals, has bottleneck value $L^\star=1$, whereas $\norm{a_1-a_4}_2/(k+1)=3$.

The second global lower bound complements~\eqref{eq:diamcut}, particularly on clustered instances where the diameter bound becomes weak due to the large denominator $p+k-1$.

{ 
\begin{theorem}\label{thm:clusterbound}
  Let $H$ be a minimum spanning tree of the terminals and let $e_1\ge e_2\ge\dots\ge e_{p-1}$
  denote the lengths of its edges, sorted in nonincreasing order, with $e_k:=0$ for $k\ge p$;
  these lengths do not depend on the choice of $H$. Then, for every $k\ge 1$,
  \begin{equation}
    \textstyle L^{\star}_{\le k}\;\ge\;\tfrac{1}{2}\,e_{k}.
    \label{eq:clustercut}
  \end{equation}
\end{theorem}

\begin{proof}
  If $k\ge p$, the claim is trivial. Hence assume $1\le k\le p-1$. Remove from $H$ a set $D$ of
  $k$ edges of lengths $e_1,\dots,e_k$ (when several edges have the same length, any such choice
  will do). As $H$ is a tree, this produces exactly $k+1\ge2$ connected components, whose
  terminal sets we denote by $C_1,\dots,C_{k+1}$. Define
  \[
    \textstyle \mathrm{sep}\;:=\;\min_{r\ne s}\;\min_{a\in C_r,\,b\in C_s}\norm{a-b}_2.
  \]
  We first show that $\mathrm{sep}=e_{k}$. A shortest edge of $D$ has length $e_k$, and its
  endpoints belong to two distinct components, so $\mathrm{sep}\le e_{k}$. For the reverse
  inequality, suppose that two terminals $a\in C_r$ and $b\in C_s$, with $r\neq s$, satisfy
  $\norm{a-b}_2<e_{k}$. By the minimax path property of minimum spanning trees, every edge on
  the path of $H$ between $a$ and $b$ has length at most $\norm{a-b}_2$
  (see, e.g., \citealp*{Ahuja2014}), hence strictly smaller than $e_k$. As every edge of $D$ has
  length at least $e_k$, no edge of that path belongs to $D$, so $a$ and $b$ remain in the same
  component, contradicting $r\neq s$. Therefore $\mathrm{sep}\ge e_{k}$, and thus
  $\mathrm{sep}=e_{k}$.

  Now suppose that a feasible Steiner tree $T$ has bottleneck value $\beta(T)<\mathrm{sep}/2$.
  For each cluster $C_r$, define the open neighborhood
  $
    B_r := \{x\in\mathbb R^d : \mathrm{dist}(x,C_r)<\mathrm{sep}/2\}$, where $\mathrm{dist}(x,C_r):=\min_{a\in C_r}\norm{a-x}_2 .
  $
  Because distinct clusters are at distance at least $\mathrm{sep}$ from each other, the sets
  $B_1,\dots,B_{k+1}$ are pairwise disjoint, and the only terminals in $B_r$ are those of $C_r$.

  Fix $r\in\{1,\dots,k+1\}$. As $T$ is connected and contains terminals both in $C_r$ and outside
  $C_r$, there is a path in $T$ from a terminal of $C_r$, which lies in $B_r$, to a terminal of
  another cluster, which does not. Hence some edge $\{u,v\}$ of $T$ has $u\in B_r$ and
  $v\notin B_r$. As
  $
    \norm{u-v}_2 \le \beta(T) < \mathrm{sep}/2,
  $ 
  $u$ cannot be a terminal: a terminal in $B_r$ belongs to $C_r$, and if $u\in C_r$ then
  $\mathrm{dist}(v,C_r)\le\norm{u-v}_2<\mathrm{sep}/2$, that is, $v\in B_r$, a contradiction.
  Hence $u$ is a Steiner point.

  Repeating the argument for every cluster yields Steiner points $u_1,\dots,u_{k+1}$ with
  $u_r\in B_r$. As the neighborhoods $B_1,\dots,B_{k+1}$ are pairwise disjoint, these Steiner
  points are distinct, so $T$ contains at least $k+1$ Steiner points, which contradicts the
  budget of at most $k$. Hence every feasible tree satisfies
  $\beta(T)\ge\mathrm{sep}/2=e_{k}/2$, and taking the minimum over all feasible trees gives the
  result.
\end{proof}

The factor $\tfrac12$ cannot be improved. For every $k\ge1$ and $d\ge2$, place $k+1$ terminals at
$a_i:=(i,0,\dots,0)$, $i=0,\dots,k$. The MST is the path through consecutive terminals, so
$e_k=1$, 
and one Steiner point at the midpoint of each edge of this path gives a tree whose $2k$ edges
all have length $\tfrac12$, so $L^\star_{\le k}=\tfrac12 e_k$. The factor $1$ is invalid already on the unit square
with $k=1$, where $e_1=1$ and a Steiner point at the center gives
$L^\star_{\le 1}\le 1/\sqrt{2}$. For $k=0$, the EBSTP is the bottleneck spanning tree problem and
$L^\star_{\le 0}=e_1$. In the plane, \citet*{bae2010exact} prove
$L^\star_{\le k}\ge e_{4k+1}$ when $p>4k+1$, from the planar degree bound of five; neither bound
dominates the other, and~\eqref{eq:clustercut} holds in every dimension.
}

By~\eqref{eq:clustercut}, we add the valid inequality $L \!\geq\! \frac{1}{2} e_{k}$ to
both formulations, {\rm \Idenone} and {\rm \Identity}. Like the diameter bound in~\eqref{eq:diamcut}, this bound on $L$ is a static constant computed in $O(p^2 d)$ time that leaves the optimal value unchanged. The two bounds are complementary: $\tfrac{1}{2}e_{k}$ is independent of $p$ and remains bounded away from zero on instances with well-separated clusters, whereas $\diam(P)/(p+k-1)$ is stronger on near-collinear configurations. We therefore impose both; equivalently, we impose $L\ge\max\{\tfrac{1}{2}e_{k},\,\diam(P)/(p+k-1)\}$. Section~\ref{sec:rootexp} measures their effect on the root~relaxation.


{ 
\subsection{The continuous relaxation of the epigraph formulation}
\label{sec:relaxation}
Challenge~\ref{c1} of Section~\ref{sec:summary} identifies the continuous relaxation of the
epigraph variable $L$ as the principal weakness of the formulations, and Section~\ref{sec:rootexp}
will show that, without the global lower bounds of Section~\ref{sec:diam}, the root bound is $0$
on every test instance with $k\ge3$ and positive on every instance with $k=1$. This section
explains both observations. 
We show that the relaxation always has an optimal solution in which all Steiner points coincide
(Proposition~\ref{prop:symrelax}); that this reduces the question of whether the relaxation
bound is $0$ to a condition on a single point, which is never satisfied for $k=1$ and can be
checked by a small second-order cone program (Theorem~\ref{thm:degenerate}); and that, when the
bound is $0$, the global lower bounds of Section~\ref{sec:diam} determine the root bound exactly
(Proposition~\ref{prop:rootLB}). The analysis also shows which of the strengthening techniques of this section can raise the
root bound.

Let $\mathrm{R}$ denote the continuous relaxation of {\rm \Idenone}
with the terminal box~\eqref{eq:x}, the candidate edge set
$E^\star$ of~\eqref{eq:elim} or $E^\star_H$ of~\eqref{eq:elimmst},
the constants~\eqref{eq:bigMdef}, and the ordering~\eqref{eq:s1}.
The root $r\in P$ and arc flows are as in Section~\ref{sec:m1}.
We replace $y_{ij}\in\{0,1\}$ by $0\le y_{ij}\le1$ and omit
the global lower bounds~\eqref{eq:diamcut} and~\eqref{eq:clustercut}.
Let $v_k$ denote the optimal value of this relaxation. Below we characterize the instances of~\eqref{eq:ebstp} for which $v_k=0$.

\begin{proposition}\label{prop:symrelax}
  Let $k\ge1$. The relaxation $\mathrm{R}$ has an optimal solution $(x,y,f,L)$ in which
  \begin{itemize}
    \item all Steiner points coincide, $x_{p+1}=\dots=x_{p+k}$;
    \vspace{-5pt}
    \item for every terminal $i\in P$, the value $y_{ij}$ is the same for all $j\in S$;
    \vspace{-5pt}
    \item the value $y_{ij}$ is the same for all Steiner--Steiner edges $[i,j]$, $i,j\in S$.
  \end{itemize}
  The optimal value of $\mathrm{R}$ is the same with and without the ordering~\eqref{eq:s1}.
  The same statements hold for the relaxation of {\rm \Identity} with all connectivity
  cuts~\eqref{eq:m2cut} and for the multi-commodity flow relaxation of Section~\ref{sec:m2}.
\end{proposition}

\begin{proof}
  Let $\mathrm{R}'$ be the relaxation $\mathrm{R}$ without~\eqref{eq:s1}. It is feasible, as
  {\rm \Idenone} is (Proposition~\ref{prop:m1}), and for every $L_0$ its feasible solutions with
  $L\le L_0$ form a compact set, so $\mathrm{R}'$ has an optimal solution.

  For a permutation $\pi$ of $S$, extend $\pi$ to $V$ by $\pi(i):=i$ for $i\in P$, and map a
  solution $(x,y,f,L)$ to $(x^\pi,y^\pi,f^\pi,L)$, where $x^\pi_{\pi(j)}:=x_j$ for $j\in S$,
  $y^\pi_{\{\pi(i),\pi(j)\}}:=y_{\{i,j\}}$ for every edge $\{i,j\}\in E^\star$, regarded here as
  an unordered pair, and $f^\pi_{\pi(i)\pi(j)}:=f_{ij}$ for $(i,j)\in A$. This map preserves
  feasibility for $\mathrm{R}'$ and the value of $L$. The candidate edge set $E^\star$ (or
  $E^\star_H$) contains every edge incident to a Steiner point, so $\pi$ maps $E^\star$ onto
  itself and $A$ onto itself. The constant $M_{ij}$ in~\eqref{eq:bigMdef} of an edge incident to
  a Steiner point depends only on its terminal endpoint, or equals $\diam(P)$ for a
  Steiner--Steiner edge, so the cone constraints~\eqref{eq:soc} are permuted among themselves.
  The bounds~\eqref{eq:x} are the same for every Steiner point, and
  constraints~\eqref{eq:edgecount}, \eqref{eq:noleaf} and~\eqref{eq:m1cap1}--\eqref{eq:m1cap}
  are permuted among themselves. Finally, the root $r$ is a terminal and is fixed by $\pi$, so
  the supply constraint~\eqref{eq:m1root} is unchanged, and the balance
  constraints~\eqref{eq:m1bal}, which have the same right-hand side at every node of
  $V\setminus\{r\}$, are permuted among themselves.

  Let $(x,y,f,L)$ be an optimal solution of $\mathrm{R}'$ and let $(\bar x,\bar y,\bar f,L)$ be
  the average of its $k!$ images under the permutations of $S$. Every constraint of
  $\mathrm{R}'$ is linear or a second-order cone constraint, so its feasible set is convex and
  the average is feasible; it is optimal because its value of $L$ is unchanged. For $j\in S$, the Steiner point with label $j$ in the image under $\pi$ is $x^\pi_j=x_{\pi^{-1}(j)}$, and each label
  $j'\in S$ equals $\pi^{-1}(j)$ for exactly $(k-1)!$ of the $k!$ permutations $\pi$, so
  $\bar x_j=\frac{1}{k!}\sum_{\pi}x_{\pi^{-1}(j)}=\frac1k\sum_{j'\in S}x_{j'}$ for every
  $j\in S$. The same counting gives $\bar y_{ij}=\frac1k\sum_{j'\in S}y_{ij'}$ for every
  $i\in P$ and $j\in S$. When $k=1$ there are no Steiner--Steiner edges, so the third property
  is vacuous. When $k\ge2$, for a Steiner--Steiner edge $\{j_1,j_2\}$, the image solution under
  $\pi$ carries the value $y_{\{\pi^{-1}(j_1),\pi^{-1}(j_2)\}}$ on it, and for every unordered
  pair $\{j',j''\}$ of Steiner labels there are exactly $2(k-2)!$ permutations $\pi$ with
  $\{\pi^{-1}(j_1),\pi^{-1}(j_2)\}=\{j',j''\}$. Hence
  \[
     \bar y_{j_1j_2}
    \;=\;\frac{2(k-2)!}{k!}\sum_{\{j',j''\}\subseteq S} y_{j'j''}
    \;=\;\left(1/\textstyle\binom{k}{2}\right)\sum_{\{j',j''\}\subseteq S} y_{j'j''},
  \]
  the mean of $y$ over all Steiner--Steiner edges, which is the same for every edge
  $\{j_1,j_2\}$. Hence $(\bar x,\bar y,\bar f,L)$ has the stated structure. As all Steiner
  points coincide, it satisfies~\eqref{eq:s1} with equality and is feasible for $\mathrm{R}$.
  The feasible set of $\mathrm{R}$ is contained in that of $\mathrm{R}'$, so the optimal value
  of $\mathrm{R}$ is at least that of $\mathrm{R}'$, and $(\bar x,\bar y,\bar f,L)$ attains it.
  The two optimal values therefore coincide, and $(\bar x,\bar y,\bar f,L)$ is optimal
  for~$\mathrm{R}$.

  The argument is unchanged for the relaxation of {\rm \Identity} with all
  cuts~\eqref{eq:m2cut}, as $\pi$ maps the family of sets $C\subseteq V\setminus\{r\}$ onto
  itself, and for the multi-commodity flow relaxation, with the commodity of node $j$ mapped
  to the commodity of node $\pi(j)$.
\end{proof}

By Proposition~\ref{prop:symrelax}, it suffices to search for an optimal solution of
$\mathrm{R}$ in which all Steiner points coincide at a single point $c$. 
This reduction yields an exact characterization, in terms of $c$, of the instances for which $v_k=0$, meaning that the relaxation provides only the trivial lower bound of $0$ on $L^\star_{\le k}$. We refer to such instances as \emph{degenerate} (for the relaxation $\mathrm{R}$; when another
relaxation is meant, we say so).

\begin{theorem}\label{thm:degenerate}
  Let $k\ge1$ and $p\ge k+1$. For $i\in P$ and $c\in\mathbb{R}^d$, let
  $\rho_i(c):=\norm{a_i-c}_2/\Delta_i$, with $\Delta_i$ defined in~\eqref{eq:DeltaiB}. Then
  $v_k=0$ if and only if there exists $c\in B$ such that
  \begin{align}
    & \rho_r(c)\le 1-\tfrac{1}{k}, \qquad
      \rho_i(c)\le 1-\tfrac{1}{k(p+k-1)}, \quad i\in P\setminus\{r\},
      \label{eq:degen1}\\
    & k\sum_{i\in P}\bigl(1-\rho_i(c)\bigr)+\tfrac{k(k-1)}{2}\;\ge\;p+k-1.
      \label{eq:degen2}
  \end{align}
  In particular, $v_1>0$ for every instance.
\end{theorem}

\begin{proof}
  Throughout, we use the structure of Proposition~\ref{prop:symrelax}: a solution of
  $\mathrm{R}$ in which $x_j=c$ for all $j\in S$, $y_{ij}=t_i$ for all $i\in P$ and $j\in S$,
  and $y_{ij}=\sigma$ for all $i,j\in S$. For such a solution, constraint~\eqref{eq:edgecount}
  reads
  \begin{equation}
    \sum_{[i,j]\in E^\star:\,i,j\in P} y_{ij}
    \;+\;k\sum_{i\in P}t_i\;+\;\sigma\,\tfrac{k(k-1)}{2}\;=\;p+k-1.
    \label{eq:degencard}
  \end{equation}

  ($\Rightarrow$) Let $v_k=0$. By Proposition~\ref{prop:symrelax}, $\mathrm{R}$ has an optimal
  solution with the structure above and $L=0$, and $c\in B$ by~\eqref{eq:x}. For a
  terminal--terminal edge $[i,j]\in E^\star$, constraint~\eqref{eq:soc} reads
  $\norm{a_i-a_j}_2\le\norm{a_i-a_j}_2(1-y_{ij})$; as the terminals are distinct, $y_{ij}=0$. For
  a terminal--Steiner edge it reads $\norm{a_i-c}_2\le\Delta_i(1-t_i)$, that is,
  $t_i\le1-\rho_i(c)$. As every terminal--terminal arc has capacity zero
  by~\eqref{eq:m1cap1}--\eqref{eq:m1cap}, the root sends its $p+k-1$ units of
  supply~\eqref{eq:m1root} over its $k$ arcs to Steiner points, each of capacity $(p+k-1)t_r$,
  so $t_r\ge1/k$; and every terminal $i\ne r$ receives its unit of demand~\eqref{eq:m1bal} over
  its $k$ arcs from Steiner points, each of capacity $(p+k-1)t_i$, so $t_i\ge1/(k(p+k-1))$.
  Together with $t_i\le1-\rho_i(c)$, this gives~\eqref{eq:degen1}.
  Finally, \eqref{eq:degencard} with $y_{ij}=0$ on terminal--terminal edges, $t_i\le1-\rho_i(c)$
  and $\sigma\le1$ gives~\eqref{eq:degen2}.

  ($\Leftarrow$) Let $c\in B$ satisfy~\eqref{eq:degen1}--\eqref{eq:degen2}. Set $L:=0$,
  $x_j:=c$ for $j\in S$, and $y_{ij}:=0$ on terminal--terminal edges. Let $\tau_r:=1/k$ and
  $\tau_i:=1/(k(p+k-1))$ for $i\ne r$. By~\eqref{eq:degen1}, the intervals
  $[\tau_i,1-\rho_i(c)]$, $i\in P$, are nonempty. As $(t,\sigma)$ ranges over
  $\prod_{i\in P}[\tau_i,1-\rho_i(c)]\times[0,1]$, the left-hand side of~\eqref{eq:degencard}
  is continuous and ranges from $k\sum_i\tau_i=1+(p-1)/(p+k-1)<2\le p+k-1$ to the left-hand
  side of~\eqref{eq:degen2}, which is at least $p+k-1$; hence some $(t,\sigma)$ in this set
  satisfies~\eqref{eq:degencard}. Fix it and set $y$ accordingly; then $0\le y\le1$. The cone
  constraints~\eqref{eq:soc} hold: on terminal--terminal edges because $y_{ij}=0$, on
  terminal--Steiner edges because $t_i\le1-\rho_i(c)$, and on Steiner--Steiner edges because
  both endpoints equal $c$. Define the flow by $f_{rj}:=(p+k-1)/k$ for $j\in S$,
  $f_{ji}:=1/k$ for $j\in S$ and $i\in P\setminus\{r\}$, and $f:=0$ on all other arcs. The root
  then sends $p+k-1$ units, each Steiner point receives $(p+k-1)/k$ units and sends $(p-1)/k$,
  for a net inflow of $1$, and each terminal $i\ne r$ receives $k\cdot\tfrac1k=1$ unit, so
  \eqref{eq:m1root}--\eqref{eq:m1bal} hold; the capacities~\eqref{eq:m1cap1}--\eqref{eq:m1cap}
  hold because $t_r\ge1/k$ and $t_i\ge1/(k(p+k-1))$. By~\eqref{eq:degencard}, the degree of
  every Steiner point is
  $\sum_{i\in P}t_i+\sigma(k-1)\ge\sum_{i\in P}t_i+\sigma(k-1)/2=(p+k-1)/k\ge2$, as $p\ge k+1$,
  so~\eqref{eq:noleaf} holds. Constraint~\eqref{eq:x} holds as $c\in B$, and~\eqref{eq:s1} holds
  with equality. Hence $(x,y,f,0)$ is feasible for $\mathrm{R}$ and $v_k=0$.

  For $k=1$, \eqref{eq:degen1} gives $\rho_r(c)\le0$, so $c=a_r$, and~\eqref{eq:degen2} becomes
  $\sum_{i\in P}(1-\rho_i(a_r))\ge p$, which requires $\rho_i(a_r)=0$, that is $a_i=a_r$, for
  every $i\in P$; this is impossible as $p\ge2$ and the terminals are distinct. Hence $v_1>0$.
\end{proof}

\begin{remark}
  The objective $\sum_i(1-\rho_i(c))$ is concave and~\eqref{eq:degen1} defines a convex set, so
  whether an instance is degenerate can be decided by solving one small second-order cone
  program. Not every instance with $k\ge2$ is degenerate.  For the regular simplex with $p$ vertices
  and unit edges in $\mathbb{R}^{p-1}$, \eqref{eq:degen1}--\eqref{eq:degen2} hold for $p=16$, $k=3$ but fail for
  $p=20$, $k=3$, where $v_3\approx0.0045$.
  On test set~C, conditions~\eqref{eq:degen1}--\eqref{eq:degen2} hold on all $840$
  instances with $k\ge3$, and solving $\mathrm{R}$ directly confirms $v_k=0$ on all of them; by Theorem~\ref{thm:degenerate}, they hold on none of the instances with $k=1$, consistent with Table~\ref{tab:rootD}, where the bare
  root gap is below $1$ at $k=1$ only.
\end{remark}

The global lower bounds~\eqref{eq:diamcut} and~\eqref{eq:clustercut} are constraints of the
form $L\ge\ell$ for a constant $\ell$ computed in preprocessing, and involve no other variable.
Their effect on the continuous relaxation is therefore the following.

\begin{proposition}\label{prop:rootLB}
  Let $\ell\le L^\star_{\le k}$ be a constant. The continuous relaxation of {\rm \Idenone} with
  the additional constraint $L\ge\ell$ has optimal value $\max\{\ell,v_k\}$.
\end{proposition}

\begin{proof}
  The hypothesis $\ell\le L^\star_{\le k}$ is needed only so that the constraint $L\ge\ell$ is
  valid for {\rm \Idenone}; the value of the relaxation does not depend on it. If $(x,y,f,L)$
  is feasible for $\mathrm{R}$, so is $(x,y,f,L')$ for every $L'\ge L$, as $L$ appears only on
  the right-hand side of~\eqref{eq:soc} and in $L\ge0$. Hence the set of values of $L$ feasible
  for $\mathrm{R}$ is $[v_k,\infty)$, and its intersection with $[\ell,\infty)$ has minimum
  $\max\{\ell,v_k\}$.
\end{proof}

With $\underline{L}:=\max\{\diam(P)/(p+k-1),\tfrac12 e_k\}$, Theorem~\ref{thm:degenerate} and
Proposition~\ref{prop:rootLB} show that, whenever~\eqref{eq:degen1}--\eqref{eq:degen2} hold, the
continuous relaxation of {\rm \Idenone} with both global lower bounds has value exactly
$\underline{L}$, so the root gap computed from the relaxation value is $1-\underline{L}/U_H$. On such instances
the root bound can only be raised by a valid inequality on $L$ or by a modification that
violates~\eqref{eq:degen1}--\eqref{eq:degen2}: the ordering~\eqref{eq:s1} does not raise it
(Proposition~\ref{prop:symrelax}), and neither does the edge elimination of
Section~\ref{sec:elim}, which only removes terminal--terminal edges, on which the relaxation
solution with $L=0$ constructed in the proof of Theorem~\ref{thm:degenerate} has $y_{ij}=0$. A stronger connectivity formulation raises the bound only on instances that are not
degenerate for that formulation, that is, where its analogue of~\eqref{eq:degen1}--\eqref{eq:degen2} fails. For the multi-commodity flow formulation of Section~\ref{sec:m2},
Proposition~\ref{prop:symrelax} holds with the same proof, and each terminal must receive a
full unit of capacity, so the first part of the proof of Theorem~\ref{thm:degenerate} gives the
stronger necessary condition $\rho_i(c)\le1-1/k$ for every $i\in P$. Fewer instances are
therefore degenerate for the multi-commodity formulation, but Proposition~\ref{prop:rootLB} applies to it as
well: its root bound equals $\underline{L}$ whenever its bare value does not exceed
$\underline{L}$, which is why the two root bounds coincide on almost every instance
(Section~\ref{sec:m2}). On test set~C, the bare multi-commodity relaxation has a value above $10^{-6}$
on $331$ of the $840$ instances with $k\ge3$, but its value never exceeds $\underline{L}$, so with both
global lower bounds the two relaxations have the same value $\underline{L}$ on all $840$.

}

\subsection{Summary of the implemented formulation and remaining challenges}
\label{sec:summary}

The implemented formulation comprises constraints \eqref{eq:soc}--\eqref{eq:x} on the candidate edge set $E^\star$ of~\eqref{eq:elim}, or on $E^\star_H$ of~\eqref{eq:elimmst} in the experiments listed in Section~\ref{sec:implementation}, both global lower bounds~\eqref{eq:diamcut} and~\eqref{eq:clustercut}, big-$M$ constants~\eqref{eq:bigMdef}, the symmetry-breaking constraints~\eqref{eq:s1}, and connectivity as in Formulation~{\rm \Idenone} or~{\rm \Identity}; the partial MIP start is $(L,x_S)=(U_H,\text{Steiner positions of }T_H)$.

Despite the strengthening techniques introduced above, four main challenges remain in solving the EBSTP exactly.
\begin{enumerate}[leftmargin=*,label={\bfseries C\arabic*},itemsep=0.4em]
  \item\label{c1} \textbf{Weakness of the bottleneck epigraph relaxation.} When $y_{ij}\in(0,1)$, the big-$M$ term $M_{ij}(1-y_{ij})$ in constraint~\eqref{eq:soc} renders the cone non-binding for any reasonable $M_{ij}$. The strengthening techniques of Sections~\ref{sec:bigM} and~\ref{sec:elim} act on the integer feasible region and the candidate edge set; they do not directly tighten the LP and second-order cone programming (SOCP) relaxations of $L$. { By Theorem~\ref{thm:degenerate}, the bare relaxation bound is exactly $0$ on all $840$ test instances with $k\ge3$, and it is positive on every instance with $k=1$ (Section~\ref{sec:rootexp}).} This relaxation weakness is analogous to that of the additive formulations MMX and FM~\citep*{fampa2016}. Only the global lower bounds~\eqref{eq:diamcut} and~\eqref{eq:clustercut} of Section~\ref{sec:diam} act directly on the relaxation of~$L$; together they remove the trivial $1.000$ gap and roughly halve the root gap for $k\ge 3$, but they do not close it.
  \item\label{c2} { \textbf{Permutation symmetry on Steiner labels.} The ordering~\eqref{eq:s1} is
    valid, but it does not eliminate the symmetry between Steiner points that share their first
    coordinate, and the solver does not detect the symmetry on its own. Designing a symmetry
    treatment that exploits the geometry of the EBSTP rather than a coordinate order remains open.}
  \item\label{c3} \textbf{Edge elimination for Steiner-incident edges.} Proposition~\ref{prop:mst}
    reduces the terminal--terminal edges to the $p-1$ edges of a terminal MST in every dimension.
    Steiner-incident edges have decision-dependent lengths and admit no analogous certificate.
  \item\label{c4} { \textbf{Growth of the formulation size with the budget.} With $E^\star_H$, the
    candidate edge set contains at most $p-1$ terminal--terminal edges and $pk+k(k-1)/2$
    Steiner-incident edges, each contributing one binary variable, one cone constraint and, in
    {\rm \Idenone}, two flow capacity constraints, so the formulation size grows as $O(pk+k^2)$.}
\end{enumerate}
Section~\ref{sec:computational} quantifies the practical impact of these four challenges, while Section~\ref{sec:conclusion} discusses the resulting directions for future research.

\FloatBarrier

\section{Computational study}
\label{sec:computational}

We address six questions in the computational study: \textbf{(A)}~validation of {\rm \Idenone} and {\rm \Identity} in the plane against the exact combinatorial baseline \ExactBCLT{}; \textbf{(B)}~comparison of {\rm \Idenone} and {\rm \Identity} across dimensions, instance families, and Steiner budgets, up to the limit of what each formulation solves; \textbf{(C)}~the behavior of {\rm \Idenone} and {\rm \Identity} on structured named instances that reach higher dimensions than the random benchmark; \textbf{(D)}~the strength of the {\rm \Idenone} root relaxation and the effect of the two global lower bounds of Section~\ref{sec:diam} on it; \textbf{(E)}~the practical effect of the symmetry-breaking scheme; \textbf{(F)}~%
the trade-off between native MISOCP and outer approximation, inside the solver's branch-and-bound and in a multi-tree loop.

\subsection{Implementation and benchmark design}
\label{sec:implementation}

\paragraph{Software and hardware} The formulations were implemented in Julia~1.12.2 using JuMP~\citep*{lubin2023} and Gurobi~13.0.2~\citep*{gurobi2025}. The cut separator of Section~\ref{sec:m2} uses the \texttt{Graphs.jl}~\citep*{graphs2021} maximum-flow oracle on the support graph with edge capacities $y_{ij}$; violated cuts at integer nodes are detected by connected-component analysis.
All experiments were run on a single workstation with a 13th-generation Intel Core i7-13650HX processor ($2.6$\,GHz base clock, 14~physical cores) and 64\,GB of RAM, running Linux under the Windows Subsystem for Linux~2; we solved every formulation on this platform.

\paragraph{Solver settings} Unless noted otherwise, runs use a solve time limit of \num{1000}\,s, a relative MIP gap tolerance of $10^{-4}$, feasibility tolerance $10^{-6}$, integer feasibility tolerance $10^{-6}$, and a fixed random seed. {\rm \Idenone}
uses two solver threads; {\rm \Identity} runs single-threaded because our lazy-cut separation reuses a per-solve callback buffer, and the serial combinatorial \ExactBCLT{} baseline likewise runs on one thread. Gurobi chooses its method for the cone constraints automatically (\texttt{MIQCPMethod}${}=-1$), except in Section~\ref{sec:oaexp}. Each formulation receives the partial MIP start of Section~\ref{sec:UH} (the bisection value $U_H$ and the Steiner positions of $T_H$; completion of the edge binaries into an incumbent is left to the solver). Both global lower bounds of Section~\ref{sec:diam}, constraints~\eqref{eq:diamcut} and~\eqref{eq:clustercut}, are enabled except in the experiments that omit one or both bounds (the ``bare'' and ``$+$ diam'' configurations of Table~\ref{tab:rootD} and {  the experiments with and without the global lower bounds } in Section~\ref{sec:rootexp}).
{ The experiments of Sections~\ref{sec:mainexp}--\ref{sec:oaexp} restrict the terminal--terminal candidate edges to $E^\star_H$ of~\eqref{eq:elimmst}. The planar validation of Section~\ref{sec:planarexp} uses $E^\star$ of~\eqref{eq:elim}. Departures from this rule are stated where they occur.}
The tolerances left symbolic in Sections~\ref{sec:m2} and~\ref{sec:elim} take the following values throughout: the maximum-flow separation threshold uses $\varepsilon=10^{-6}$, and the edge-elimination comparison of~\eqref{eq:elim} uses $\varepsilon=10^{-8}$.

\paragraph{Instance families} Two families are used. \emph{Uniform}: terminals drawn independently from $[0,1]^d$. \emph{Clustered}: terminals drawn from a mixture of four compact Gaussian clusters with random centers and spread~$0.2$.

\paragraph{Benchmark test sets}

We use three test sets:

\begin{itemize}[leftmargin=*,nosep]
  \item Test set~A (planar validation): $d=2$, $p\in\{10,15,20,25\}$, $k=1$, seeds $1$ to $5$, both families; $40$~instances total.
  \item Test set~B (smaller set): $d\in\{2,3\}$, $p\in\{10,15,20,25\}$, $k\in\{1,3,5\}$, seeds $1$ to $5$, both families; $240$~instances total.
  \item Test set~C (larger set): $d\in\{2,3,4\}$, $p\in\{10,15,20,25,30,40,50\}$,
  $k\in\{1,3,5,7,9\}$, seeds $1$ to $5$, both families; $1050$~instances total, $10$ for each
  combination of $d$, $p$ and~$k$.
\end{itemize}

{ On test set~B, {\rm \Idenone} in its default configuration leaves one instance, with $d=2$, $p=25$
and $k=5$, at the time limit, and most of its variants solve it. Where a variant of {\rm \Idenone} solves one
instance more on test set~B, it is this one, and where it changes a comparison of total solve times we also
report the total over the other $239$ instances.}

\paragraph{Reporting conventions} For each row we report the number of instances solved to proven optimality within the time limit (``solv.''), the median solve time (wall-clock) ``$t_{\med}$'' and the shifted geometric mean solve time ``$t_{\sgm}$'' (shift $s=10$\,s~\citep*{achterberg2007}) in seconds, the median branch-and-bound node count ``nodes'' over solved instances, and the maximum remaining optimality gap ``gap$^*$'' over the unsolved instances (equal to $0\%$ when all are solved); the optimality gap of an instance is $(\mathit{UB}-\mathit{LB})/\mathit{UB}$, where $\mathit{UB}$ is the smaller of its best feasible objective and $U_H$, and $\mathit{LB}$ its best dual bound, reported as a percentage~\citep*{pinto2023,maculan2025}. The shifted geometric mean is insensitive to outliers in the long tail of solve times while penalizing failures: instances that hit the time limit enter the mean at the time-limit value of \num{1000}\,s. In each table, boldface marks the unique best on an informative metric within the comparison group (typically a row, or a $k$-block in Table~\ref{tab:rootD}); when the best value is shared by more than one method, no entry in that group is bolded, as the comparison is then uninformative, in particular for sub-second ties on solve time and for ties on solve count when every method solves all instances in the row. Median node counts are computed over each method's solved instances; when solve counts differ, the node-count comparison is indicative only, and we therefore use boldface for a node-count entry only in rows where the methods solve equally many instances.

\subsection{Validation against \ExactBCLT{} in the plane}
\label{sec:planarexp}

The first experiment compares the optimal value reported by {\rm \Idenone} and~{\rm \Identity} to the value found by the combinatorial planar algorithm \ExactBCLT{} described in~\citep*{bae2009exact}, on the $40$ instances of test set~A.

Table~\ref{tab:planar} reports the median solve times of all three methods. We omit columns for solve count and instance-level discrepancy: every method solves every instance, and all discrepancies fall within the configured MIP gap tolerance of $10^{-4}$. At this resolution, the discrepancy figures (typically of order $10^{-6}$ or smaller) reflect solver convergence parameters rather than properties of the formulations.

\begin{table}[!ht]
  \centering
   
  \caption{Planar timing comparison on test set~A: all three methods solve
    the $40$ instances, with the conic objectives agreeing with \ExactBCLT{}
    within $4\times10^{-6}$. Solve times are medians (s) over five seeds;
    ``$<0.05$'' marks a sub-$0.05$\,s cell.}
  \label{tab:planar}
  \small
  \setlength{\tabcolsep}{6pt}
  \renewcommand*{\arraystretch}{1.3}
  \begin{tabular}{
      S[table-format=1.0] S[table-format=2.0] S[table-format=1.0] S[table-format=1.0]
      l
    S[table-format=2.1] S[table-format=2.1] S[table-format=2.1]}
    \toprule
    {$d$} & {$p$} & {$k$} & {\#inst} & {family} & {\ExactBCLT{}} & {\Mone{}} & {\Mtwo{}}\\
    \midrule
2 & 10 & 1 & 5 & clustered & {$<0.05$} & {$<0.05$} & {$<0.05$}\\
2 & 15 & 1 & 5 & clustered & 0.1 & {$<0.05$} & 0.1\\
2 & 20 & 1 & 5 & clustered & {$<0.05$} & {$<0.05$} & 0.2\\
2 & 25 & 1 & 5 & clustered & 0.1 & {$<0.05$} & 0.2\\
2 & 10 & 1 & 5 & uniform & 0.1 & {$<0.05$} & {$<0.05$}\\
2 & 15 & 1 & 5 & uniform & {$<0.05$} & 0.1 & {$<0.05$}\\
2 & 20 & 1 & 5 & uniform & {$<0.05$} & 0.1 & 0.2\\
2 & 25 & 1 & 5 & uniform & 0.1 & 0.1 & 0.2\\
    \bottomrule
  \end{tabular}
\end{table}

All three methods are fast on these planar instances: every median in Table~\ref{tab:planar} is below $0.25$\,s. We do not interpret these timings as evidence of strict outperformance, since our \ExactBCLT{} implementation lacks the low-level optimizations of a dedicated computational-geometry toolkit; the timing comparison does show that the conic formulations are not penalized by the planar specialization of the baseline. The experiment confirms the correctness of {\rm \Idenone} and {\rm \Identity}.

{ 
We also compared the three methods at $k=2$, on $30$ planar instances with $p\in\{6,8,10\}$. In
this comparison, {\rm \Idenone} and {\rm \Identity} use the candidate edge set $E^\star_H$
of~\eqref{eq:elimmst}. Our implementation of \ExactBCLT{} solves the $10$ instances with $p=6$, in
$106$\,s on average, and reaches the time limit of \num{1000}\,s on all instances with
$p\in\{8,10\}$. {\rm \Idenone} and {\rm \Identity} solve all $30$ instances, in at most $0.27$\,s
and $0.37$\,s, respectively. On the instances that all three methods solve, the optimal values
agree within $3\times10^{-6}$. As noted above, our implementation of the combinatorial algorithm
lacks the low-level optimizations of a dedicated computational-geometry code, so this comparison
measures implementations rather than algorithms.
}

\subsection{Main comparison of \Mone{} and \Mtwo{}}
\label{sec:mainexp}

The main experiment runs {\rm \Idenone} and {\rm \Identity} on the $1050$ instances of test set~C, both with the candidate edge set $E^\star_H$ of~\eqref{eq:elimmst}.

\begin{table}[!ht]
  \centering
   
  \caption{\Mone{} versus \Mtwo{} on the $1050$ instances of test set~C, both with the MST restriction, aggregated by ambient dimension, by dimension and instance family, and by dimension and Steiner budget~$k$. ``solv.''/``\%'' count proven-optimal solves within \num{1000}\,s; $t_{\med}$/$t_{\sgm}$ are the median and shifted-geometric-mean solve times (s, shift~$10$); ``nodes'' the median branch-and-bound node count over solved instances; ``gap$^*$'' the worst gap over unsolved instances; ``all'' aggregates the omitted parameters. Bold marks the per-row best; node counts are bold only on equal solved sets.}
  \label{tab:main}
  \footnotesize
  \setlength{\tabcolsep}{2pt}
  \renewcommand*{\arraystretch}{1.3}
  \begin{tabular*}{\textwidth}{@{\extracolsep{\fill}}
      c l c S[table-format=3.0]
      *2{S[table-format=3.0] S[table-format=3.1] S[table-format=4.1] S[table-format=4.1]
    S[table-format=6.0] S[table-format=2.2]}}
    \toprule
    & & & & \multicolumn{6}{c}{\Mone{}} & \multicolumn{6}{c}{\Mtwo{}}\\
    \cmidrule(lr){5-10}\cmidrule(lr){11-16}
    {$d$} & {family} & {$k$} & {\#inst} & {solv.} & {\%} & {$t_{\med}$} & {$t_{\sgm}$} & {nodes} & {gap$^*$}
    & {solv.} & {\%} & {$t_{\med}$} & {$t_{\sgm}$} & {nodes} & {gap$^*$}\\
    \midrule
\multicolumn{16}{@{}l}{\itshape By ambient dimension}\\
2 & {all} & {all} & 350 & \best{280} & \best{80.0} & \best{21.8} & \best{54.4} & 3071 & \best{75.63} & 253 & 72.3 & 34.7 & 73.0 & 1839 & 78.56\\
3 & {all} & {all} & 350 & \best{244} & \best{69.7} & \best{41.0} & \best{78.5} & 2927 & \best{72.86} & 181 & 51.7 & 721.6 & 192.2 & 2937 & 79.94\\
4 & {all} & {all} & 350 & \best{209} & \best{59.7} & \best{169.5} & \best{129.7} & 6418 & \best{71.83} & 132 & 37.7 & 1000.0 & 289.5 & 1576 & 77.55\\
\midrule
\multicolumn{16}{@{}l}{\itshape By dimension and instance family}\\
2 & clustered & {all} & 175 & \best{119} & \best{68.0} & \best{53.1} & \best{80.0} & 2642 & \best{75.63} & 111 & 63.4 & 81.3 & 96.4 & 1618 & 78.56\\
2 & uniform & {all} & 175 & \best{161} & \best{92.0} & \best{16.8} & \best{36.1} & 3551 & \best{42.67} & 142 & 81.1 & 24.7 & 54.8 & 1876 & 59.52\\
3 & clustered & {all} & 175 & \best{107} & \best{61.1} & \best{89.0} & \best{106.7} & 3066 & \best{72.86} & 86 & 49.1 & 1000.0 & 222.7 & 4458 & 79.94\\
3 & uniform & {all} & 175 & \best{137} & \best{78.3} & \best{23.0} & \best{57.1} & 2788 & \best{46.11} & 95 & 54.3 & 380.6 & 165.7 & 2599 & 72.43\\
4 & clustered & {all} & 175 & \best{95} & \best{54.3} & \best{479.8} & \best{161.1} & 6081 & \best{71.83} & 61 & 34.9 & 1000.0 & 315.2 & 1822 & 77.55\\
4 & uniform & {all} & 175 & \best{114} & \best{65.1} & \best{74.8} & \best{104.1} & 6448 & \best{68.88} & 71 & 40.6 & 1000.0 & 265.7 & 1016 & 73.31\\
\midrule
\multicolumn{16}{@{}l}{\itshape By dimension and Steiner budget}\\
2 & {all} & 1 & 70 & 70 & 100.0 & {\bestmath{<0.05}} & \best{0.1} & \best{1} & 0 & 70 & 100.0 & 0.1 & 0.2 & 56 & 0\\
2 & {all} & 3 & 70 & 70 & 100.0 & 6.7 & 7.8 & 1449 & 0 & 70 & 100.0 & \best{3.7} & \best{6.3} & \best{1081} & 0\\
2 & {all} & 5 & 70 & 63 & 90.0 & \best{19.5} & \best{42.2} & 6311 & 50.70 & 63 & 90.0 & 37.7 & 62.8 & \best{5246} & \best{50.68}\\
2 & {all} & 7 & 70 & \best{48} & \best{68.6} & \best{137.1} & \best{186.5} & 40136 & \best{56.85} & 41 & 58.6 & 403.4 & 339.2 & 16695 & 72.75\\
2 & {all} & 9 & 70 & \best{29} & \best{41.4} & 1000.0 & \best{591.4} & 122239 & \best{75.63} & 9 & 12.9 & 1000.0 & 928.0 & 37925 & 78.56\\
3 & {all} & 1 & 70 & 70 & 100.0 & {\bestmath{<0.05}} & \best{0.1} & \best{19} & 0 & 70 & 100.0 & 0.1 & 0.4 & 84 & 0\\
3 & {all} & 3 & 70 & 70 & 100.0 & \best{8.3} & \best{10.5} & \best{2596} & 0 & 70 & 100.0 & 41.8 & 46.0 & 6264 & 0\\
3 & {all} & 5 & 70 & \best{62} & \best{88.6} & \best{33.6} & \best{59.4} & 12050 & \best{50.93} & 39 & 55.7 & 762.8 & 560.9 & 18593 & 77.88\\
3 & {all} & 7 & 70 & \best{38} & \best{54.3} & \best{596.0} & \best{374.9} & 44850 & \best{70.03} & 2 & 2.9 & 1000.0 & 991.8 & 26230 & 79.94\\
3 & {all} & 9 & 70 & \best{4} & \best{5.7} & 1000.0 & \best{964.7} & 114849 & \best{72.86} & 0 & 0.0 & 1000.0 & 1000.0 & {--} & 74.95\\
4 & {all} & 1 & 70 & 70 & 100.0 & \best{0.1} & \best{0.5} & \best{331} & 0 & 70 & 100.0 & 0.8 & 1.6 & 338 & 0\\
4 & {all} & 3 & 70 & \best{70} & \best{100.0} & \best{16.4} & \best{24.0} & 7762 & \best{0} & 56 & 80.0 & 189.1 & 204.7 & 11352 & 59.49\\
4 & {all} & 5 & 70 & \best{55} & \best{78.6} & \best{130.7} & \best{157.7} & 28896 & \best{57.74} & 6 & 8.6 & 1000.0 & 936.2 & 15566 & 77.55\\
4 & {all} & 7 & 70 & \best{14} & \best{20.0} & 1000.0 & \best{874.1} & 108408 & \best{71.83} & 0 & 0.0 & 1000.0 & 1000.0 & {--} & 75.79\\
4 & {all} & 9 & 70 & 0 & 0.0 & 1000.0 & 1000.0 & {--} & \best{70.05} & 0 & 0.0 & 1000.0 & 1000.0 & {--} & 74.08\\
    \bottomrule
  \end{tabular*}
\end{table}

{ 
The by-dimension block of Table~\ref{tab:main} shows three results. First, {\rm \Idenone} is
faster than {\rm \Identity} in median solve time in every dimension: the ratios of median solve
times are $1.6$, $17.6$ and $5.9$ at $d=2,3,4$, the last bounded by the time limit, which is the
median solve time of {\rm\Identity} at $d=4$. The shifted geometric mean tracks the same ordering in
every dimension ($54.4$ versus $73.0$, $78.5$ versus $192.2$, and $129.7$ versus $289.5$\,s). Second, the
solve rates diverge with dimension: {\rm\Idenone} solves $280$, $244$ and $209$ of the $350$ instances at
$d=2,3,4$, and {\rm\Identity} $253$, $181$ and $132$. Third, the instances left unsolved end with
smaller gaps for {\rm\Idenone}: their median gap is $20.2$, $40.9$ and $41.6\%$ at $d=2,3,4$, against
$43.8$, $50.3$ and $49.6\%$ for {\rm\Identity}. { {\rm \Idenone} runs on two solver threads and
{\rm \Identity} on one (Section~\ref{sec:implementation}). Run on one thread, {\rm \Idenone} solves $696$ of
the $1050$ instances, against $733$ on two, and still solves more instances than {\rm \Identity} in every
dimension: $266$, $236$ and $194$ at $d=2,3,4$, against $253$, $181$ and $132$. It solves $132$
instances that {\rm \Identity} does not, against $2$ ($p<10^{-36}$), and on the instances both solve,
{\rm \Identity} is slower by median factors of $1.8$, $7.4$ and $7.5$. The thread setting therefore does
not explain the difference between the formulations.}

The breakdown by family shows that the clustered instances are the harder ones for both
formulations in every dimension, and that the advantage of {\rm\Idenone} holds in both families: at
$d=3$ and $d=4$ it solves $42$ and $43$ more uniform instances than {\rm\Identity}, and $21$ and $34$ more
clustered ones. Its median solve time is shorter by factors of $1.5$, $16.5$ and $13.4$ on the
uniform instances at $d=2,3,4$, and of $1.5$, $11.2$ and $2.1$ on the clustered ones.

The budget breakdown shows that difficulty grows { quickly} with the Steiner budget~$k$. At $k=1$,
both formulations solve every instance. At $d=4$, $k=5$, {\rm \Idenone} solves $55$ of the $70$
instances, whereas {\rm \Identity} solves $6$, with its median solve time at the \num{1000}\,s limit
against $130.7$\,s for {\rm \Idenone}. At $d=3$, $k=7$, the counts are $38$ against $2$.
}

{ Table~\ref{tab:matrix} reports, for each combination of $d$, $p$ and~$k$, how many of
its $10$ instances each formulation solves. {\rm\Idenone} solves every instance with $k\le3$;
beyond that, the number of solved instances falls with the budget and the dimension, and in most
series also with the number of terminals. With
$k=9$, {\rm\Idenone} solves $29$ of the $70$ planar instances, $4$ at $d=3$, all with $p=10$, and none at
$d=4$; {\rm\Identity} solves $9$, none and none. {\rm\Idenone} solves more instances than
{\rm\Identity} in $42$ of the $105$ combinations, fewer in one, and as many in the other $62$. The
instances that {\rm\Idenone} leaves unsolved end with a median gap of $39.5\%$, and those of
{\rm\Identity} with $49.7\%$.}

{ 
\begin{table}[!ht]
  \centering
   
  \caption{Number of the $10$ instances of each combination of test set~C (seeds $1$ to $5$, both
    families) that each formulation, with the MST restriction, solves to proven optimality within
    \num{1000}\,s. Bold marks the formulation that solves more instances of a combination.}
  \label{tab:matrix}
  \small
  \setlength{\tabcolsep}{3.5pt}
  \begin{tabular}{cc ccccccc ccccccc}
    \toprule
    & & \multicolumn{7}{c}{\Mone{}} & \multicolumn{7}{c}{\Mtwo{}}\\
    \cmidrule(lr){3-9}\cmidrule(lr){10-16}
    $d$ & $k$ & $p=10$ & $15$ & $20$ & $25$ & $30$ & $40$ & $50$ & $p=10$ & $15$ & $20$ & $25$ & $30$ & $40$ & $50$\\
    \midrule
2 & 1 & 10 & 10 & 10 & 10 & 10 & 10 & 10 & 10 & 10 & 10 & 10 & 10 & 10 & 10\\
2 & 3 & 10 & 10 & 10 & 10 & 10 & 10 & 10 & 10 & 10 & 10 & 10 & 10 & 10 & 10\\
2 & 5 & 10 & 10 & 10 & 9 & 9 & \best{9} & 6 & 10 & 10 & 10 & \best{10} & 9 & 8 & 6\\
2 & 7 & 10 & \best{10} & \best{7} & \best{6} & \best{6} & \best{6} & \best{3} & 10 & 8 & 6 & 5 & 5 & 5 & 2\\
2 & 9 & \best{8} & \best{6} & \best{5} & \best{3} & \best{5} & \best{2} & 0 & 4 & 2 & 1 & 2 & 0 & 0 & 0\\
\midrule
3 & 1 & 10 & 10 & 10 & 10 & 10 & 10 & 10 & 10 & 10 & 10 & 10 & 10 & 10 & 10\\
3 & 3 & 10 & 10 & 10 & 10 & 10 & 10 & 10 & 10 & 10 & 10 & 10 & 10 & 10 & 10\\
3 & 5 & 10 & 10 & \best{10} & \best{10} & \best{9} & \best{8} & \best{5} & 10 & 10 & 7 & 5 & 5 & 1 & 1\\
3 & 7 & \best{9} & \best{8} & \best{7} & \best{4} & \best{5} & \best{3} & \best{2} & 2 & 0 & 0 & 0 & 0 & 0 & 0\\
3 & 9 & \best{4} & 0 & 0 & 0 & 0 & 0 & 0 & 0 & 0 & 0 & 0 & 0 & 0 & 0\\
\midrule
4 & 1 & 10 & 10 & 10 & 10 & 10 & 10 & 10 & 10 & 10 & 10 & 10 & 10 & 10 & 10\\
4 & 3 & 10 & 10 & 10 & \best{10} & \best{10} & \best{10} & \best{10} & 10 & 10 & 10 & 9 & 8 & 4 & 5\\
4 & 5 & \best{10} & \best{10} & \best{9} & \best{9} & \best{7} & \best{5} & \best{5} & 3 & 2 & 0 & 1 & 0 & 0 & 0\\
4 & 7 & \best{8} & \best{3} & \best{1} & \best{1} & 0 & \best{1} & 0 & 0 & 0 & 0 & 0 & 0 & 0 & 0\\
4 & 9 & 0 & 0 & 0 & 0 & 0 & 0 & 0 & 0 & 0 & 0 & 0 & 0 & 0 & 0\\
    \bottomrule
  \end{tabular}
\end{table}
}

\begin{figure}[htbp]
  \centering
\begin{tikzpicture}
\begin{semilogxaxis}
  [
    steineraxis,
    width = 0.82\linewidth,
    height = 0.45\linewidth,
    xlabel = {$\tau$ (performance ratio relative to faster method)},
    ylabel = {$\rho_s(\tau)$ (fraction of instances)},
    xmin = 1, xmax = 1e4,
    ymin = 0, ymax = 1.0,
    grid = major,
    grid style = {line width = .1pt, draw = gray!20},
    major grid style = {line width = .2pt, draw = gray!50},
    minor tick num = 8,
    log basis x = 10,
    legend pos = south east,
    legend style = {font = \small, black, draw = gray!60, fill = white},
    tick label style = {black, font = \small},
    label style = {black, font = \small},
    every axis plot/.append style = {thick},
  ]
  \addplot+[const plot, no marks, color = btterminal, solid, ultra thick] coordinates {
    (1,0.0) (1,0.0010) (1,0.0019) (1,0.0029)
    (1,0.0038) (1,0.0048) (1,0.0057) (1,0.0067)
    (1,0.0076) (1,0.0086) (1,0.0095) (1,0.0105)
    (1,0.0114) (1,0.0124) (1,0.0133) (1,0.0143)
    (1,0.0152) (1,0.0162) (1,0.0171) (1,0.0181)
    (1,0.0190) (1,0.0200) (1,0.0210) (1,0.0219)
    (1,0.0229) (1,0.0238) (1,0.0248) (1,0.0257)
    (1,0.0267) (1,0.0276) (1,0.0286) (1,0.0295)
    (1,0.0305) (1,0.0314) (1,0.0324) (1,0.0333)
    (1,0.0343) (1,0.0352) (1,0.0362) (1,0.0371)
    (1,0.0381) (1,0.0390) (1,0.0400) (1,0.0410)
    (1,0.0419) (1,0.0429) (1,0.0438) (1,0.0448)
    (1,0.0457) (1,0.0467) (1,0.0476) (1,0.0486)
    (1,0.0495) (1,0.0505) (1,0.0514) (1,0.0524)
    (1,0.0533) (1,0.0543) (1,0.0552) (1,0.0562)
    (1,0.0571) (1,0.0581) (1,0.0590) (1,0.0600)
    (1,0.0610) (1,0.0619) (1,0.0629) (1,0.0638)
    (1,0.0648) (1,0.0657) (1,0.0667) (1,0.0676)
    (1,0.0686) (1,0.0695) (1,0.0705) (1,0.0714)
    (1,0.0724) (1,0.0733) (1,0.0743) (1,0.0752)
    (1,0.0762) (1,0.0771) (1,0.0781) (1,0.0790)
    (1,0.0800) (1,0.0810) (1,0.0819) (1,0.0829)
    (1,0.0838) (1,0.0848) (1,0.0857) (1,0.0867)
    (1,0.0876) (1,0.0886) (1,0.0895) (1,0.0905)
    (1,0.0914) (1,0.0924) (1,0.0933) (1,0.0943)
    (1,0.0952) (1,0.0962) (1,0.0971) (1,0.0981)
    (1,0.0990) (1,0.1000) (1,0.1010) (1,0.1019)
    (1,0.1029) (1,0.1038) (1,0.1048) (1,0.1057)
    (1,0.1067) (1,0.1076) (1,0.1086) (1,0.1095)
    (1,0.1105) (1,0.1114) (1,0.1124) (1,0.1133)
    (1,0.1143) (1,0.1152) (1,0.1162) (1,0.1171)
    (1,0.1181) (1,0.1190) (1,0.1200) (1,0.1210)
    (1,0.1219) (1,0.1229) (1,0.1238) (1,0.1248)
    (1,0.1257) (1,0.1267) (1,0.1276) (1,0.1286)
    (1,0.1295) (1,0.1305) (1,0.1314) (1,0.1324)
    (1,0.1333) (1,0.1343) (1,0.1352) (1,0.1362)
    (1,0.1371) (1,0.1381) (1,0.1390) (1,0.1400)
    (1,0.1410) (1,0.1419) (1,0.1429) (1,0.1438)
    (1,0.1448) (1,0.1457) (1,0.1467) (1,0.1476)
    (1,0.1486) (1,0.1495) (1,0.1505) (1,0.1514)
    (1,0.1524) (1,0.1533) (1,0.1543) (1,0.1552)
    (1,0.1562) (1,0.1571) (1,0.1581) (1,0.1590)
    (1,0.1600) (1,0.1610) (1,0.1619) (1,0.1629)
    (1,0.1638) (1,0.1648) (1,0.1657) (1,0.1667)
    (1,0.1676) (1,0.1686) (1,0.1695) (1,0.1705)
    (1,0.1714) (1,0.1724) (1,0.1733) (1,0.1743)
    (1,0.1752) (1,0.1762) (1,0.1771) (1,0.1781)
    (1,0.1790) (1,0.1800) (1,0.1810) (1,0.1819)
    (1,0.1829) (1,0.1838) (1,0.1848) (1,0.1857)
    (1,0.1867) (1,0.1876) (1,0.1886) (1,0.1895)
    (1,0.1905) (1,0.1914) (1,0.1924) (1,0.1933)
    (1,0.1943) (1,0.1952) (1,0.1962) (1,0.1971)
    (1,0.1981) (1,0.1990) (1,0.2000) (1,0.2010)
    (1,0.2019) (1,0.2029) (1,0.2038) (1,0.2048)
    (1,0.2057) (1,0.2067) (1,0.2076) (1,0.2086)
    (1,0.2095) (1,0.2105) (1,0.2114) (1,0.2124)
    (1,0.2133) (1,0.2143) (1,0.2152) (1,0.2162)
    (1,0.2171) (1,0.2181) (1,0.2190) (1,0.2200)
    (1,0.2210) (1,0.2219) (1,0.2229) (1,0.2238)
    (1,0.2248) (1,0.2257) (1,0.2267) (1,0.2276)
    (1,0.2286) (1,0.2295) (1,0.2305) (1,0.2314)
    (1,0.2324) (1,0.2333) (1,0.2343) (1,0.2352)
    (1,0.2362) (1,0.2371) (1,0.2381) (1,0.2390)
    (1,0.2400) (1,0.2410) (1,0.2419) (1,0.2429)
    (1,0.2438) (1,0.2448) (1,0.2457) (1,0.2467)
    (1,0.2476) (1,0.2486) (1,0.2495) (1,0.2505)
    (1,0.2514) (1,0.2524) (1,0.2533) (1,0.2543)
    (1,0.2552) (1,0.2562) (1,0.2571) (1,0.2581)
    (1,0.2590) (1,0.2600) (1,0.2610) (1,0.2619)
    (1,0.2629) (1,0.2638) (1,0.2648) (1,0.2657)
    (1,0.2667) (1,0.2676) (1,0.2686) (1,0.2695)
    (1,0.2705) (1,0.2714) (1,0.2724) (1,0.2733)
    (1,0.2743) (1,0.2752) (1,0.2762) (1,0.2771)
    (1,0.2781) (1,0.2790) (1,0.2800) (1,0.2810)
    (1,0.2819) (1,0.2829) (1,0.2838) (1,0.2848)
    (1,0.2857) (1,0.2867) (1,0.2876) (1,0.2886)
    (1,0.2895) (1,0.2905) (1,0.2914) (1,0.2924)
    (1,0.2933) (1,0.2943) (1,0.2952) (1,0.2962)
    (1,0.2971) (1,0.2981) (1,0.2990) (1,0.3000)
    (1,0.3010) (1,0.3019) (1,0.3029) (1,0.3038)
    (1,0.3048) (1,0.3057) (1,0.3067) (1,0.3076)
    (1,0.3086) (1,0.3095) (1,0.3105) (1,0.3114)
    (1,0.3124) (1,0.3133) (1,0.3143) (1,0.3152)
    (1,0.3162) (1,0.3171) (1,0.3181) (1,0.3190)
    (1,0.3200) (1,0.3210) (1,0.3219) (1,0.3229)
    (1,0.3238) (1,0.3248) (1,0.3257) (1,0.3267)
    (1,0.3276) (1,0.3286) (1,0.3295) (1,0.3305)
    (1,0.3314) (1,0.3324) (1,0.3333) (1,0.3343)
    (1,0.3352) (1,0.3362) (1,0.3371) (1,0.3381)
    (1,0.3390) (1,0.3400) (1,0.3410) (1,0.3419)
    (1,0.3429) (1,0.3438) (1,0.3448) (1,0.3457)
    (1,0.3467) (1,0.3476) (1,0.3486) (1,0.3495)
    (1,0.3505) (1,0.3514) (1,0.3524) (1,0.3533)
    (1,0.3543) (1,0.3552) (1,0.3562) (1,0.3571)
    (1,0.3581) (1,0.3590) (1,0.3600) (1,0.3610)
    (1,0.3619) (1,0.3629) (1,0.3638) (1,0.3648)
    (1,0.3657) (1,0.3667) (1,0.3676) (1,0.3686)
    (1,0.3695) (1,0.3705) (1,0.3714) (1,0.3724)
    (1,0.3733) (1,0.3743) (1,0.3752) (1,0.3762)
    (1,0.3771) (1,0.3781) (1,0.3790) (1,0.3800)
    (1,0.3810) (1,0.3819) (1,0.3829) (1,0.3838)
    (1,0.3848) (1,0.3857) (1,0.3867) (1,0.3876)
    (1,0.3886) (1,0.3895) (1,0.3905) (1,0.3914)
    (1,0.3924) (1,0.3933) (1,0.3943) (1,0.3952)
    (1,0.3962) (1,0.3971) (1,0.3981) (1,0.3990)
    (1,0.4000) (1,0.4010) (1,0.4019) (1,0.4029)
    (1,0.4038) (1,0.4048) (1,0.4057) (1,0.4067)
    (1,0.4076) (1,0.4086) (1,0.4095) (1,0.4105)
    (1,0.4114) (1,0.4124) (1,0.4133) (1,0.4143)
    (1,0.4152) (1,0.4162) (1,0.4171) (1,0.4181)
    (1,0.4190) (1,0.4200) (1,0.4210) (1,0.4219)
    (1,0.4229) (1,0.4238) (1,0.4248) (1,0.4257)
    (1,0.4267) (1,0.4276) (1,0.4286) (1,0.4295)
    (1,0.4305) (1,0.4314) (1,0.4324) (1,0.4333)
    (1,0.4343) (1,0.4352) (1,0.4362) (1,0.4371)
    (1,0.4381) (1,0.4390) (1,0.4400) (1,0.4410)
    (1,0.4419) (1,0.4429) (1,0.4438) (1,0.4448)
    (1,0.4457) (1,0.4467) (1,0.4476) (1,0.4486)
    (1,0.4495) (1,0.4505) (1,0.4514) (1,0.4524)
    (1,0.4533) (1,0.4543) (1,0.4552) (1,0.4562)
    (1,0.4571) (1,0.4581) (1,0.4590) (1,0.4600)
    (1,0.4610) (1,0.4619) (1,0.4629) (1,0.4638)
    (1,0.4648) (1,0.4657) (1,0.4667) (1,0.4676)
    (1,0.4686) (1,0.4695) (1,0.4705) (1,0.4714)
    (1,0.4724) (1,0.4733) (1,0.4743) (1,0.4752)
    (1,0.4762) (1,0.4771) (1,0.4781) (1,0.4790)
    (1,0.4800) (1,0.4810) (1,0.4819) (1,0.4829)
    (1,0.4838) (1,0.4848) (1,0.4857) (1,0.4867)
    (1,0.4876) (1,0.4886) (1,0.4895) (1,0.4905)
    (1,0.4914) (1,0.4924) (1,0.4933) (1,0.4943)
    (1,0.4952) (1,0.4962) (1,0.4971) (1,0.4981)
    (1,0.4990) (1,0.5000) (1,0.5010) (1,0.5019)
    (1,0.5029) (1,0.5038) (1,0.5048) (1,0.5057)
    (1,0.5067) (1,0.5076) (1,0.5086) (1,0.5095)
    (1,0.5105) (1,0.5114) (1,0.5124) (1,0.5133)
    (1,0.5143) (1,0.5152) (1,0.5162) (1,0.5171)
    (1,0.5181) (1,0.5190) (1,0.5200) (1,0.5210)
    (1,0.5219) (1,0.5229) (1,0.5238) (1,0.5248)
    (1,0.5257) (1,0.5267) (1,0.5276) (1,0.5286)
    (1,0.5295) (1,0.5305) (1,0.5314) (1,0.5324)
    (1,0.5333) (1,0.5343) (1,0.5352) (1,0.5362)
    (1,0.5371) (1,0.5381) (1,0.5390) (1,0.5400)
    (1,0.5410) (1,0.5419) (1,0.5429) (1,0.5438)
    (1,0.5448) (1,0.5457) (1,0.5467) (1,0.5476)
    (1,0.5486) (1,0.5495) (1,0.5505) (1,0.5514)
    (1,0.5524) (1,0.5533) (1,0.5543) (1,0.5552)
    (1,0.5562) (1,0.5571) (1,0.5581) (1,0.5590)
    (1,0.5600) (1,0.5610) (1,0.5619) (1,0.5629)
    (1,0.5638) (1,0.5648) (1,0.5657) (1,0.5667)
    (1,0.5676) (1,0.5686) (1,0.5695) (1,0.5705)
    (1,0.5714) (1,0.5724) (1,0.5733) (1,0.5743)
    (1,0.5752) (1,0.5762) (1,0.5771) (1,0.5781)
    (1,0.5790) (1,0.5800) (1,0.5810) (1,0.5819)
    (1,0.5829) (1,0.5838) (1,0.5848) (1,0.5857)
    (1,0.5867) (1,0.5876) (1,0.5886) (1,0.5895)
    (1,0.5905) (1,0.5914) (1,0.5924) (1,0.5933)
    (1,0.5943) (1,0.5952) (1,0.5962) (1,0.5971)
    (1,0.5981) (1,0.5990) (1,0.6000) (1,0.6010)
    (1,0.6019) (1,0.6029) (1,0.6038) (1,0.6048)
    (1,0.6057) (1.015,0.6067) (1.033,0.6076) (1.036,0.6086)
    (1.039,0.6095) (1.048,0.6105) (1.054,0.6114) (1.058,0.6124)
    (1.113,0.6133) (1.121,0.6143) (1.121,0.6152) (1.133,0.6162)
    (1.136,0.6171) (1.136,0.6181) (1.178,0.6190) (1.178,0.6200)
    (1.188,0.6210) (1.189,0.6219) (1.196,0.6229) (1.196,0.6238)
    (1.213,0.6248) (1.223,0.6257) (1.224,0.6267) (1.229,0.6276)
    (1.245,0.6286) (1.267,0.6295) (1.293,0.6305) (1.296,0.6314)
    (1.3,0.6324) (1.301,0.6333) (1.31,0.6343) (1.314,0.6352)
    (1.324,0.6362) (1.329,0.6371) (1.329,0.6381) (1.35,0.6390)
    (1.378,0.6400) (1.383,0.6410) (1.469,0.6419) (1.498,0.6429)
    (1.508,0.6438) (1.512,0.6448) (1.52,0.6457) (1.52,0.6467)
    (1.539,0.6476) (1.542,0.6486) (1.592,0.6495) (1.611,0.6505)
    (1.651,0.6514) (1.657,0.6524) (1.685,0.6533) (1.698,0.6543)
    (1.72,0.6552) (1.789,0.6562) (1.79,0.6571) (1.823,0.6581)
    (1.838,0.6590) (1.855,0.6600) (1.884,0.6610) (1.895,0.6619)
    (1.909,0.6629) (1.942,0.6638) (1.978,0.6648) (2.014,0.6657)
    (2.016,0.6667) (2.055,0.6676) (2.103,0.6686) (2.115,0.6695)
    (2.115,0.6705) (2.131,0.6714) (2.142,0.6724) (2.182,0.6733)
    (2.227,0.6743) (2.258,0.6752) (2.284,0.6762) (2.37,0.6771)
    (2.57,0.6781) (2.601,0.6790) (2.637,0.6800) (2.708,0.6810)
    (2.787,0.6819) (2.808,0.6829) (2.827,0.6838) (2.949,0.6848)
    (2.997,0.6857) (3.193,0.6867) (3.259,0.6876) (3.269,0.6886)
    (3.32,0.6895) (3.518,0.6905) (3.576,0.6914) (3.882,0.6924)
    (4.811,0.6933) (4.921,0.6943) (4.941,0.6952) (5.619,0.6962)
    (6.016,0.6971) (211.9,0.6981) (10000,0.6981)
  };
  \addlegendentry{\Mone{}}
  \addplot+[const plot, no marks, color = btbottleneck, solid, ultra thick] coordinates {
    (1,0.0) (1,0.0010) (1,0.0019) (1,0.0029)
    (1,0.0038) (1,0.0048) (1,0.0057) (1,0.0067)
    (1,0.0076) (1,0.0086) (1,0.0095) (1,0.0105)
    (1,0.0114) (1,0.0124) (1,0.0133) (1,0.0143)
    (1,0.0152) (1,0.0162) (1,0.0171) (1,0.0181)
    (1,0.0190) (1,0.0200) (1,0.0210) (1,0.0219)
    (1,0.0229) (1,0.0238) (1,0.0248) (1,0.0257)
    (1,0.0267) (1,0.0276) (1,0.0286) (1,0.0295)
    (1,0.0305) (1,0.0314) (1,0.0324) (1,0.0333)
    (1,0.0343) (1,0.0352) (1,0.0362) (1,0.0371)
    (1,0.0381) (1,0.0390) (1,0.0400) (1,0.0410)
    (1,0.0419) (1,0.0429) (1,0.0438) (1,0.0448)
    (1,0.0457) (1,0.0467) (1,0.0476) (1,0.0486)
    (1,0.0495) (1,0.0505) (1,0.0514) (1,0.0524)
    (1,0.0533) (1,0.0543) (1,0.0552) (1,0.0562)
    (1,0.0571) (1,0.0581) (1,0.0590) (1,0.0600)
    (1,0.0610) (1,0.0619) (1,0.0629) (1,0.0638)
    (1,0.0648) (1,0.0657) (1,0.0667) (1,0.0676)
    (1,0.0686) (1,0.0695) (1,0.0705) (1,0.0714)
    (1,0.0724) (1,0.0733) (1,0.0743) (1,0.0752)
    (1,0.0762) (1,0.0771) (1,0.0781) (1,0.0790)
    (1,0.0800) (1,0.0810) (1,0.0819) (1,0.0829)
    (1,0.0838) (1,0.0848) (1,0.0857) (1,0.0867)
    (1,0.0876) (1,0.0886) (1,0.0895) (1,0.0905)
    (1,0.0914) (1,0.0924) (1,0.0933) (1.003,0.0943)
    (1.005,0.0952) (1.015,0.0962) (1.019,0.0971) (1.02,0.0981)
    (1.03,0.0990) (1.033,0.1000) (1.04,0.1010) (1.057,0.1019)
    (1.062,0.1029) (1.077,0.1038) (1.079,0.1048) (1.082,0.1057)
    (1.113,0.1067) (1.123,0.1076) (1.13,0.1086) (1.134,0.1095)
    (1.134,0.1105) (1.157,0.1114) (1.164,0.1124) (1.168,0.1133)
    (1.171,0.1143) (1.187,0.1152) (1.193,0.1162) (1.195,0.1171)
    (1.2,0.1181) (1.211,0.1190) (1.22,0.1200) (1.234,0.1210)
    (1.243,0.1219) (1.258,0.1229) (1.272,0.1238) (1.278,0.1248)
    (1.282,0.1257) (1.285,0.1267) (1.291,0.1276) (1.301,0.1286)
    (1.31,0.1295) (1.315,0.1305) (1.321,0.1314) (1.356,0.1324)
    (1.364,0.1333) (1.366,0.1343) (1.4,0.1352) (1.412,0.1362)
    (1.414,0.1371) (1.42,0.1381) (1.428,0.1390) (1.431,0.1400)
    (1.434,0.1410) (1.434,0.1419) (1.435,0.1429) (1.444,0.1438)
    (1.455,0.1448) (1.463,0.1457) (1.464,0.1467) (1.47,0.1476)
    (1.481,0.1486) (1.484,0.1495) (1.502,0.1505) (1.503,0.1514)
    (1.514,0.1524) (1.523,0.1533) (1.563,0.1543) (1.566,0.1552)
    (1.572,0.1562) (1.578,0.1571) (1.58,0.1581) (1.584,0.1590)
    (1.585,0.1600) (1.587,0.1610) (1.588,0.1619) (1.589,0.1629)
    (1.591,0.1638) (1.618,0.1648) (1.624,0.1657) (1.629,0.1667)
    (1.645,0.1676) (1.65,0.1686) (1.658,0.1695) (1.665,0.1705)
    (1.672,0.1714) (1.682,0.1724) (1.686,0.1733) (1.699,0.1743)
    (1.702,0.1752) (1.708,0.1762) (1.724,0.1771) (1.725,0.1781)
    (1.752,0.1790) (1.759,0.1800) (1.79,0.1810) (1.796,0.1819)
    (1.799,0.1829) (1.808,0.1838) (1.817,0.1848) (1.831,0.1857)
    (1.838,0.1867) (1.895,0.1876) (1.896,0.1886) (1.909,0.1895)
    (1.916,0.1905) (1.917,0.1914) (1.951,0.1924) (1.962,0.1933)
    (1.982,0.1943) (1.982,0.1952) (2.014,0.1962) (2.018,0.1971)
    (2.032,0.1981) (2.072,0.1990) (2.082,0.2000) (2.082,0.2010)
    (2.085,0.2019) (2.092,0.2029) (2.115,0.2038) (2.118,0.2048)
    (2.123,0.2057) (2.142,0.2067) (2.17,0.2076) (2.176,0.2086)
    (2.178,0.2095) (2.198,0.2105) (2.203,0.2114) (2.204,0.2124)
    (2.204,0.2133) (2.214,0.2143) (2.238,0.2152) (2.241,0.2162)
    (2.242,0.2171) (2.266,0.2181) (2.29,0.2190) (2.307,0.2200)
    (2.34,0.2210) (2.34,0.2219) (2.344,0.2229) (2.37,0.2238)
    (2.395,0.2248) (2.4,0.2257) (2.42,0.2267) (2.428,0.2276)
    (2.444,0.2286) (2.456,0.2295) (2.457,0.2305) (2.461,0.2314)
    (2.472,0.2324) (2.473,0.2333) (2.479,0.2343) (2.492,0.2352)
    (2.497,0.2362) (2.516,0.2371) (2.526,0.2381) (2.531,0.2390)
    (2.535,0.2400) (2.54,0.2410) (2.553,0.2419) (2.579,0.2429)
    (2.586,0.2438) (2.634,0.2448) (2.651,0.2457) (2.651,0.2467)
    (2.676,0.2476) (2.679,0.2486) (2.686,0.2495) (2.7,0.2505)
    (2.702,0.2514) (2.702,0.2524) (2.72,0.2533) (2.726,0.2543)
    (2.731,0.2552) (2.743,0.2562) (2.749,0.2571) (2.794,0.2581)
    (2.796,0.2590) (2.804,0.2600) (2.81,0.2610) (2.839,0.2619)
    (2.854,0.2629) (2.872,0.2638) (2.901,0.2648) (2.935,0.2657)
    (2.986,0.2667) (3.008,0.2676) (3.033,0.2686) (3.035,0.2695)
    (3.05,0.2705) (3.083,0.2714) (3.107,0.2724) (3.112,0.2733)
    (3.114,0.2743) (3.123,0.2752) (3.134,0.2762) (3.158,0.2771)
    (3.159,0.2781) (3.189,0.2790) (3.191,0.2800) (3.2,0.2810)
    (3.225,0.2819) (3.227,0.2829) (3.254,0.2838) (3.272,0.2848)
    (3.273,0.2857) (3.299,0.2867) (3.299,0.2876) (3.314,0.2886)
    (3.315,0.2895) (3.319,0.2905) (3.378,0.2914) (3.379,0.2924)
    (3.393,0.2933) (3.454,0.2943) (3.476,0.2952) (3.554,0.2962)
    (3.604,0.2971) (3.627,0.2981) (3.627,0.2990) (3.729,0.3000)
    (3.73,0.3010) (3.732,0.3019) (3.839,0.3029) (3.852,0.3038)
    (3.854,0.3048) (3.9,0.3057) (4.044,0.3067) (4.124,0.3076)
    (4.174,0.3086) (4.194,0.3095) (4.209,0.3105) (4.21,0.3114)
    (4.228,0.3124) (4.27,0.3133) (4.296,0.3143) (4.329,0.3152)
    (4.33,0.3162) (4.346,0.3171) (4.357,0.3181) (4.367,0.3190)
    (4.379,0.3200) (4.423,0.3210) (4.457,0.3219) (4.46,0.3229)
    (4.462,0.3238) (4.517,0.3248) (4.551,0.3257) (4.588,0.3267)
    (4.592,0.3276) (4.609,0.3286) (4.679,0.3295) (4.701,0.3305)
    (4.71,0.3314) (4.739,0.3324) (4.744,0.3333) (4.798,0.3343)
    (4.82,0.3352) (5.15,0.3362) (5.155,0.3371) (5.187,0.3381)
    (5.188,0.3390) (5.231,0.3400) (5.271,0.3410) (5.323,0.3419)
    (5.342,0.3429) (5.367,0.3438) (5.37,0.3448) (5.474,0.3457)
    (5.485,0.3467) (5.536,0.3476) (5.615,0.3486) (5.626,0.3495)
    (5.724,0.3505) (5.796,0.3514) (5.824,0.3524) (5.89,0.3533)
    (5.906,0.3543) (5.955,0.3552) (5.977,0.3562) (6.069,0.3571)
    (6.092,0.3581) (6.104,0.3590) (6.308,0.3600) (6.313,0.3610)
    (6.33,0.3619) (6.353,0.3629) (6.359,0.3638) (6.387,0.3648)
    (6.655,0.3657) (6.693,0.3667) (6.728,0.3676) (6.731,0.3686)
    (6.768,0.3695) (6.787,0.3705) (6.855,0.3714) (6.878,0.3724)
    (6.885,0.3733) (7.312,0.3743) (7.365,0.3752) (7.436,0.3762)
    (7.458,0.3771) (7.52,0.3781) (8.029,0.3790) (8.09,0.3800)
    (8.142,0.3810) (8.219,0.3819) (8.268,0.3829) (8.3,0.3838)
    (8.322,0.3848) (8.334,0.3857) (8.344,0.3867) (8.433,0.3876)
    (8.454,0.3886) (8.459,0.3895) (8.471,0.3905) (8.582,0.3914)
    (8.662,0.3924) (8.761,0.3933) (8.989,0.3943) (9.222,0.3952)
    (9.227,0.3962) (9.233,0.3971) (9.266,0.3981) (9.366,0.3990)
    (9.54,0.4000) (9.541,0.4010) (9.557,0.4019) (9.739,0.4029)
    (9.762,0.4038) (9.825,0.4048) (9.875,0.4057) (9.94,0.4067)
    (10,0.4076) (10.01,0.4086) (10.03,0.4095) (10.07,0.4105)
    (10.29,0.4114) (10.41,0.4124) (10.47,0.4133) (10.48,0.4143)
    (10.72,0.4152) (10.79,0.4162) (10.99,0.4171) (11.03,0.4181)
    (11.24,0.4190) (11.25,0.4200) (11.29,0.4210) (11.38,0.4219)
    (11.47,0.4229) (11.5,0.4238) (11.58,0.4248) (11.7,0.4257)
    (11.84,0.4267) (11.89,0.4276) (11.95,0.4286) (11.95,0.4295)
    (11.96,0.4305) (12.02,0.4314) (12.28,0.4324) (12.4,0.4333)
    (12.43,0.4343) (12.65,0.4352) (12.7,0.4362) (12.8,0.4371)
    (12.85,0.4381) (13.03,0.4390) (13.03,0.4400) (13.21,0.4410)
    (13.31,0.4419) (13.42,0.4429) (13.54,0.4438) (13.54,0.4448)
    (13.72,0.4457) (13.96,0.4467) (14.02,0.4476) (14.09,0.4486)
    (14.14,0.4495) (14.32,0.4505) (14.76,0.4514) (14.93,0.4524)
    (15.47,0.4533) (15.51,0.4543) (15.58,0.4552) (15.91,0.4562)
    (16.24,0.4571) (16.54,0.4581) (16.62,0.4590) (17.37,0.4600)
    (17.45,0.4610) (17.77,0.4619) (17.8,0.4629) (17.84,0.4638)
    (17.97,0.4648) (18.12,0.4657) (18.19,0.4667) (18.68,0.4676)
    (19.15,0.4686) (19.18,0.4695) (19.43,0.4705) (20.06,0.4714)
    (20.35,0.4724) (20.35,0.4733) (20.54,0.4743) (20.6,0.4752)
    (21.81,0.4762) (21.86,0.4771) (22.24,0.4781) (22.27,0.4790)
    (22.93,0.4800) (23.04,0.4810) (23.32,0.4819) (23.84,0.4829)
    (25.65,0.4838) (26.89,0.4848) (27.79,0.4857) (28.42,0.4867)
    (28.46,0.4876) (28.47,0.4886) (28.87,0.4895) (29.12,0.4905)
    (29.12,0.4914) (29.74,0.4924) (30.33,0.4933) (30.51,0.4943)
    (30.73,0.4952) (31.2,0.4962) (31.37,0.4971) (32.03,0.4981)
    (32.53,0.4990) (33.68,0.5000) (34.01,0.5010) (35.84,0.5019)
    (38.16,0.5029) (38.6,0.5038) (38.94,0.5048) (39.13,0.5057)
    (39.23,0.5067) (39.55,0.5076) (40.18,0.5086) (40.35,0.5095)
    (41.21,0.5105) (41.32,0.5114) (42.29,0.5124) (42.44,0.5133)
    (42.61,0.5143) (43.24,0.5152) (43.82,0.5162) (47.53,0.5171)
    (47.7,0.5181) (47.93,0.5190) (48.33,0.5200) (48.64,0.5210)
    (51.08,0.5219) (52.16,0.5229) (60.33,0.5238) (63.48,0.5248)
    (64.57,0.5257) (65.05,0.5267) (70.62,0.5276) (70.69,0.5286)
    (71.88,0.5295) (75.77,0.5305) (82.63,0.5314) (87.67,0.5324)
    (94.8,0.5333) (94.85,0.5343) (181.6,0.5352) (246.6,0.5362)
    (269.8,0.5371) (5527,0.5381) (9722,0.5390) (10000,0.5390)
  };
  \addlegendentry{\Mtwo{}}
\end{semilogxaxis}
\end{tikzpicture}
   
  \caption{Performance profile~\citep*{dolan2002} for \Mone{} and \Mtwo{}, both with the MST restriction, on the
    $1050$ instances of test set~C. For each instance $i\in\mathcal{I}$ and method $s$,
    let $r_{i,s}:=t_{i,s}/\min_{s'}t_{i,s'}$ be the ratio of the solve
    time of $s$ to the solve time of the faster method on $i$, with $r_{i,s}=+\infty$
    when $s$ fails to solve $i$ within the time limit. The performance
    profile $\rho_s(\tau):=\bigl|\{i\in\mathcal{I}\colon r_{i,s}\le\tau\}\bigr|/|\mathcal{I}|$, with $|\mathcal{I}|=1050$, is
    the fraction of instances on which $s$ is within a factor of $\tau$ of
    the faster method (horizontal axis on a logarithmic scale). \Mone{} is
    the blue curve and \Mtwo{} the red curve. Failures
  contribute zero to $\rho_s(\tau)$ for every finite~$\tau$.}
  \label{fig:perfprof}
\end{figure}
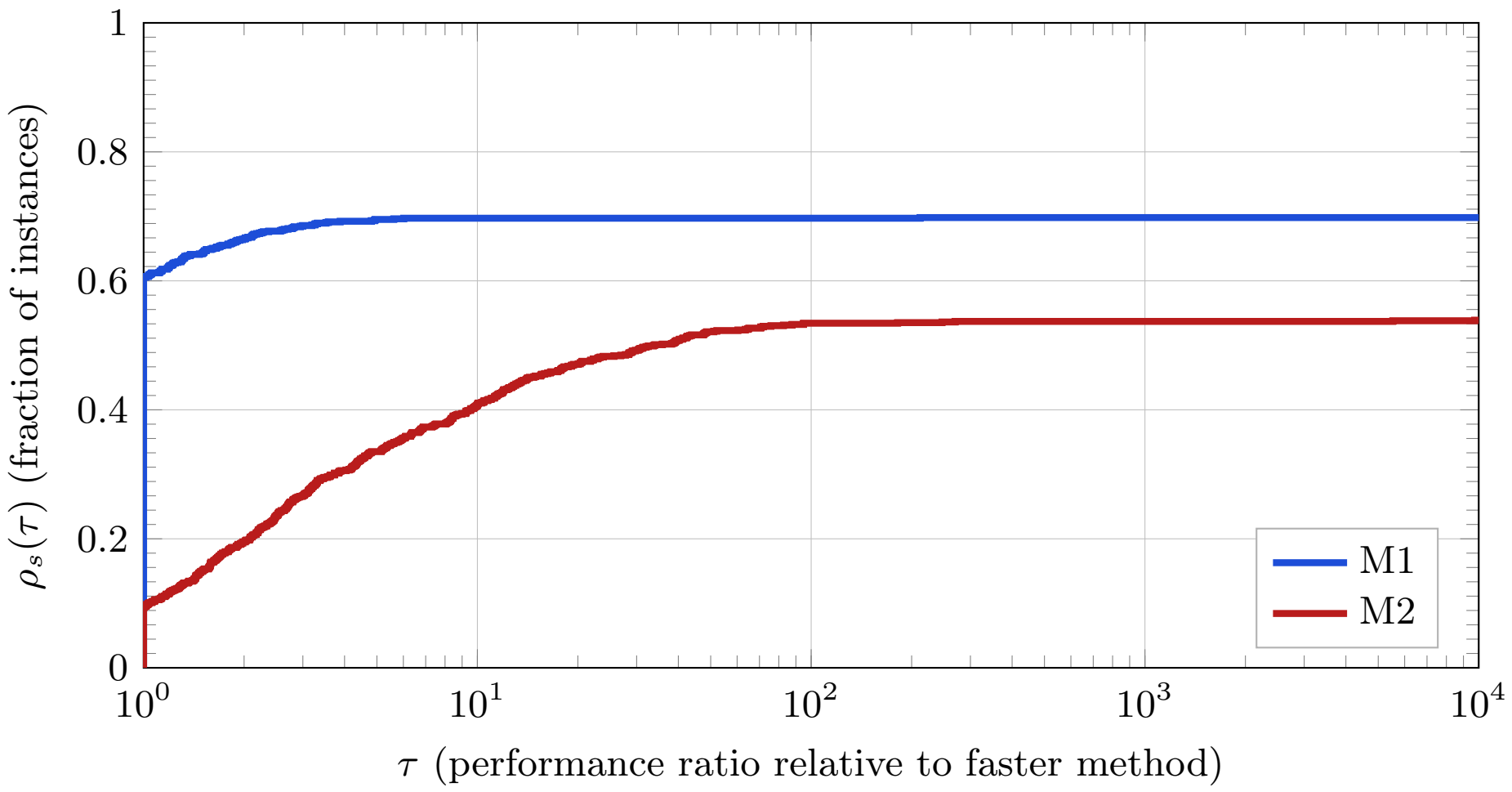

The performance profile of Figure~\ref{fig:perfprof} { summarizes Table~\ref{tab:main}.} { The intercept at $\tau=1$ gives the fraction of instances on
which each method is faster: $\rho_{\Mone{}}(1)=0.606$ and $\rho_{\Mtwo{}}(1)=0.093$ (the
two values do not sum to~$1$ because $316$ of the $1050$ instances are not solved by either method
within the time limit). { These intercepts use the raw solve times, so each of the $94$ instances whose two solve times
differ by less than $0.05$\,s counts for the faster method, $51$ for {\rm \Idenone} and $43$ for
{\rm \Identity}; the summary below counts them as ties.}
At $\tau=10$, $\rho_{\Mone{}}=0.70$ against $\rho_{\Mtwo{}}=0.41$, and
the curves saturate at $0.70$ and $0.54$, so the gap between the formulations is one of both
solve time and solve rate, with {\rm \Idenone} solving $16$ percentage points more of the
benchmark.}

{ A paired McNemar test on the $1050$ instances confirms that the difference in solve
rates is significant: {\rm \Idenone} solves $168$ instances that {\rm \Identity} does not, while
the converse happens on a single instance ($p<10^{-40}$). The difference grows with the dimension
and the budget: it is $23$ against $0$ at $(d,k)=(3,5)$, $36$ against $0$ at $(3,7)$ and $49$ against $0$
at $(4,5)$.}

\paragraph{Summary} { {\rm \Idenone} is faster than {\rm \Identity} in median solve time in every
dimension, is the faster method on $585$ instances against $55$ (solve times within the $0.05$\,s
reporting resolution count as ties, on $94$ instances, and $316$ instances are solved by neither),
and solves $168$ instances that {\rm \Identity} does not, against $1$, with an advantage that grows with
the dimension and the Steiner budget. {\rm \Idenone} is therefore the recommended exact formulation
across the test sets.}

\paragraph{Effect of the MST restriction} { On the $360$ instances of test set~C with
$p\le25$ and $k\le5$, restricting the terminal--terminal edges to $H$
(Proposition~\ref{prop:mst}) removes a median of $47\%$ of the terminal--terminal edges that
survive the elimination~\eqref{eq:elim} and changes no optimal value on the instances solved with and without it. For {\rm \Idenone} it
solves $358$ instances against $348$, losing none, reduces the total solve time by $55\%$, and has a
median speedup of $1.40$ on the instances both versions solve. For {\rm \Identity} it barely changes
the number of solved instances ($9$ gained and $8$ lost). Table~\ref{tab:msteffect} reports
{\rm \Idenone} with and without the restriction. Without the restriction for either formulation,
the ratios of median solve times of {\rm \Identity} to {\rm \Idenone} are $1.3$, $2.7$ and
$13.3$ at $d=2,3,4$, and {\rm \Idenone} solves $35$ instances that {\rm \Identity} does not,
against $3$.}
On 270 larger instances with $p\in\{30,40,50\}$ and $k\le5$, the restriction gains {\rm\Idenone} 36 solved instances, loses none, and reduces its total solve time
by 43\%.

\begin{table}[!ht]
  \centering
   
  \caption{\Mone{} versus \Mone{} with the MST restriction on the $360$ instances of test set~C with $p\le25$ and $k\le5$, aggregated by ambient dimension, by dimension and instance family, and by dimension and Steiner budget~$k$. ``solv.''/``\%'' count proven-optimal solves within \num{1000}\,s; $t_{\med}$/$t_{\sgm}$ are the median and shifted-geometric-mean solve times (s, shift~$10$); ``nodes'' the median branch-and-bound node count over solved instances; ``gap$^*$'' the worst gap over unsolved instances; ``all'' aggregates the omitted parameters. Bold marks the per-row best; node counts are bold only on equal solved sets.}
  \label{tab:msteffect}
  \footnotesize
  \setlength{\tabcolsep}{2pt}
  \renewcommand*{\arraystretch}{1.3}
  \begin{tabular*}{\textwidth}{@{\extracolsep{\fill}}
      c l c S[table-format=3.0]
      *2{S[table-format=3.0] S[table-format=3.1] S[table-format=4.1] S[table-format=3.1]
    S[table-format=5.0] S[table-format=2.2]}}
    \toprule
    & & & & \multicolumn{6}{c}{\Mone{}} & \multicolumn{6}{c}{\Mone{} with MST restriction}\\
    \cmidrule(lr){5-10}\cmidrule(lr){11-16}
    {$d$} & {family} & {$k$} & {\#inst} & {solv.} & {\%} & {$t_{\med}$} & {$t_{\sgm}$} & {nodes} & {gap$^*$}
    & {solv.} & {\%} & {$t_{\med}$} & {$t_{\sgm}$} & {nodes} & {gap$^*$}\\
    \midrule
\multicolumn{16}{@{}l}{\itshape By ambient dimension}\\
2 & {all} & {all} & 120 & 118 & 98.3 & 1.7 & 6.5 & 1046 & 18.00 & \best{119} & \best{99.2} & \best{1.2} & \best{4.9} & 1011 & \best{1.14}\\
3 & {all} & {all} & 120 & 118 & 98.3 & 5.2 & 10.0 & 1786 & 21.30 & \best{120} & \best{100.0} & \best{2.8} & \best{6.1} & 1534 & \best{0}\\
4 & {all} & {all} & 120 & 112 & 93.3 & 6.9 & 19.9 & 3602 & 54.25 & \best{119} & \best{99.2} & \best{6.5} & \best{13.5} & 3083 & \best{28.80}\\
\midrule
\multicolumn{16}{@{}l}{\itshape By dimension and instance family}\\
2 & clustered & {all} & 60 & 58 & 96.7 & 1.7 & 10.0 & 1452 & 18.00 & \best{59} & \best{98.3} & \best{1.3} & \best{6.7} & 1355 & \best{1.14}\\
2 & uniform & {all} & 60 & 60 & 100.0 & 1.7 & 3.6 & 965 & 0 & 60 & 100.0 & \best{1.1} & \best{3.3} & \best{760} & 0\\
3 & clustered & {all} & 60 & 58 & 96.7 & 4.6 & 15.4 & 2896 & 21.30 & \best{60} & \best{100.0} & \best{2.2} & \best{8.1} & 1940 & \best{0}\\
3 & uniform & {all} & 60 & 60 & 100.0 & 5.3 & 5.7 & 1116 & 0 & 60 & 100.0 & \best{3.6} & \best{4.3} & \best{1090} & 0\\
4 & clustered & {all} & 60 & 52 & 86.7 & 7.6 & 32.9 & 4216 & 54.25 & \best{59} & \best{98.3} & \best{6.1} & \best{19.0} & 5350 & \best{28.80}\\
4 & uniform & {all} & 60 & 60 & 100.0 & 6.9 & 10.8 & 3244 & 0 & 60 & 100.0 & \best{6.7} & \best{9.0} & \best{2896} & 0\\
\midrule
\multicolumn{16}{@{}l}{\itshape By dimension and Steiner budget}\\
2 & {all} & 1 & 40 & 40 & 100.0 & {$<0.05$} & 0.1 & 1 & 0 & 40 & 100.0 & {$<0.05$} & {\bestmath{<0.05}} & 1 & 0\\
2 & {all} & 3 & 40 & 40 & 100.0 & 3.4 & 4.6 & 1194 & 0 & 40 & 100.0 & \best{1.3} & \best{2.3} & \best{1094} & 0\\
2 & {all} & 5 & 40 & 38 & 95.0 & \best{9.7} & 20.4 & 4138 & 18.00 & \best{39} & \best{97.5} & 10.5 & \best{16.9} & 3715 & \best{1.14}\\
3 & {all} & 1 & 40 & 40 & 100.0 & 0.1 & 0.1 & 12 & 0 & 40 & 100.0 & {\bestmath{<0.05}} & {\bestmath{<0.05}} & \best{1} & 0\\
3 & {all} & 3 & 40 & 40 & 100.0 & 5.6 & 6.2 & 2438 & 0 & 40 & 100.0 & \best{2.8} & \best{3.6} & \best{1534} & 0\\
3 & {all} & 5 & 40 & 38 & 95.0 & 19.3 & 38.8 & 8847 & 21.30 & \best{40} & \best{100.0} & \best{13.7} & \best{20.5} & 8445 & \best{0}\\
4 & {all} & 1 & 40 & 40 & 100.0 & 0.2 & 0.3 & 627 & 0 & 40 & 100.0 & \best{0.1} & \best{0.1} & \best{290} & 0\\
4 & {all} & 3 & 40 & 38 & 95.0 & 6.9 & 15.1 & 3794 & 12.46 & \best{40} & \best{100.0} & \best{6.5} & \best{7.8} & 3119 & \best{0}\\
4 & {all} & 5 & 40 & 34 & 85.0 & 61.8 & 92.8 & 18458 & 54.25 & \best{39} & \best{97.5} & \best{45.2} & \best{62.3} & 21727 & \best{28.80}\\
    \bottomrule
  \end{tabular*}
\end{table}

\subsection{Structured named instances}
\label{sec:namedexp}

\begin{figure}[htbp]
  \centering
  \begin{minipage}[t]{0.32\linewidth}\centering
\begin{tikzpicture}
\begin{axis}
  [
    steineraxis3d,
    width = 4.6000cm,
    height = 4.6000cm,
    xmin = -1.0400, xmax = 1.0400,
    ymin = -1.0400, ymax = 1.0400,
    zmin = -1.0400, zmax = 1.0400,
  ]
  \addplot3[btwire, no marks, forget plot] coordinates {(1.0000,0.0000,0.0000) (0.0000,1.0000,0.0000)};
  \addplot3[btwire, no marks, forget plot] coordinates {(1.0000,0.0000,0.0000) (0.0000,-1.0000,0.0000)};
  \addplot3[btwire, no marks, forget plot] coordinates {(1.0000,0.0000,0.0000) (0.0000,0.0000,1.0000)};
  \addplot3[btwire, no marks, forget plot] coordinates {(1.0000,0.0000,0.0000) (0.0000,0.0000,-1.0000)};
  \addplot3[btwire, no marks, forget plot] coordinates {(-1.0000,0.0000,0.0000) (0.0000,1.0000,0.0000)};
  \addplot3[btwire, no marks, forget plot] coordinates {(-1.0000,0.0000,0.0000) (0.0000,-1.0000,0.0000)};
  \addplot3[btwire, no marks, forget plot] coordinates {(-1.0000,0.0000,0.0000) (0.0000,0.0000,1.0000)};
  \addplot3[btwire, no marks, forget plot] coordinates {(-1.0000,0.0000,0.0000) (0.0000,0.0000,-1.0000)};
  \addplot3[btwire, no marks, forget plot] coordinates {(0.0000,1.0000,0.0000) (0.0000,0.0000,1.0000)};
  \addplot3[btwire, no marks, forget plot] coordinates {(0.0000,1.0000,0.0000) (0.0000,0.0000,-1.0000)};
  \addplot3[btwire, no marks, forget plot] coordinates {(0.0000,-1.0000,0.0000) (0.0000,0.0000,1.0000)};
  \addplot3[btwire, no marks, forget plot] coordinates {(0.0000,-1.0000,0.0000) (0.0000,0.0000,-1.0000)};
  \addplot3[btbottleneck3d, no marks, forget plot] coordinates {(-1.0000,0.0000,0.0000) (-0.0000,0.0000,0.0000)};
  \addplot3[btbottleneck3d, no marks, forget plot] coordinates {(0.0000,-1.0000,0.0000) (-0.0000,0.0000,0.0000)};
  \addplot3[btbottleneck3d, no marks, forget plot] coordinates {(0.0000,0.0000,-1.0000) (-0.0000,0.0000,0.0000)};
  \addplot3[btbottleneck3d, no marks, forget plot] coordinates {(0.0000,0.0000,1.0000) (-0.0000,0.0000,0.0000)};
  \addplot3[btbottleneck3d, no marks, forget plot] coordinates {(0.0000,1.0000,0.0000) (-0.0000,0.0000,0.0000)};
  \addplot3[btbottleneck3d, no marks, forget plot] coordinates {(1.0000,0.0000,0.0000) (-0.0000,0.0000,0.0000)};
  \addplot3[only marks, mark = *, mark size = 2.4pt, mark options = {fill = btterminal, draw = btterminal}, forget plot] coordinates {(1.0000,0.0000,0.0000) (-1.0000,0.0000,0.0000) (0.0000,1.0000,0.0000) (0.0000,-1.0000,0.0000) (0.0000,0.0000,1.0000) (0.0000,0.0000,-1.0000) };
  \addplot3[only marks, mark = triangle*, mark size = 3.4pt, mark options = {fill = btsteiner, draw = btsteiner}, forget plot] coordinates {(-0.0000,0.0000,0.0000) };
\end{axis}
\end{tikzpicture}\\[2pt]
    {\small (a) \texttt{octahedron}, $k=1$}
  \end{minipage}\hfill
  \begin{minipage}[t]{0.32\linewidth}\centering
\begin{tikzpicture}
\begin{axis}
  [
    steineraxis3d,
    width = 4.6000cm,
    height = 4.6000cm,
    xmin = -1.0400, xmax = 1.0400,
    ymin = -1.0400, ymax = 1.0400,
    zmin = -1.0400, zmax = 1.0400,
  ]
  \addplot3[btwire, no marks, forget plot] coordinates {(1.0000,0.0000,0.0000) (0.0000,1.0000,0.0000)};
  \addplot3[btwire, no marks, forget plot] coordinates {(1.0000,0.0000,0.0000) (0.0000,-1.0000,0.0000)};
  \addplot3[btwire, no marks, forget plot] coordinates {(1.0000,0.0000,0.0000) (0.0000,0.0000,1.0000)};
  \addplot3[btwire, no marks, forget plot] coordinates {(1.0000,0.0000,0.0000) (0.0000,0.0000,-1.0000)};
  \addplot3[btwire, no marks, forget plot] coordinates {(-1.0000,0.0000,0.0000) (0.0000,1.0000,0.0000)};
  \addplot3[btwire, no marks, forget plot] coordinates {(-1.0000,0.0000,0.0000) (0.0000,-1.0000,0.0000)};
  \addplot3[btwire, no marks, forget plot] coordinates {(-1.0000,0.0000,0.0000) (0.0000,0.0000,1.0000)};
  \addplot3[btwire, no marks, forget plot] coordinates {(-1.0000,0.0000,0.0000) (0.0000,0.0000,-1.0000)};
  \addplot3[btwire, no marks, forget plot] coordinates {(0.0000,1.0000,0.0000) (0.0000,0.0000,1.0000)};
  \addplot3[btwire, no marks, forget plot] coordinates {(0.0000,1.0000,0.0000) (0.0000,0.0000,-1.0000)};
  \addplot3[btwire, no marks, forget plot] coordinates {(0.0000,-1.0000,0.0000) (0.0000,0.0000,1.0000)};
  \addplot3[btwire, no marks, forget plot] coordinates {(0.0000,-1.0000,0.0000) (0.0000,0.0000,-1.0000)};
  \addplot3[btedge3d, no marks, forget plot] coordinates {(-1.0000,0.0000,0.0000) (-0.2403,-0.2403,-0.2403)};
  \addplot3[btedge3d, no marks, forget plot] coordinates {(0.0000,-1.0000,0.0000) (-0.2403,-0.2403,-0.2403)};
  \addplot3[btbottleneck3d, no marks, forget plot] coordinates {(0.0000,0.0000,-1.0000) (-0.2403,-0.2403,-0.2403)};
  \addplot3[btbottleneck3d, no marks, forget plot] coordinates {(-0.2403,-0.2403,-0.2403) (0.2403,0.2403,0.2403)};
  \addplot3[btedge3d, no marks, forget plot] coordinates {(0.0000,0.0000,1.0000) (0.2403,0.2403,0.2403)};
  \addplot3[btedge3d, no marks, forget plot] coordinates {(0.0000,1.0000,0.0000) (0.2403,0.2403,0.2403)};
  \addplot3[btedge3d, no marks, forget plot] coordinates {(1.0000,0.0000,0.0000) (0.2403,0.2403,0.2403)};
  \addplot3[only marks, mark = *, mark size = 2.4pt, mark options = {fill = btterminal, draw = btterminal}, forget plot] coordinates {(1.0000,0.0000,0.0000) (-1.0000,0.0000,0.0000) (0.0000,1.0000,0.0000) (0.0000,-1.0000,0.0000) (0.0000,0.0000,1.0000) (0.0000,0.0000,-1.0000) };
  \addplot3[only marks, mark = triangle*, mark size = 3.4pt, mark options = {fill = btsteiner, draw = btsteiner}, forget plot] coordinates {(-0.2403,-0.2403,-0.2403) (0.2403,0.2403,0.2403) };
\end{axis}
\end{tikzpicture}\\[2pt]
    {\small (b) \texttt{octahedron}, $k=2$}
  \end{minipage}\hfill
  \begin{minipage}[t]{0.32\linewidth}\centering
\begin{tikzpicture}
\begin{axis}
  [
    steineraxis3d,
    width = 4.6000cm,
    height = 4.6000cm,
    xmin = -1.0400, xmax = 1.0400,
    ymin = -1.0400, ymax = 1.0400,
    zmin = -1.0400, zmax = 1.0400,
  ]
  \addplot3[btwire, no marks, forget plot] coordinates {(1.0000,0.0000,0.0000) (0.0000,1.0000,0.0000)};
  \addplot3[btwire, no marks, forget plot] coordinates {(1.0000,0.0000,0.0000) (0.0000,-1.0000,0.0000)};
  \addplot3[btwire, no marks, forget plot] coordinates {(1.0000,0.0000,0.0000) (0.0000,0.0000,1.0000)};
  \addplot3[btwire, no marks, forget plot] coordinates {(1.0000,0.0000,0.0000) (0.0000,0.0000,-1.0000)};
  \addplot3[btwire, no marks, forget plot] coordinates {(-1.0000,0.0000,0.0000) (0.0000,1.0000,0.0000)};
  \addplot3[btwire, no marks, forget plot] coordinates {(-1.0000,0.0000,0.0000) (0.0000,-1.0000,0.0000)};
  \addplot3[btwire, no marks, forget plot] coordinates {(-1.0000,0.0000,0.0000) (0.0000,0.0000,1.0000)};
  \addplot3[btwire, no marks, forget plot] coordinates {(-1.0000,0.0000,0.0000) (0.0000,0.0000,-1.0000)};
  \addplot3[btwire, no marks, forget plot] coordinates {(0.0000,1.0000,0.0000) (0.0000,0.0000,1.0000)};
  \addplot3[btwire, no marks, forget plot] coordinates {(0.0000,1.0000,0.0000) (0.0000,0.0000,-1.0000)};
  \addplot3[btwire, no marks, forget plot] coordinates {(0.0000,-1.0000,0.0000) (0.0000,0.0000,1.0000)};
  \addplot3[btwire, no marks, forget plot] coordinates {(0.0000,-1.0000,0.0000) (0.0000,0.0000,-1.0000)};
  \addplot3[btbottleneck3d, no marks, forget plot] coordinates {(-1.0000,0.0000,0.0000) (-0.3248,0.3248,0.0839)};
  \addplot3[btbottleneck3d, no marks, forget plot] coordinates {(0.0000,-1.0000,0.0000) (-0.0000,-0.3151,0.3151)};
  \addplot3[btbottleneck3d, no marks, forget plot] coordinates {(0.0000,0.0000,-1.0000) (0.3248,-0.0839,-0.3248)};
  \addplot3[btbottleneck3d, no marks, forget plot] coordinates {(-0.3248,0.3248,0.0839) (-0.0000,-0.3151,0.3151)};
  \addplot3[btbottleneck3d, no marks, forget plot] coordinates {(-0.0000,-0.3151,0.3151) (0.3248,-0.0839,-0.3248)};
  \addplot3[btbottleneck3d, no marks, forget plot] coordinates {(0.0000,0.0000,1.0000) (-0.0000,-0.3151,0.3151)};
  \addplot3[btbottleneck3d, no marks, forget plot] coordinates {(0.0000,1.0000,0.0000) (-0.3248,0.3248,0.0839)};
  \addplot3[btbottleneck3d, no marks, forget plot] coordinates {(1.0000,0.0000,0.0000) (0.3248,-0.0839,-0.3248)};
  \addplot3[only marks, mark = *, mark size = 2.4pt, mark options = {fill = btterminal, draw = btterminal}, forget plot] coordinates {(1.0000,0.0000,0.0000) (-1.0000,0.0000,0.0000) (0.0000,1.0000,0.0000) (0.0000,-1.0000,0.0000) (0.0000,0.0000,1.0000) (0.0000,0.0000,-1.0000) };
  \addplot3[only marks, mark = triangle*, mark size = 3.4pt, mark options = {fill = btsteiner, draw = btsteiner}, forget plot] coordinates {(-0.3248,0.3248,0.0839) (-0.0000,-0.3151,0.3151) (0.3248,-0.0839,-0.3248) };
\end{axis}
\end{tikzpicture}\\[2pt]
    {\small (c) \texttt{octahedron}, $k=3$}
  \end{minipage}

  \vspace{1em}
  \begin{minipage}[t]{0.32\linewidth}\centering
\begin{tikzpicture}
\begin{axis}
  [
    steineraxis3d,
    width = 4.6000cm,
    height = 4.6000cm,
    xmin = -1.6828, xmax = 1.6828,
    ymin = -1.6828, ymax = 1.6828,
    zmin = -1.6828, zmax = 1.6828,
  ]
  \addplot3[btwire, no marks, forget plot] coordinates {(0.0000,1.0000,1.6180) (0.0000,-1.0000,1.6180)};
  \addplot3[btwire, no marks, forget plot] coordinates {(0.0000,1.0000,1.6180) (1.0000,1.6180,0.0000)};
  \addplot3[btwire, no marks, forget plot] coordinates {(0.0000,1.0000,1.6180) (-1.0000,1.6180,0.0000)};
  \addplot3[btwire, no marks, forget plot] coordinates {(0.0000,1.0000,1.6180) (1.6180,0.0000,1.0000)};
  \addplot3[btwire, no marks, forget plot] coordinates {(0.0000,1.0000,1.6180) (-1.6180,0.0000,1.0000)};
  \addplot3[btwire, no marks, forget plot] coordinates {(0.0000,1.0000,-1.6180) (0.0000,-1.0000,-1.6180)};
  \addplot3[btwire, no marks, forget plot] coordinates {(0.0000,1.0000,-1.6180) (1.0000,1.6180,0.0000)};
  \addplot3[btwire, no marks, forget plot] coordinates {(0.0000,1.0000,-1.6180) (-1.0000,1.6180,0.0000)};
  \addplot3[btwire, no marks, forget plot] coordinates {(0.0000,1.0000,-1.6180) (1.6180,0.0000,-1.0000)};
  \addplot3[btwire, no marks, forget plot] coordinates {(0.0000,1.0000,-1.6180) (-1.6180,0.0000,-1.0000)};
  \addplot3[btwire, no marks, forget plot] coordinates {(0.0000,-1.0000,1.6180) (1.0000,-1.6180,0.0000)};
  \addplot3[btwire, no marks, forget plot] coordinates {(0.0000,-1.0000,1.6180) (-1.0000,-1.6180,0.0000)};
  \addplot3[btwire, no marks, forget plot] coordinates {(0.0000,-1.0000,1.6180) (1.6180,0.0000,1.0000)};
  \addplot3[btwire, no marks, forget plot] coordinates {(0.0000,-1.0000,1.6180) (-1.6180,0.0000,1.0000)};
  \addplot3[btwire, no marks, forget plot] coordinates {(0.0000,-1.0000,-1.6180) (1.0000,-1.6180,0.0000)};
  \addplot3[btwire, no marks, forget plot] coordinates {(0.0000,-1.0000,-1.6180) (-1.0000,-1.6180,0.0000)};
  \addplot3[btwire, no marks, forget plot] coordinates {(0.0000,-1.0000,-1.6180) (1.6180,0.0000,-1.0000)};
  \addplot3[btwire, no marks, forget plot] coordinates {(0.0000,-1.0000,-1.6180) (-1.6180,0.0000,-1.0000)};
  \addplot3[btwire, no marks, forget plot] coordinates {(1.0000,1.6180,0.0000) (-1.0000,1.6180,0.0000)};
  \addplot3[btwire, no marks, forget plot] coordinates {(1.0000,1.6180,0.0000) (1.6180,0.0000,1.0000)};
  \addplot3[btwire, no marks, forget plot] coordinates {(1.0000,1.6180,0.0000) (1.6180,0.0000,-1.0000)};
  \addplot3[btwire, no marks, forget plot] coordinates {(1.0000,-1.6180,0.0000) (-1.0000,-1.6180,0.0000)};
  \addplot3[btwire, no marks, forget plot] coordinates {(1.0000,-1.6180,0.0000) (1.6180,0.0000,1.0000)};
  \addplot3[btwire, no marks, forget plot] coordinates {(1.0000,-1.6180,0.0000) (1.6180,0.0000,-1.0000)};
  \addplot3[btwire, no marks, forget plot] coordinates {(-1.0000,1.6180,0.0000) (-1.6180,0.0000,1.0000)};
  \addplot3[btwire, no marks, forget plot] coordinates {(-1.0000,1.6180,0.0000) (-1.6180,0.0000,-1.0000)};
  \addplot3[btwire, no marks, forget plot] coordinates {(-1.0000,-1.6180,0.0000) (-1.6180,0.0000,1.0000)};
  \addplot3[btwire, no marks, forget plot] coordinates {(-1.0000,-1.6180,0.0000) (-1.6180,0.0000,-1.0000)};
  \addplot3[btwire, no marks, forget plot] coordinates {(1.6180,0.0000,1.0000) (1.6180,0.0000,-1.0000)};
  \addplot3[btwire, no marks, forget plot] coordinates {(-1.6180,0.0000,1.0000) (-1.6180,0.0000,-1.0000)};
  \addplot3[btbottleneck3d, no marks, forget plot] coordinates {(-1.6180,0.0000,-1.0000) (0.0000,0.0000,0.0000)};
  \addplot3[btbottleneck3d, no marks, forget plot] coordinates {(-1.0000,-1.6180,0.0000) (0.0000,0.0000,0.0000)};
  \addplot3[btbottleneck3d, no marks, forget plot] coordinates {(0.0000,-1.0000,-1.6180) (0.0000,0.0000,0.0000)};
  \addplot3[btbottleneck3d, no marks, forget plot] coordinates {(-1.6180,0.0000,1.0000) (0.0000,0.0000,0.0000)};
  \addplot3[btbottleneck3d, no marks, forget plot] coordinates {(0.0000,1.0000,-1.6180) (0.0000,0.0000,0.0000)};
  \addplot3[btbottleneck3d, no marks, forget plot] coordinates {(-1.0000,1.6180,0.0000) (0.0000,0.0000,0.0000)};
  \addplot3[btbottleneck3d, no marks, forget plot] coordinates {(1.0000,-1.6180,0.0000) (0.0000,0.0000,0.0000)};
  \addplot3[btbottleneck3d, no marks, forget plot] coordinates {(0.0000,-1.0000,1.6180) (0.0000,0.0000,0.0000)};
  \addplot3[btbottleneck3d, no marks, forget plot] coordinates {(1.6180,0.0000,-1.0000) (0.0000,0.0000,0.0000)};
  \addplot3[btbottleneck3d, no marks, forget plot] coordinates {(0.0000,1.0000,1.6180) (0.0000,0.0000,0.0000)};
  \addplot3[btbottleneck3d, no marks, forget plot] coordinates {(1.0000,1.6180,0.0000) (0.0000,0.0000,0.0000)};
  \addplot3[btbottleneck3d, no marks, forget plot] coordinates {(1.6180,0.0000,1.0000) (0.0000,0.0000,0.0000)};
  \addplot3[only marks, mark = *, mark size = 2.4pt, mark options = {fill = btterminal, draw = btterminal}, forget plot] coordinates {(0.0000,1.0000,1.6180) (0.0000,1.0000,-1.6180) (0.0000,-1.0000,1.6180) (0.0000,-1.0000,-1.6180) (1.0000,1.6180,0.0000) (1.0000,-1.6180,0.0000) (-1.0000,1.6180,0.0000) (-1.0000,-1.6180,0.0000) (1.6180,0.0000,1.0000) (1.6180,0.0000,-1.0000) (-1.6180,0.0000,1.0000) (-1.6180,0.0000,-1.0000) };
  \addplot3[only marks, mark = triangle*, mark size = 3.4pt, mark options = {fill = btsteiner, draw = btsteiner}, forget plot] coordinates {(0.0000,0.0000,0.0000) };
\end{axis}
\end{tikzpicture}\\[2pt]
    {\small (d) \texttt{icosahedron}, $k=1$}
  \end{minipage}\hfill
  \begin{minipage}[t]{0.32\linewidth}\centering
\begin{tikzpicture}
\begin{axis}
  [
    steineraxis3d,
    width = 4.6000cm,
    height = 4.6000cm,
    xmin = -1.6828, xmax = 1.6828,
    ymin = -1.6828, ymax = 1.6828,
    zmin = -1.6828, zmax = 1.6828,
  ]
  \addplot3[btwire, no marks, forget plot] coordinates {(0.0000,1.0000,1.6180) (0.0000,-1.0000,1.6180)};
  \addplot3[btwire, no marks, forget plot] coordinates {(0.0000,1.0000,1.6180) (1.0000,1.6180,0.0000)};
  \addplot3[btwire, no marks, forget plot] coordinates {(0.0000,1.0000,1.6180) (-1.0000,1.6180,0.0000)};
  \addplot3[btwire, no marks, forget plot] coordinates {(0.0000,1.0000,1.6180) (1.6180,0.0000,1.0000)};
  \addplot3[btwire, no marks, forget plot] coordinates {(0.0000,1.0000,1.6180) (-1.6180,0.0000,1.0000)};
  \addplot3[btwire, no marks, forget plot] coordinates {(0.0000,1.0000,-1.6180) (0.0000,-1.0000,-1.6180)};
  \addplot3[btwire, no marks, forget plot] coordinates {(0.0000,1.0000,-1.6180) (1.0000,1.6180,0.0000)};
  \addplot3[btwire, no marks, forget plot] coordinates {(0.0000,1.0000,-1.6180) (-1.0000,1.6180,0.0000)};
  \addplot3[btwire, no marks, forget plot] coordinates {(0.0000,1.0000,-1.6180) (1.6180,0.0000,-1.0000)};
  \addplot3[btwire, no marks, forget plot] coordinates {(0.0000,1.0000,-1.6180) (-1.6180,0.0000,-1.0000)};
  \addplot3[btwire, no marks, forget plot] coordinates {(0.0000,-1.0000,1.6180) (1.0000,-1.6180,0.0000)};
  \addplot3[btwire, no marks, forget plot] coordinates {(0.0000,-1.0000,1.6180) (-1.0000,-1.6180,0.0000)};
  \addplot3[btwire, no marks, forget plot] coordinates {(0.0000,-1.0000,1.6180) (1.6180,0.0000,1.0000)};
  \addplot3[btwire, no marks, forget plot] coordinates {(0.0000,-1.0000,1.6180) (-1.6180,0.0000,1.0000)};
  \addplot3[btwire, no marks, forget plot] coordinates {(0.0000,-1.0000,-1.6180) (1.0000,-1.6180,0.0000)};
  \addplot3[btwire, no marks, forget plot] coordinates {(0.0000,-1.0000,-1.6180) (-1.0000,-1.6180,0.0000)};
  \addplot3[btwire, no marks, forget plot] coordinates {(0.0000,-1.0000,-1.6180) (1.6180,0.0000,-1.0000)};
  \addplot3[btwire, no marks, forget plot] coordinates {(0.0000,-1.0000,-1.6180) (-1.6180,0.0000,-1.0000)};
  \addplot3[btwire, no marks, forget plot] coordinates {(1.0000,1.6180,0.0000) (-1.0000,1.6180,0.0000)};
  \addplot3[btwire, no marks, forget plot] coordinates {(1.0000,1.6180,0.0000) (1.6180,0.0000,1.0000)};
  \addplot3[btwire, no marks, forget plot] coordinates {(1.0000,1.6180,0.0000) (1.6180,0.0000,-1.0000)};
  \addplot3[btwire, no marks, forget plot] coordinates {(1.0000,-1.6180,0.0000) (-1.0000,-1.6180,0.0000)};
  \addplot3[btwire, no marks, forget plot] coordinates {(1.0000,-1.6180,0.0000) (1.6180,0.0000,1.0000)};
  \addplot3[btwire, no marks, forget plot] coordinates {(1.0000,-1.6180,0.0000) (1.6180,0.0000,-1.0000)};
  \addplot3[btwire, no marks, forget plot] coordinates {(-1.0000,1.6180,0.0000) (-1.6180,0.0000,1.0000)};
  \addplot3[btwire, no marks, forget plot] coordinates {(-1.0000,1.6180,0.0000) (-1.6180,0.0000,-1.0000)};
  \addplot3[btwire, no marks, forget plot] coordinates {(-1.0000,-1.6180,0.0000) (-1.6180,0.0000,1.0000)};
  \addplot3[btwire, no marks, forget plot] coordinates {(-1.0000,-1.6180,0.0000) (-1.6180,0.0000,-1.0000)};
  \addplot3[btwire, no marks, forget plot] coordinates {(1.6180,0.0000,1.0000) (1.6180,0.0000,-1.0000)};
  \addplot3[btwire, no marks, forget plot] coordinates {(-1.6180,0.0000,1.0000) (-1.6180,0.0000,-1.0000)};
  \addplot3[btbottleneck3d, no marks, forget plot] coordinates {(-1.6180,0.0000,-1.0000) (-0.0000,-0.4470,-0.7233)};
  \addplot3[btbottleneck3d, no marks, forget plot] coordinates {(-1.0000,-1.6180,0.0000) (-0.0000,-0.4470,-0.7233)};
  \addplot3[btedge3d, no marks, forget plot] coordinates {(0.0000,-1.0000,-1.6180) (-0.0000,-0.4470,-0.7233)};
  \addplot3[btbottleneck3d, no marks, forget plot] coordinates {(0.0000,1.0000,-1.6180) (-0.0000,-0.4470,-0.7233)};
  \addplot3[btbottleneck3d, no marks, forget plot] coordinates {(-1.6180,0.0000,1.0000) (-0.0000,0.4470,0.7232)};
  \addplot3[btbottleneck3d, no marks, forget plot] coordinates {(1.0000,-1.6180,0.0000) (-0.0000,-0.4470,-0.7233)};
  \addplot3[btedge3d, no marks, forget plot] coordinates {(-0.0000,-0.4470,-0.7233) (-0.0000,0.4470,0.7232)};
  \addplot3[btbottleneck3d, no marks, forget plot] coordinates {(-1.0000,1.6180,0.0000) (-0.0000,0.4470,0.7232)};
  \addplot3[btbottleneck3d, no marks, forget plot] coordinates {(1.6180,0.0000,-1.0000) (-0.0000,-0.4470,-0.7233)};
  \addplot3[btbottleneck3d, no marks, forget plot] coordinates {(0.0000,-1.0000,1.6180) (-0.0000,0.4470,0.7232)};
  \addplot3[btedge3d, no marks, forget plot] coordinates {(0.0000,1.0000,1.6180) (-0.0000,0.4470,0.7232)};
  \addplot3[btbottleneck3d, no marks, forget plot] coordinates {(1.0000,1.6180,0.0000) (-0.0000,0.4470,0.7232)};
  \addplot3[btbottleneck3d, no marks, forget plot] coordinates {(1.6180,0.0000,1.0000) (-0.0000,0.4470,0.7232)};
  \addplot3[only marks, mark = *, mark size = 2.4pt, mark options = {fill = btterminal, draw = btterminal}, forget plot] coordinates {(0.0000,1.0000,1.6180) (0.0000,1.0000,-1.6180) (0.0000,-1.0000,1.6180) (0.0000,-1.0000,-1.6180) (1.0000,1.6180,0.0000) (1.0000,-1.6180,0.0000) (-1.0000,1.6180,0.0000) (-1.0000,-1.6180,0.0000) (1.6180,0.0000,1.0000) (1.6180,0.0000,-1.0000) (-1.6180,0.0000,1.0000) (-1.6180,0.0000,-1.0000) };
  \addplot3[only marks, mark = triangle*, mark size = 3.4pt, mark options = {fill = btsteiner, draw = btsteiner}, forget plot] coordinates {(-0.0000,-0.4470,-0.7233) (-0.0000,0.4470,0.7232) };
\end{axis}
\end{tikzpicture}\\[2pt]
    {\small (e) \texttt{icosahedron}, $k=2$}
  \end{minipage}\hfill
  \begin{minipage}[t]{0.32\linewidth}\centering
\begin{tikzpicture}
\begin{axis}
  [
    steineraxis3d,
    width = 4.6000cm,
    height = 4.6000cm,
    xmin = -1.6828, xmax = 1.6828,
    ymin = -1.6828, ymax = 1.6828,
    zmin = -1.6828, zmax = 1.6828,
  ]
  \addplot3[btwire, no marks, forget plot] coordinates {(0.0000,1.0000,1.6180) (0.0000,-1.0000,1.6180)};
  \addplot3[btwire, no marks, forget plot] coordinates {(0.0000,1.0000,1.6180) (1.0000,1.6180,0.0000)};
  \addplot3[btwire, no marks, forget plot] coordinates {(0.0000,1.0000,1.6180) (-1.0000,1.6180,0.0000)};
  \addplot3[btwire, no marks, forget plot] coordinates {(0.0000,1.0000,1.6180) (1.6180,0.0000,1.0000)};
  \addplot3[btwire, no marks, forget plot] coordinates {(0.0000,1.0000,1.6180) (-1.6180,0.0000,1.0000)};
  \addplot3[btwire, no marks, forget plot] coordinates {(0.0000,1.0000,-1.6180) (0.0000,-1.0000,-1.6180)};
  \addplot3[btwire, no marks, forget plot] coordinates {(0.0000,1.0000,-1.6180) (1.0000,1.6180,0.0000)};
  \addplot3[btwire, no marks, forget plot] coordinates {(0.0000,1.0000,-1.6180) (-1.0000,1.6180,0.0000)};
  \addplot3[btwire, no marks, forget plot] coordinates {(0.0000,1.0000,-1.6180) (1.6180,0.0000,-1.0000)};
  \addplot3[btwire, no marks, forget plot] coordinates {(0.0000,1.0000,-1.6180) (-1.6180,0.0000,-1.0000)};
  \addplot3[btwire, no marks, forget plot] coordinates {(0.0000,-1.0000,1.6180) (1.0000,-1.6180,0.0000)};
  \addplot3[btwire, no marks, forget plot] coordinates {(0.0000,-1.0000,1.6180) (-1.0000,-1.6180,0.0000)};
  \addplot3[btwire, no marks, forget plot] coordinates {(0.0000,-1.0000,1.6180) (1.6180,0.0000,1.0000)};
  \addplot3[btwire, no marks, forget plot] coordinates {(0.0000,-1.0000,1.6180) (-1.6180,0.0000,1.0000)};
  \addplot3[btwire, no marks, forget plot] coordinates {(0.0000,-1.0000,-1.6180) (1.0000,-1.6180,0.0000)};
  \addplot3[btwire, no marks, forget plot] coordinates {(0.0000,-1.0000,-1.6180) (-1.0000,-1.6180,0.0000)};
  \addplot3[btwire, no marks, forget plot] coordinates {(0.0000,-1.0000,-1.6180) (1.6180,0.0000,-1.0000)};
  \addplot3[btwire, no marks, forget plot] coordinates {(0.0000,-1.0000,-1.6180) (-1.6180,0.0000,-1.0000)};
  \addplot3[btwire, no marks, forget plot] coordinates {(1.0000,1.6180,0.0000) (-1.0000,1.6180,0.0000)};
  \addplot3[btwire, no marks, forget plot] coordinates {(1.0000,1.6180,0.0000) (1.6180,0.0000,1.0000)};
  \addplot3[btwire, no marks, forget plot] coordinates {(1.0000,1.6180,0.0000) (1.6180,0.0000,-1.0000)};
  \addplot3[btwire, no marks, forget plot] coordinates {(1.0000,-1.6180,0.0000) (-1.0000,-1.6180,0.0000)};
  \addplot3[btwire, no marks, forget plot] coordinates {(1.0000,-1.6180,0.0000) (1.6180,0.0000,1.0000)};
  \addplot3[btwire, no marks, forget plot] coordinates {(1.0000,-1.6180,0.0000) (1.6180,0.0000,-1.0000)};
  \addplot3[btwire, no marks, forget plot] coordinates {(-1.0000,1.6180,0.0000) (-1.6180,0.0000,1.0000)};
  \addplot3[btwire, no marks, forget plot] coordinates {(-1.0000,1.6180,0.0000) (-1.6180,0.0000,-1.0000)};
  \addplot3[btwire, no marks, forget plot] coordinates {(-1.0000,-1.6180,0.0000) (-1.6180,0.0000,1.0000)};
  \addplot3[btwire, no marks, forget plot] coordinates {(-1.0000,-1.6180,0.0000) (-1.6180,0.0000,-1.0000)};
  \addplot3[btwire, no marks, forget plot] coordinates {(1.6180,0.0000,1.0000) (1.6180,0.0000,-1.0000)};
  \addplot3[btwire, no marks, forget plot] coordinates {(-1.6180,0.0000,1.0000) (-1.6180,0.0000,-1.0000)};
  \addplot3[btbottleneck3d, no marks, forget plot] coordinates {(-1.6180,0.0000,-1.0000) (-0.7596,-0.4694,0.2901)};
  \addplot3[btedge3d, no marks, forget plot] coordinates {(-1.0000,-1.6180,0.0000) (-0.7596,-0.4694,0.2901)};
  \addplot3[btedge3d, no marks, forget plot] coordinates {(-1.6180,0.0000,1.0000) (-0.7596,-0.4694,0.2901)};
  \addplot3[btbottleneck3d, no marks, forget plot] coordinates {(0.0000,-1.0000,-1.6180) (0.7564,-0.4956,-0.2782)};
  \addplot3[btbottleneck3d, no marks, forget plot] coordinates {(0.0000,-1.0000,1.6180) (-0.7596,-0.4694,0.2901)};
  \addplot3[btbottleneck3d, no marks, forget plot] coordinates {(-0.7596,-0.4694,0.2901) (0.7564,-0.4956,-0.2782)};
  \addplot3[btbottleneck3d, no marks, forget plot] coordinates {(-0.7596,-0.4694,0.2901) (-0.0281,0.9457,0.0000)};
  \addplot3[btbottleneck3d, no marks, forget plot] coordinates {(0.0000,1.0000,-1.6180) (-0.0281,0.9457,0.0000)};
  \addplot3[btedge3d, no marks, forget plot] coordinates {(1.0000,-1.6180,0.0000) (0.7564,-0.4956,-0.2782)};
  \addplot3[btedge3d, no marks, forget plot] coordinates {(-1.0000,1.6180,0.0000) (-0.0281,0.9457,0.0000)};
  \addplot3[btedge3d, no marks, forget plot] coordinates {(1.6180,0.0000,-1.0000) (0.7564,-0.4956,-0.2782)};
  \addplot3[btbottleneck3d, no marks, forget plot] coordinates {(0.0000,1.0000,1.6180) (-0.0281,0.9457,0.0000)};
  \addplot3[btbottleneck3d, no marks, forget plot] coordinates {(1.6180,0.0000,1.0000) (0.7564,-0.4956,-0.2782)};
  \addplot3[btedge3d, no marks, forget plot] coordinates {(1.0000,1.6180,0.0000) (-0.0281,0.9457,0.0000)};
  \addplot3[only marks, mark = *, mark size = 2.4pt, mark options = {fill = btterminal, draw = btterminal}, forget plot] coordinates {(0.0000,1.0000,1.6180) (0.0000,1.0000,-1.6180) (0.0000,-1.0000,1.6180) (0.0000,-1.0000,-1.6180) (1.0000,1.6180,0.0000) (1.0000,-1.6180,0.0000) (-1.0000,1.6180,0.0000) (-1.0000,-1.6180,0.0000) (1.6180,0.0000,1.0000) (1.6180,0.0000,-1.0000) (-1.6180,0.0000,1.0000) (-1.6180,0.0000,-1.0000) };
  \addplot3[only marks, mark = triangle*, mark size = 3.4pt, mark options = {fill = btsteiner, draw = btsteiner}, forget plot] coordinates {(-0.7596,-0.4694,0.2901) (-0.0281,0.9457,0.0000) (0.7564,-0.4956,-0.2782) };
\end{axis}
\end{tikzpicture}\\[2pt]
    {\small (f) \texttt{icosahedron}, $k=3$}
  \end{minipage}
  \caption{Optimal bottleneck Steiner trees on the octahedron (top row)
    and icosahedron (bottom row) as the Steiner budget grows. At $k=1$ the
    single Steiner point coincides with the centroid, the optimal
    bottleneck value equals the circumradius, and every spoke ties for the
    maximum, so all edges are highlighted; for $k\ge 2$ the tree
    reorganizes and the bottleneck value drops, from $1.0000$ to $0.8323$ and
    $0.7539$ on the octahedron and from $1.9021$ to $1.7013$ and $1.6192$
    on the icosahedron. Filled blue circles are terminals, orange triangles
    are Steiner points, and the thick red edges mark the
    maximum-length (bottleneck) edges; the thin light-gray lines trace the
    polytope's own edges for perspective, and a nearer tree edge is drawn
  over a farther one at every crossing.}
  \label{fig:named3d}
\end{figure}
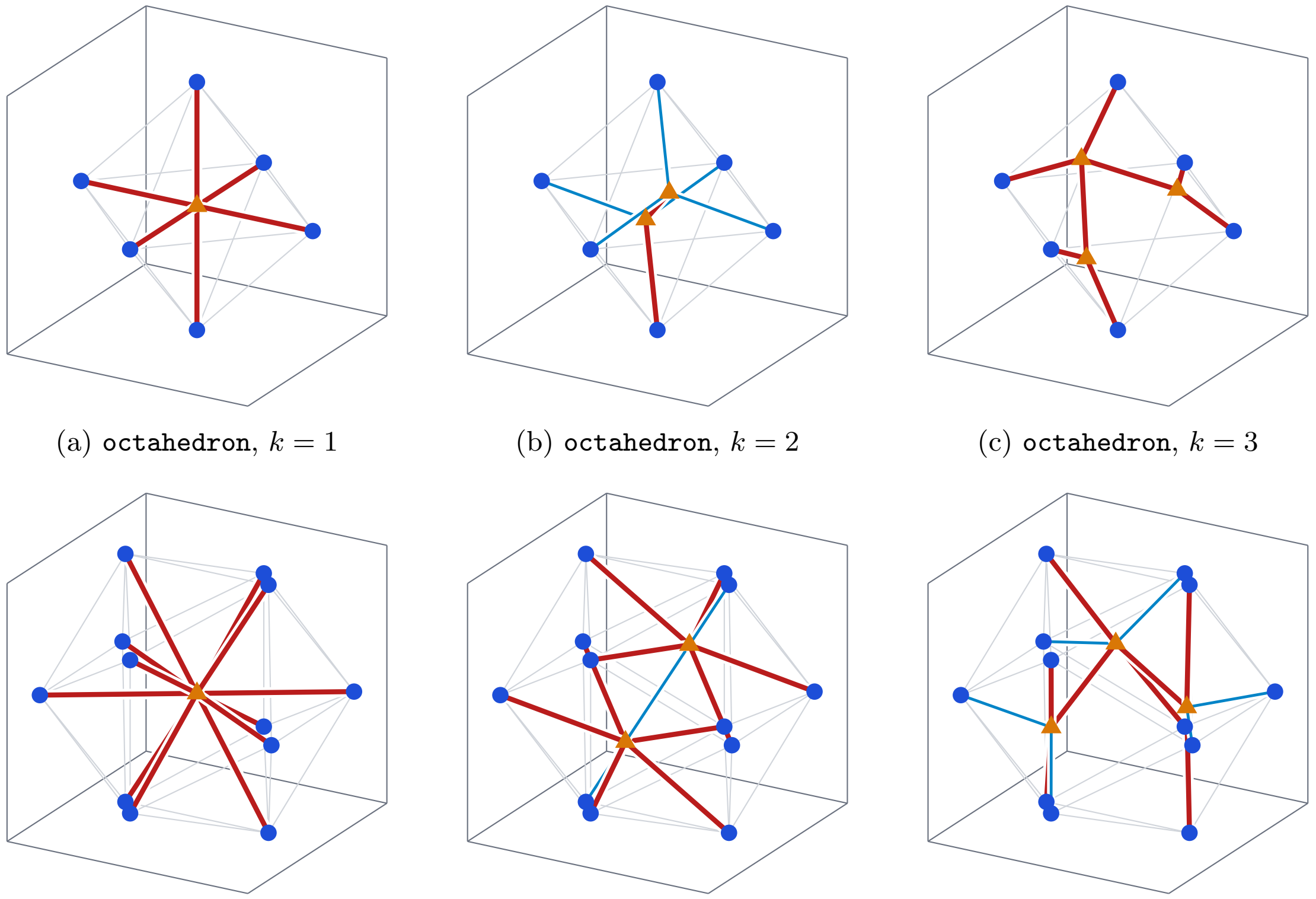

Having established {\rm \Idenone} as the faster and more robust formulation on the random benchmark, we evaluate both formulations on structured, reproducible inputs in higher dimensions than the random benchmark reaches: $30$ named instances built from regular polytopes, namely the five Platonic solids in $\mathbb{R}^3$, the simplices \texttt{nsimp-4/5} (in $\mathbb{R}^4$, $\mathbb{R}^5$), and the cross-polytopes \texttt{nocta-4/5/6} (in $\mathbb{R}^4$, $\mathbb{R}^5$, $\mathbb{R}^6$), following the higher-dimensional instance families of \citet*{smith1992}, each at $k\in\{1,2,3\}$. Every instance is solved with both formulations, with the candidate edge set $E^\star_H$ of~\eqref{eq:elimmst}, under a three-hour time limit. Table~\ref{tab:named} reports the results at the hardest budget $k=3$, where the two formulations separate. At $k\in\{1,2\}$, both formulations solve all $20$ instance-budget pairs, with maximum solve times of $3.1$\,s for {\rm \Idenone} and $11.5$\,s for {\rm \Identity}. Figure~\ref{fig:named3d} illustrates how the optimal tree reorganizes as the budget grows on two of the Platonic solids. For four of the five Platonic solids the single Steiner point at $k=1$ coincides with the centroid; for the dodecahedron the spanning tree on the vertices already attains the optimal value ($1.2361$, the polytope edge length, below the $\sqrt{3}$ of the centroid star).

\begin{table}[!ht]
  \centering
   
  \caption{Per-instance results of \Mone{} and \Mtwo{}, both with the MST restriction, on the $10$ named instances at the hardest budget $k=3$, under a three-hour (\num{10800}\,s) limit; columns follow the conic-Steiner convention~\citep*{pinto2023,maculan2025}. ``TL'' marks a run that reached the limit. For such a run, obj.\ is the better of its incumbent and $U_H$, the value its gap uses.}
  \label{tab:named}
  \small
  \setlength{\tabcolsep}{3.5pt}
  \renewcommand*{\arraystretch}{1.2}
  \begin{tabular}{l S[table-format=1.0] S[table-format=2.0] l S[table-format=1.4] S[table-format=5.1] S[table-format=5.0] S[table-format=2.2]}
    \toprule
    instance & {$d$} & {$p$} & {method} & {obj.} & {solve time (s)} & {nodes} & {gap (\%)}\\
    \midrule
\texttt{tetrahedron} & 3 & 4 & {\rm \Idenone} & 1.4142 & \best{0.0} & \best{5} & 0\\
 & {} & {} & {\rm \Identity} & 1.4142 & 0.1 & 20 & 0\\
\texttt{cube} & 3 & 8 & {\rm \Idenone} & 1.4142 & \best{1.1} & 1120 & 0\\
 & {} & {} & {\rm \Identity} & 1.4142 & 1.5 & \best{681} & 0\\
\texttt{octahedron} & 3 & 6 & {\rm \Idenone} & 0.7539 & \best{0.3} & \best{980} & 0\\
 & {} & {} & {\rm \Identity} & 0.7539 & 3.7 & 1082 & 0\\
\texttt{icosahedron} & 3 & 12 & {\rm \Idenone} & 1.6192 & \best{4.3} & \best{1284} & 0\\
 & {} & {} & {\rm \Identity} & 1.6192 & 214.2 & 19923 & 0\\
\texttt{dodecahedron} & 3 & 20 & {\rm \Idenone} & 1.2361 & \best{8.7} & \best{3146} & 0\\
 & {} & {} & {\rm \Identity} & 1.2361 & 108.8 & 16838 & 0\\
\texttt{nsimp-4} & 4 & 5 & {\rm \Idenone} & 0.7071 & \best{0.1} & 293 & 0\\
 & {} & {} & {\rm \Identity} & 0.7071 & 0.6 & \best{197} & 0\\
\texttt{nocta-4} & 4 & 8 & {\rm \Idenone} & 0.8165 & \best{4.8} & \best{2616} & 0\\
 & {} & {} & {\rm \Identity} & 0.8165 & 91.7 & 7812 & 0\\
\texttt{nsimp-5} & 5 & 6 & {\rm \Idenone} & 0.7229 & \best{3.4} & \best{2316} & 0\\
 & {} & {} & {\rm \Identity} & 0.7229 & 2466.6 & 11730 & 0\\
\texttt{nocta-5} & 5 & 10 & {\rm \Idenone} & \best{0.8660} & \best{300.5} & 14296 & \best{0}\\
 & {} & {} & {\rm \Identity} & 1.0981 & {TL} & 20816 & 23.19\\
\texttt{nocta-6} & 6 & 12 & {\rm \Idenone} & \best{0.8660} & \best{150.4} & 10487 & \best{0}\\
 & {} & {} & {\rm \Identity} & 1.4141 & {TL} & 23098 & 42.26\\
    \bottomrule
  \end{tabular}
\end{table}

{\rm \Idenone} solves all $30$ named instances to proven optimality; only two, \texttt{nocta-5} and \texttt{nocta-6} at $k=3$, require more than $10$\,s, at $300.5$\,s and $150.4$\,s, respectively, and every other {\rm \Idenone} solve finishes in at most $8.7$\,s. {\rm \Identity} solves $28$ of the $30$ instances but hits the three-hour limit on exactly those two hardest instances, terminating with optimality gaps of $23.19\%$ on \texttt{nocta-5} and $42.26\%$ on \texttt{nocta-6}; on both it reaches the limit with a feasible incumbent ($1.0981$ and $1.4141$, respectively) but cannot certify optimality. Consistent with the main comparison, {\rm \Identity} remains much slower at $k=3$ even on instances it solves: it takes $214.2$\,s on \texttt{icosahedron} against $4.3$\,s for {\rm \Idenone}, and $2466.6$\,s on \texttt{nsimp-5} against $3.4$\,s. On every solved instance the two formulations agree on the optimal value to at least four decimals.

\subsection{Root-relaxation strength}
\label{sec:rootexp}

Section~\ref{sec:mainexp} establishes that {\rm \Identity} is  slower than {\rm \Idenone} for $d\ge 3$ and explores larger search trees in most cells with $d\ge 3$, although its root relaxation is not weaker (Sections~\ref{sec:m2} and~\ref{sec:relaxation}). We examine a relaxation weakness shared by both formulations, the bottleneck epigraph of Challenge~\ref{c1}, quantifying it on the root relaxation of {\rm \Idenone}, whose static conic formulation gives a root bound ({\rm \Identity} accrues its bound through lazily separated cuts, so it has no comparable single static root), and measuring how much the two global lower bounds of Section~\ref{sec:diam} raise that bound.

We measure the gap between the heuristic upper bound $U_H$ and the dual bound $LB_{\mathrm{root}}$ obtained at the end of root processing,
\begin{equation}
  \operatorname{rootGap}_{U_H}\;:=\;\frac{U_H - LB_{\mathrm{root}}}{U_H}.
  \label{eq:rootgap}
\end{equation}
Smaller values indicate a stronger relaxation; normalization by $U_H$ makes gaps comparable across instance scales. { At the higher budgets, every instance is degenerate at the solver's precision: the median root gap of the bare relaxation is $1.000$ at the reported precision for $k\ge 3$, and the two global lower bounds of Section~\ref{sec:diam} bring this median down to between $0.48$ and $0.55$ without closing it (Table~\ref{tab:rootD}).}

\begin{table}[htbp]
  \centering
   
  \caption{Root-relaxation strength of \Mone{} with the MST restriction on the $1050$ instances of test set~C, stratified by
    Steiner budget $k$ ($210$ instances per row): mean and median of the
    root gap~\eqref{eq:rootgap}. ``bare'' uses neither global
    bound; ``$+$ diam'' adds the diameter bound~\eqref{eq:diamcut}; ``$+$ diam
    $+$ clust'' also adds the clustering bound~\eqref{eq:clustercut}, $e_k/2$ (the
    default). The better (smaller) median per $k$-block is bold.}
  \label{tab:rootD}
  \small
  \setlength{\tabcolsep}{4pt}
  \renewcommand*{\arraystretch}{1.3}
  \begin{tabular}{
      S[table-format=1.0] S[table-format=3.0]
    *3{S[table-format=1.3] S[table-format=1.3]}}
    \toprule
    & & \multicolumn{2}{c}{bare} & \multicolumn{2}{c}{$+$ diam}
    & \multicolumn{2}{c}{$+$ diam $+$ clust}\\
    \cmidrule(lr){3-4}\cmidrule(lr){5-6}\cmidrule(lr){7-8}
    {$k$} & {\#inst} & {mean} & {median} & {mean} & {median} & {mean} & {median}\\
    \midrule
1 & 210 & 0.854 & 0.860 & 0.821 & 0.818 & 0.358 & \best{0.418}\\
3 & 210 & 1.000 & 1.000 & 0.826 & 0.844 & 0.470 & \best{0.479}\\
5 & 210 & 1.000 & 1.000 & 0.809 & 0.827 & 0.563 & \best{0.507}\\
7 & 210 & 1.000 & 1.000 & 0.797 & 0.814 & 0.565 & \best{0.511}\\
9 & 210 & 1.000 & 1.000 & 0.789 & 0.803 & 0.572 & \best{0.547}\\
    \bottomrule
  \end{tabular}
\end{table}

{ Table~\ref{tab:rootD} shows that the bare relaxation gives an informative bound only at $k=1$, where its median gap is $0.860$. For $k\ge3$, every instance yields a root gap of $1.000$ at the reported precision (the root dual bounds lie below the $10^{-6}$ solver feasibility tolerance): the dual bound is numerically indistinguishable from the trivial lower bound $0$ that arises from setting $L=0$ in the cone constraints. Theorem~\ref{thm:degenerate} shows that this is not a numerical artifact: conditions~\eqref{eq:degen1}--\eqref{eq:degen2} hold on these instances, so the relaxation value is exactly $0$, and by Proposition~\ref{prop:rootLB} the root bound with the global lower bounds equals $\max\{\diam(P)/(p+k-1),\tfrac12 e_k\}$ on them.}

\paragraph{Interpretation} The min-max objective coupled to the big-$M$ epigraph is the source of the relaxation weakness for $k\geq 3$: as described in Challenge~\ref{c1}, at fractional $y_{ij}$ the big-$M$ term $M_{ij}(1-y_{ij})$ deactivates the cone constraint~\eqref{eq:soc}{ , and Theorem~\ref{thm:degenerate} characterizes exactly when this drives the relaxation value to $0$}. { On these instances the specialized constants and candidate edge sets
still admit the zero-valued construction of Theorem~\ref{thm:degenerate}.
The global lower bounds raise its value to $\underline L$.
A partial MIP start does not change the continuous feasible region.}

\paragraph{Effect of the global lower bounds} { Because the diameter bound~\eqref{eq:diamcut} is strictly positive, it removes the trivial 
$1.000$ gap entirely: the median falls from $1.000$ to between $0.803$ and $0.844$ for $k\ge3$,
while at $k=1$ it improves the median from $0.860$ to $0.818$. Adding the clustering
bound~\eqref{eq:clustercut} of Theorem~\ref{thm:clusterbound} is substantially stronger
still: it brings the median gap down to $0.418$ at $k=1$ and to between $0.479$ and $0.547$ for
$k\ge3$, roughly halving the root gap of the bare relaxation.} { The tighter root relaxation also
speeds up the full branch-and-bound search. On test set~B, with the MST restriction, {\rm \Idenone} solves $239$ of the $240$ instances
with and without the bounds, and enabling both bounds makes it faster by a median factor of $1.47$, and of $1.56$ on the
instances that take at least $10$\,s without them, and reduces the total solve time from \num{6330} to \num{4745}\,s, by $25\%$. On
the high-budget cells it lowers the median solve time from $28.4$ to $18.6$\,s at $d=3$, $k=5$ and
from $16.2$ to $9.6$\,s at $d=2$, $k=5$. Relative to the diameter bound alone, adding the
clustering bound gives a median speedup of $1.32$.}

{ On $46$ combinations of test set~C near the limit of {\rm\Idenone}, $460$ instances in
total, both bounds gain $23$ solved instances and lose $6$ against
no bounds (McNemar $p=0.002$), with a median speedup of $1.56$ on the instances both configurations
solve. Against the diameter bound alone, both bounds gain $13$ solved instances and lose $9$
($p=0.52$, not significant), with a median speedup of $1.31$ on the instances both configurations
solve.}

A bound that \emph{closes} the $k\ge 3$ gap would require separating the exponential family of cluster inequalities rather than the single worst-case clustering used here; we discuss this, and the reason a per-edge conic reformulation cannot help, as directions for future work in Section~\ref{sec:conclusion}.

\subsection{Symmetry breaking}
\label{sec:symmexp}

{ 
We isolate the practical effect of the
symmetry-breaking constraints~\eqref{eq:s1}: the fifth experiment solves the $240$ instances of
test set~B with and without them, within both formulations with the MST restriction
(Table~\ref{tab:symm}). With the ordering, {\rm \Idenone} solves $239$ instances against $219$
without it, and {\rm \Identity} $230$ against $193$; no instance is solved only without the
ordering. The effect concentrates at $k=5$, where {\rm \Idenone} solves $39$ against $30$
instances at $d=2$ and $40$ against $29$ at $d=3$, and {\rm \Identity} $30$ against $8$ at $d=3$.
The total solve time of {\rm \Idenone} falls from \num{31754} to \num{4745}\,s, and on the
instances both versions solve that take at least $10$\,s without the ordering, the ordering speeds {\rm \Idenone} up by
a median factor of $6.1$. Without~\eqref{eq:s1}, setting the symmetry detection of the solver to
its most aggressive level leaves the node count unchanged on all $219$ instances both settings
solve, so the solver does not find this symmetry on its own.

The ordering matters most where {\rm \Idenone} is near its limit. On the $420$ instances of test
set~C with $d\ge3$ and $k\ge5$, {\rm \Idenone} solves $74$ instances without~\eqref{eq:s1} against
$173$ with it, and none that it misses with the ordering. At $k=7$ and $k=9$ it solves none without
the ordering, against $56$ with it, and on the $74$ instances both versions solve it explores a median
of $4.2$ times as many nodes without it.
}

\begin{table}[htbp]
  \centering
   
  \caption{Symmetry breaking: the ordering~\eqref{eq:s1} (S1) against no symmetry breaking
    (none), within \Mone{} and \Mtwo{} with the MST restriction on test set~B, by ambient dimension and Steiner
    budget~$k$; columns as in Table~\ref{tab:main} without the percentage column ($40$
    instances per row). ``all'' aggregates $p\in\{10,15,20,25\}$ and both families;
    ``$<0.05$'' marks a sub-$0.05$\,s cell.}
  \label{tab:symm}
  \small
  \setlength{\tabcolsep}{3.5pt}
  \renewcommand*{\arraystretch}{1.3}
  \begin{tabular*}{\textwidth}{@{\extracolsep{\fill}}
      S[table-format=1.0] c S[table-format=1.0] S[table-format=2.0]
      *2{S[table-format=2.0] S[table-format=4.1] S[table-format=3.1]
    S[table-format=5.0] S[table-format=3.2]}}
    \toprule
    & & & & \multicolumn{5}{c}{S1} & \multicolumn{5}{c}{none}\\
    \cmidrule(lr){5-9}\cmidrule(lr){10-14}
    {$d$} & {$p$} & {$k$} & {\#inst} & {solv.} & {$t_{\med}$} & {$t_{\sgm}$} & {nodes} & {gap$^*$} & {solv.} & {$t_{\med}$} & {$t_{\sgm}$} & {nodes} & {gap$^*$}\\
    \midrule
\multicolumn{14}{@{}l}{\itshape \Mone{}}\\
2 & {all} & 1 & 40 & 40 & {$<0.05$} & {$<0.05$} & 1 & 0 & 40 & {$<0.05$} & {$<0.05$} & 1 & 0\\
2 & {all} & 3 & 40 & 40 & \best{2.0} & \best{2.8} & \best{1094} & 0 & 40 & 2.5 & 3.2 & 1462 & 0\\
2 & {all} & 5 & 40 & \best{39} & \best{9.6} & \best{17.4} & 3715 & \best{1.14} & 30 & 188.3 & 171.7 & 48310 & 48.63\\
3 & {all} & 1 & 40 & 40 & {$<0.05$} & 0.1 & 1 & 0 & 40 & {$<0.05$} & 0.1 & 1 & 0\\
3 & {all} & 3 & 40 & 40 & \best{3.6} & \best{4.6} & \best{1534} & 0 & 40 & 5.4 & 7.1 & 1836 & 0\\
3 & {all} & 5 & 40 & \best{40} & \best{18.6} & \best{27.3} & 8445 & \best{0} & 29 & 188.5 & 203.3 & 31957 & 55.12\\
\midrule
\multicolumn{14}{@{}l}{\itshape \Mtwo{}}\\
2 & {all} & 1 & 40 & 40 & 0.1 & 0.1 & 37 & 0 & 40 & {\bestmath{<0.05}} & 0.1 & 37 & 0\\
2 & {all} & 3 & 40 & 40 & \best{3.2} & \best{3.4} & \best{786} & 0 & 40 & 4.3 & 5.6 & 987 & 0\\
2 & {all} & 5 & 40 & \best{40} & \best{20.6} & \best{32.2} & 4040 & \best{0} & 25 & 195.8 & 274.8 & 12239 & 68.72\\
3 & {all} & 1 & 40 & 40 & 0.1 & \best{0.1} & 38 & 0 & 40 & 0.1 & 0.2 & 38 & 0\\
3 & {all} & 3 & 40 & 40 & \best{18.1} & \best{22.0} & \best{3078} & 0 & 40 & 51.5 & 51.7 & 6072 & 0\\
3 & {all} & 5 & 40 & \best{30} & \best{491.3} & \best{492.7} & 18808 & 77.24 & 8 & 1000.0 & 894.2 & 18186 & \best{75.45}\\
    \bottomrule
  \end{tabular*}
\end{table}

\subsection{\texorpdfstring{Outer approximation}{Outer approximation}}
\label{sec:oaexp}

{ 
The cone constraints~\eqref{eq:soc} of {\rm \Idenone} can be replaced by linear supporting
inequalities generated on demand~\citep*{duran1986,quesada1992}. Gurobi can add these inequalities
inside its own branch-and-bound tree (\texttt{MIQCPMethod}${}=1$); a multi-tree loop instead solves
a mixed-integer linear master problem, adds the inequalities that its solution violates, and
solves the master problem again. We compared both with {\rm \Idenone} under Gurobi's default and
under its other method, a continuous conic relaxation at every node (\texttt{MIQCPMethod}${}=0$),
on the $240$ instances of test set~B, all with the MST restriction $E^\star_H$
of~\eqref{eq:elimmst}. The default solves $239$ instances, as does \texttt{MIQCPMethod}${}=1$;
\texttt{MIQCPMethod}${}=0$ solves $217$ and the multi-tree loop $215$, and neither of the last two
solves an instance that the default misses. The difference concentrates at $d=3$, $k=5$, where the
multi-tree loop solves $19$ of the $40$ instances against $40$ for the default and leaves a
largest gap of $50.25\%$. All variants agree with the default on the optimal value within a
relative tolerance of $10^{-4}$ on every instance both solve.

On the instances that each variant and the default both solve, \texttt{MIQCPMethod}${}=1$ matches the default in
median solve time and, on the $133$ instances that take at least $1$\,s, explores a median
of $1.24$ times as many nodes; \texttt{MIQCPMethod}${}=0$ is slower by a median factor of $2.0$ and the
multi-tree loop by $5.9$.  Outer approximation inside the solver's search tree thus
solves as many instances as the default, while the multi-tree form is not competitive, and we do
not pursue it further.
}

\FloatBarrier


\section{Conclusion}
\label{sec:conclusion}

We proposed two exact mixed-integer second-order cone formulations
for the Euclidean bottleneck Steiner tree problem (EBSTP) in $\mathbb{R}^d$. We also developed the structural strengthening techniques that make these formulations computationally effective and evaluated their performance { on three test sets of randomly generated instances and on a set of instances built from regular polytopes}. The computational study includes the $1050$ instances of the main comparison, with up to $50$ terminals and $9$ Steiner points, $240$ in the experiments that isolate individual components, $40$ planar validation instances, and $30$ structured named instances. The main findings are summarized below.
\begin{enumerate}[leftmargin=*,itemsep=0.3em]
  \item {\rm \Idenone} and {\rm \Identity} reproduce the optimal value found by the combinatorial planar algorithm of \citet*{bae2009exact} on all $40$ planar test instances of test set~A; all per-instance discrepancies fall within the configured MIP optimality tolerance of~$10^{-4}$.
  \item Across the $1050$ instances of test set~C in $d\in\{2,3,4\}$, both with the MST restriction, {\rm \Idenone} is faster than {\rm \Identity} in median solve time in every dimension, { by factors of $1.6$, $17.6$ and $5.9$ at $d=2,3,4$, and is the faster method on $585$ against $55$ instances; it also solves $168$ instances that {\rm\Identity} does not, against $1$, a difference that a paired McNemar test finds significant ($p<10^{-40}$). The difference in solved instances grows with the dimension and the Steiner budget: {\rm\Idenone} solves $27$, $63$ and $77$ more instances than {\rm\Identity} at $d=2,3,4$, and at $(d,k)=(4,5)$, {\rm\Identity} solves $6$ of the $70$ instances and {\rm\Idenone} $55$, and with $k=9$, {\rm\Idenone} solves $33$ instances against $9$.}
  \item { Every test instance with $k\ge3$ is degenerate at the solver's precision, and by Theorem~\ref{thm:degenerate} no instance with $k=1$ is; this is the principal Challenge~\ref{c1} of the bottleneck-epigraph approach. The two proved global lower bounds, the diameter bound and the clustering bound, are the only strengthening techniques that act directly on the relaxation of $L$; together they replace the trivial root bound by a positive one and roughly halve the root gap (to a median of $0.42$ at $k=1$ and between $0.48$ and $0.55$ for $k$ from $3$ to $9$), but do not close it.}
  \item The coordinate ordering~S1 yields $20$ additional solved instances for {\rm\Idenone} and $37$ for {\rm\Identity} on test set~B, and $99$ for {\rm\Idenone} on the $420$ instances of test set~C with $d\ge3$ and $k\ge5$; the solver does not detect the symmetry it breaks.
  Outer approximation inside the solver's own search tree solves as many instances as
  the default, at the same median solve time, and a multi-tree loop is not competitive
  (Section~\ref{sec:oaexp}). On the structured named instances, {\rm \Idenone} solves all $30$ instances within minutes whereas {\rm \Identity} times out on the two hardest.
\end{enumerate}
The empirical evidence identifies {\rm \Idenone} with the structural strengthening techniques developed in this paper as the recommended exact baseline for the EBSTP in $\mathbb{R}^d$.

The main open problem is the bottleneck epigraph relaxation (Challenge~\ref{c1}): the two global lower bounds of Section~\ref{sec:diam} only halve the residual root gap, so closing it would require a different mechanism. { A per-edge perspective or disjunctive reformulation does not directly address this weakness}: every edge shares the same $x_i$, $x_j$, and $L$, so no per-edge switch isolates $L$ in a perspective reformulation~\citep*{gunluk2010perspective}, and the two-term disjunctive conic cut~\citep*{kilinckarzan2015twoterm} does not apply directly to the shifted-apex cones that the disjunction on $y_{ij}$ produces. A possible approach is to separate the exponential family of cluster inequalities generalizing~\eqref{eq:clustercut}, a separation problem we expect to be no easier than the underlying problem, which is \textsf{NP}-hard~\citep*{wang2002,wang2002bsk}.

The remaining challenges suggest two directions for future research.
\begin{itemize}[leftmargin=*,itemsep=0.3em]
  \item Designing a symmetry-breaking scheme that exploits the geometry of the EBSTP (for instance, ordering Steiner points along a principal axis of the terminal set) rather than the coordinate ordering~\eqref{eq:s1}, with connections to the isomorphism-elimination techniques surveyed by \citet*{margot2010} (Challenge~\ref{c2}). A lexicographic order on two coordinates would have to distinguish a strictly smaller first coordinate from an equal one, and the strict branch is a condition that a mixed-integer formulation cannot model exactly without a tolerance.
  \item { Developing a tuned MMX-style bottleneck MINLP, solved by a global optimizer, to enable a direct comparison that the existing additive conic formulations~\citep*{pinto2023,maculan2025} cannot provide.
    Deriving, validating, and tuning such a formulation
    requires a substantial separate study and lies beyond the scope of the present paper. The conic and MINLP formulations of \citet*{pinto2023} and the review of \citet*{maculan2025} target the \emph{additive} Euclidean Steiner objective; because they minimize total edge length rather than the longest edge, a head-to-head comparison on a shared instance and objective is not meaningful without re-deriving those formulations for the bottleneck criterion.} A related direction is to connect the empirical dependence of {\rm \Idenone} on the Steiner budget $k$ to the fixed-parameter-tractable framework of \citet{bandyapadhyay2024} (Challenge~\ref{c4}).
\end{itemize}

\medskip

\noindent {\bf Acknowledgments.} { We are grateful to an anonymous reviewer for observing
that the exchange argument of \citet*{bae2010exact} extends to $\mathbb{R}^d$, which gives
Proposition~\ref{prop:mst} and the restriction of the terminal--terminal candidate edges to a
minimum spanning tree in every dimension.} M. Fampa was supported in part by CNPq grant 307167/2022-4.

{ 
\section*{Declaration of generative AI and AI-assisted technologies in the manuscript preparation process}

During the preparation of this work the authors used Claude (Anthropic) in order to
improve the language of the manuscript and, in the revised version, to assist in drafting
revised passages and in checking bibliographic references. The formulations, proofs,
implementation, experimental design and analysis are the authors' own. After using this
tool, the authors reviewed and edited the content as needed and take full responsibility
for the content of the publication.
}

\end{document}